\documentclass[leqno]{article}

\usepackage[cp1250]{inputenc}
\usepackage{authblk}
\usepackage[english]{babel}

\usepackage{a4wide}
 
\usepackage{amsfonts}
\usepackage{graphicx}
\usepackage{epstopdf}
\usepackage{bm}
\usepackage{subfig}
\usepackage{amsthm} 
\usepackage{amsmath}
\usepackage{color}
\usepackage{picture}
\usepackage{tikz}
\usepackage{adjustbox}
\usepackage{url}
\newcommand{\dx}{\ d \boldsymbol{x}}

\newcommand{\ds}{\ d s}

\newtheorem{testcase}{Test case}
\title{\textbf{Numerical assessment of two-level domain decomposition preconditioners for incompressible Stokes and elasticity equations}} \author[1]{Gabriel R. Barrenechea}
\author[1]{Micha\l{} Bosy \thanks{Corresponding author: E-mail: michal.bosy@strath.ac.uk}}
\author[1]{Victorita Dolean}
\affil[1]{\small \textit{Department of Mathematics and Statistics, University of Strathclyde, 26 Richmond Street, G1 1XH Glasgow, United Kingdom}}

\date{\today}

\providecommand{\keywords}[1]{\textbf{Key words.} #1}

\begin{document}
\maketitle

\ 

\begin{abstract}
Solving the linear elasticity and Stokes equations by an optimal domain decomposition method derived algebraically involves the use of non standard interface conditions. The one-level domain decomposition preconditioners are based on the solution of local problems. This has the undesired consequence that the results are not scalable, it means that the number of iterations needed to reach convergence increases with the number of subdomains. This is the reason why in this work we introduce, and test numerically, two-level preconditioners. Such preconditioners use a coarse space in their construction. We consider the nearly incompressible elasticity problems and Stokes equations, and discretise them by using two finite element methods, namely, the hybrid discontinuous Galerkin and Taylor-Hood discretisations.
\end{abstract}

 % \\
% 
%  \textbf{Mathematics Subject Classification (2000):} 65F10, 65N22, 65N30, 65N55}

% \thispagestyle{empty}
% 
% \newpage
% 
% \tableofcontents
% 
% \newpage

\section{Introduction}
\label{sec:introduction}

\begin{figure}[!ht]
\centering
    \subfloat[Taylor-Hood]{%
      \includegraphics[width=0.45\textwidth]{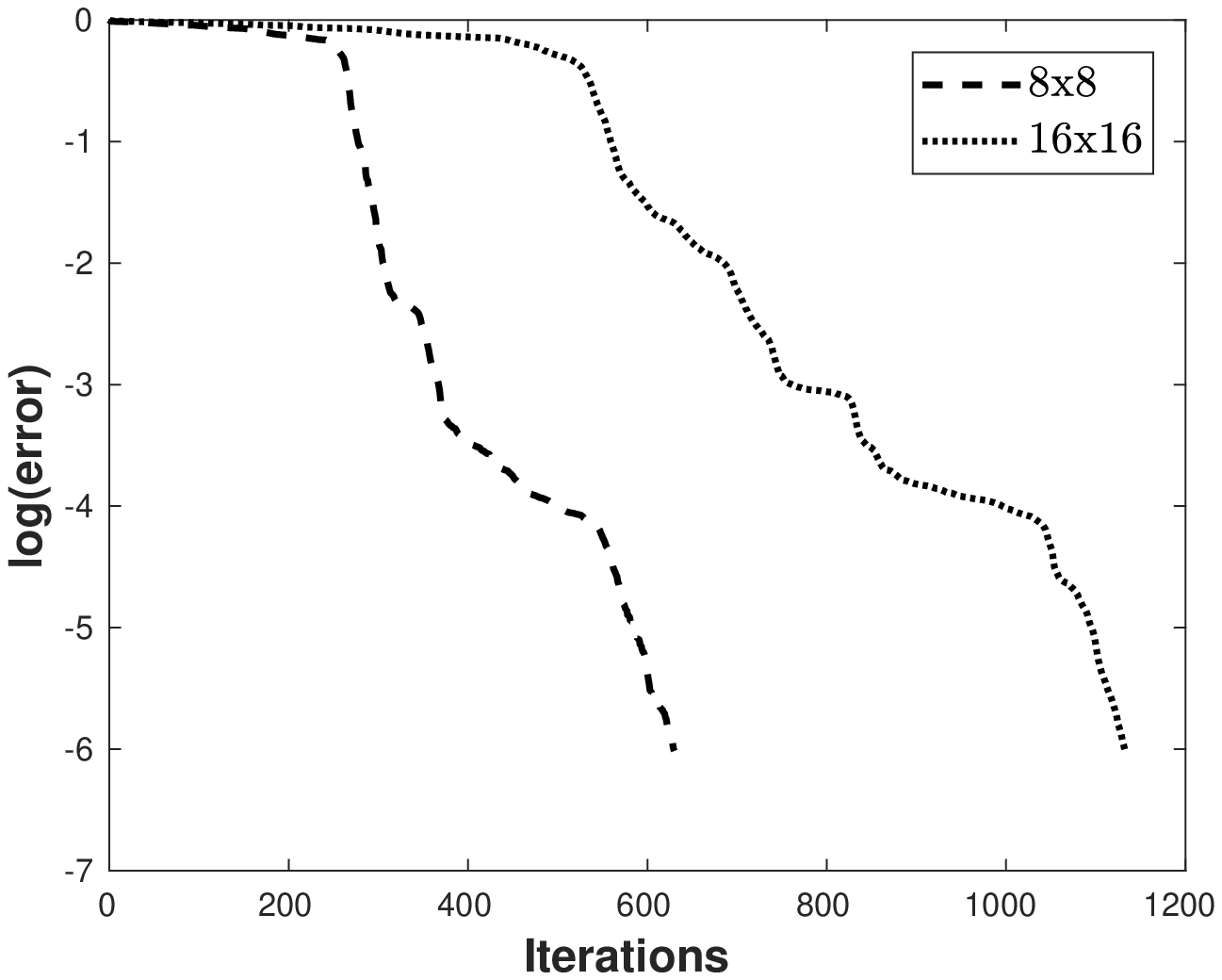}
    }
    \subfloat[hdG]{%
      \includegraphics[width=0.45\textwidth]{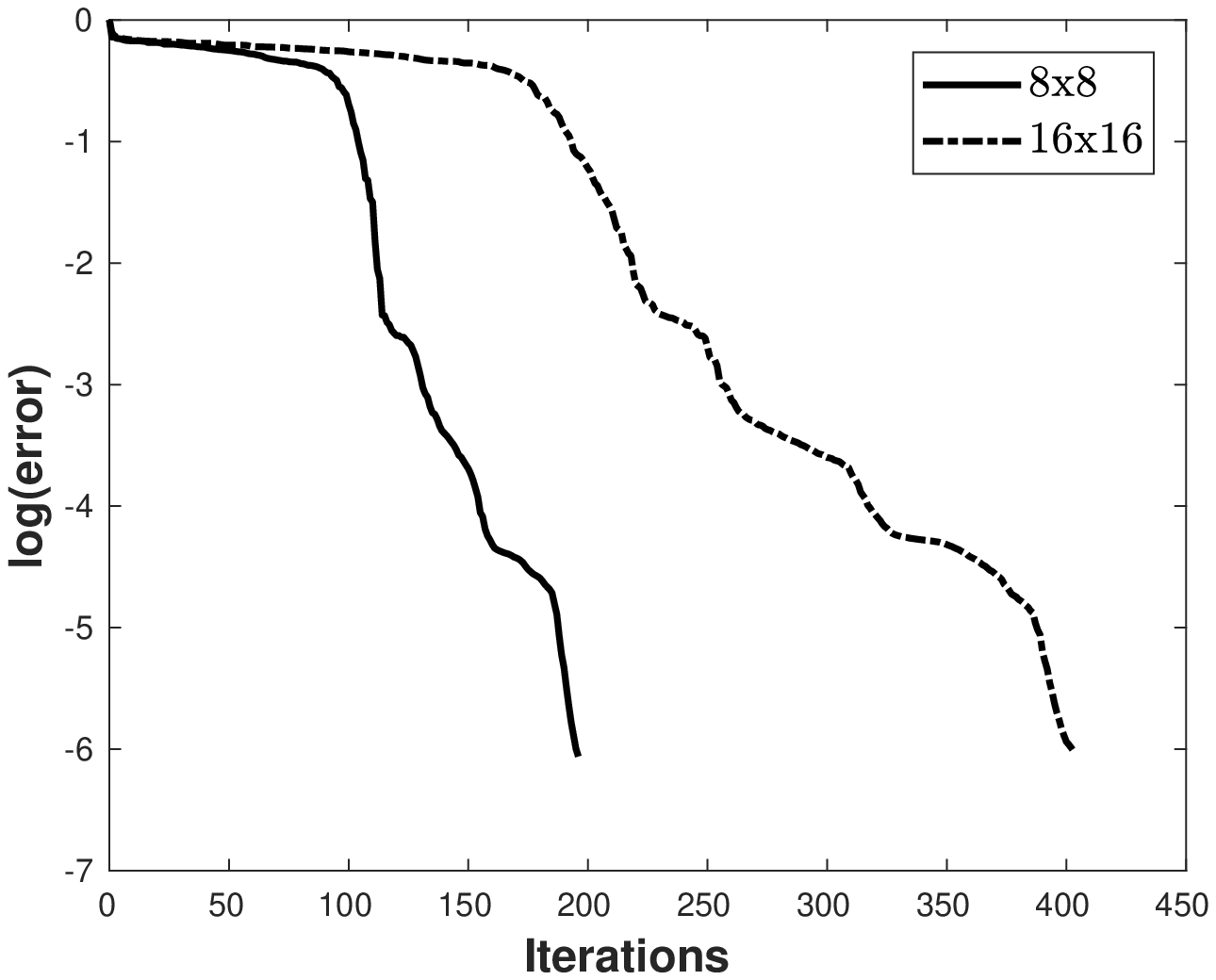}
    }
 	\caption{Increase of the size of the plateau region for increasing number of subdomains.}
  \label{fig:plateau}
  \end{figure}

In~\cite{barrenechea2016hybrid} the one-level domain decomposition methods for Stokes equations were introduced in conjunction with the non standard interface conditions. Although it can be observed there the lack of scalability with respect to the number of subdomains. It means that by splitting the problem in a larger number of subdomains leads to the increase of size of the plateau region in the convergence of an iterative method (see Figure~\ref{fig:plateau}) when using the one-level domain decomposition methods. This is caused by the lack of global information, as subdomains can only communicate with their neighbours. Hence, when the number of subdomains increases in one direction, the length of the plateau also increases. Even in cases when the local problems are of the same size, the iteration count grows with the increase of the number of subdomain. This can be also observed in all experiments in this manuscript in case of one-level methods.

The remedy for this is the use of a second level in the preconditioner or a coarse space correction that adds the necessary global information. Two-level algorithms have been analysed for several classes of problems in~\cite{Toselli:2005:DDM}. %The analysis is based on the idea of {\it stable splitting} which leads to convergence estimates independent of the parameters of the problems when the number of degrees of freedom per domain is kept constant. 
The key point of these kind of methods is to choose the appropriate coarse space. The classical coarse space introduced by Nicolaides in~\cite{MR881370} is defined by vectors that support is in each subdomain. Hence the coarse space has the size equal to a number of subdomains. A coarse space construction, named spectral coarse space, was motivated by the complexity of the problems that classical coarse space performance was not satisfying. This construction allows to enrich a bigger size of the coarse space, but can be also reduced to the classical one.
% In many applications we have to deal with strongly heterogeneous problems. For this reason, building two-level methods which are also robust with respect to the physical coefficients is very important. 
This idea was introduced for the first time in~\cite{brezina1999iterative} in the case of multigrid methods. It relies on solving local generalised eigenvalue problems allowing to choose suitable vectors for the coarse space.

For overlapping domain decomposition preconditioners, a similar idea was introduced in the case of Darcy equations in~\cite{MR2718268,MR2728702}.
%} by solving local generalised eigenvalue problems in the subdomains. Their solutions were used to identify vanishing eigenvalues and adding the associated eigenvectors to the coarse space. Unfortunately, since the right hand side of each problem involves a local mass matrix, this can lead to a large coarse space. That is why, the same authors redefined this approach by using a multiscale partition of unity in~\cite{.
The authors of~\cite{MR2738920} consider also the heterogeneous Darcy equation and presented a different generalised eigenvalue problem based on local Dirichlet-to-Neumann maps. %This approach comes down to solve problems only on the interfaces of the subdomains. 
The method has been analysed in~\cite{MR3033238} and proved to be very robust in the case of small overlaps. The same idea was extended numerically to the heterogeneous Helmholtz problem in~\cite{MR3209915}. The authors of~\cite{MR3432853} apply the coarse space associated with low-frequency eigenfunctions of the subdomain Dirichlet-to-Neumann maps for the generalisation of the optimised Schwarz methods, named 2-Lagrange multiplier methods.%The optimisation is associated in usage the local Robin boundary value problems instead of the local Dirichlet or Neumann ones as it is in the previous cases.

The first attempt to extend this spectral approach to general symmetric positive definite problems was made in~\cite{MR2916377} as an extension of~\cite{MR2718268,MR2728702}. Since some of the assumptions of the previous framework are hard to fulfil, authors of~\cite{MR3175183} proposed slightly different approach for symmetric positive definite problems. Their idea of constructing partition of unity operator associated with degrees of freedom allows to work with various finite element spaces. An overview of different kinds of two-level methods can be found in~\cite[Chapters 5 and 7]{MR3450068}.

Despite the fact that all these approaches provide satisfying results, there is no universal treatment to build efficient coarse spaces in the case of non definite problems such as Stokes equations. The spectral coarse spaces that we use in this work are inspired by those proposed in~\cite{haferssas:hal-01278347}. The authors introduced and tested numerically symmetrised two-level preconditioners for overlapping algorithms which use Robin interface conditions between the subdomains (see~\eqref{eq:coarse_elasticity_Robin_BC} for details). They have applied these preconditioners to solve saddle point problems such as nearly incompressible elasticity and Stokes discretised by Taylor-Hood finite elements. In our case, we use non standard interface conditions. Therefore the use of spectral coarse spaces could lead to an important gain.

In this work, we test this improvement in case of nearly incompressible elasticity and Stokes equations that are discussed in Section~\ref{sec:problem}. As the discretisations we use the Taylor-Hood~\cite[Chapter II, Section 4.2]{MR851383} and hybrid discontinuous Galerkin method~\cite{MR2485455,MR2485446} presented in Section~\ref{sec:methods}. In Section~\ref{sec:dd} we introduce the two-level domain decomposition preconditioners. Sections~\ref{sec:2D_numerics} and~\ref{sec:3D_numerics} present the two and three dimensional numerical experiments, respectively. Finally, a summary is outlined in Section~\ref{sec:conclusions}.

%%%%%%%%%%%%%%%%%%%%%%%%%%%%%%%%%%%%%%%%%%%%%%%%%%%%%%%%%%%%%%%%%%%%%%%%%%%%%%%%%%%%%%%%%%%%%%%%%%%%%%%%%%%%
% Problem
%%%%%%%%%%%%%%%%%%%%%%%%%%%%%%%%%%%%%%%%%%%%%%%%%%%%%%%%%%%%%%%%%%%%%%%%%%%%%%%%%%%%%%%%%%%%%%%%%%%%%%%%%%%%

\section{The differential equations}
\label{sec:problem}

Let $\Omega$ be an open polygon in $\mathbb{R}^2$ or an open Lipschitz polyhedron in $\mathbb{R}^3$, with Lipschitz boundary $\Gamma := \partial \Omega$. We use $d = 2,3$ to denote the dimension of the space. We use bold for tensor or vector variables. 
%%%%%%%%%%%%%%%%%%%%%%%%%%%%%%%%%%%%%%%%%%%%%%%%%%
%% Definition of normal and tangential part
In addition we denote normal and tangential components as follows ${u_n} := \boldsymbol{u} \cdot \boldsymbol{n}$, and $u_t := \boldsymbol{u} - u_n  \boldsymbol{n}$, where $\boldsymbol{n}$ is the outward unit normal vector to the
boundary $\Gamma$ . 

For $D \subset \Omega$, we use the standard $L^2(D)$ space
 %with following norm
% \begin{eqnarray*}
% \|f\|_D^2 := \int_D f^2 \dx & \mbox{for all } f \in L^2(D).
% \end{eqnarray*}
and $C^0(\bar{D})$ denotes the set of all continuous functions on the closure of a set $D$.
%%%%%%%%%%%%%%%%%%%%%%%%%%%%%%%%%%%%%%%%%%%%%%%%%%
%% Spaces
Let us define following Sobolev spaces
\begin{align*}
H^m(D) & :=  \left\{v \in L^2(D): \ \forall \ |\boldsymbol{\alpha}| \leq m \ \partial^{\boldsymbol{\alpha}} {v} \in L^2(D)\right\} \mbox{ for } m \in \mathbb{N}, \\
H^{\frac{1}{2}}(\partial D) & :=  \left\{\tilde{v} \in L^2(\partial D): \ \exists v \in H^1(D) \quad \tilde{v} = \mbox{tr}(v) \right\}, \\
H\left(div, D\right) & := \left\{\boldsymbol{v} \in [L^2(D)]^d: \ \nabla \cdot \boldsymbol{v} \in L^2(D)\right\},
\end{align*} 
where, for $\boldsymbol{\alpha} = (\alpha_1,..., \alpha_d) \in \mathbb{N}^d$ and
$|\boldsymbol{\alpha}| = \sum_{i=1}^{d}\alpha_i$ we denote $\partial^{\boldsymbol{\alpha}} = \frac{\partial^{|\boldsymbol{\alpha}|}}{\partial x_1^{\alpha_1} ... \partial x_d^{\alpha_d}}$, and $\mbox{tr}: H^1(\Omega) \rightarrow H^{\frac{1}{2}}(\partial \Omega)$ is the trace operator.
In addition, we %will use following standard semi-norm and norm for the Sobolev space $H^m(D)$ for $m \in \mathbb{N}$ 
% \begin{align*}
% |f|_{H^m(D)}^2 := \sum_{|\boldsymbol{\alpha}| = m} \|\partial^{\boldsymbol{\alpha}} f\|_D^2 && \|f\|_{H^m(D)}^2 := \sum_{k = 0}^m |f|_{H^k(D)}^2 && \mbox{for all } f \in H^m(D).
% \end{align*}
% We 
use the following notation of the space including boundary and average conditions
\begin{align*}
L^2_0(D) &:= \left\{v  \in L^2(D): \ \int_D v \dx = 0\right\}, \\
H^1_{\tilde{\Gamma}}(D) &:= \left\{v  \in H^1(D): \ v = 0 \mbox{ on } \tilde{\Gamma}\right\},
\end{align*} 
where $\tilde{\Gamma} \subset \partial D$. If $\tilde{\Gamma} = \partial D$, then $H^1_{\tilde{\Gamma}}(D)$ is denoted $H^1_0(D)$.

Now we present the two differential problems considered in this work.

\subsection{Stokes equation}
\label{sec:stokes}
%%%%%%%%%%%%%%%%%%%%%%%%%%%%%%%%%%%%%%%%%%%%%%%%%%%%%%%%%%%%%%%%%%%%%%%%%%%%%%%%%%%%%%%%%%%%%%%%%%%%%%%%%%%%
% Problem
Let us start with $d$-dimensional, $d=2,3$, Stokes problem 
\begin{equation}
\label{eq:stokes}
\left\{
\begin{array}{rclclr}
 -\nu \Delta \boldsymbol{u} & \ + & \nabla p & = & \boldsymbol{f} & \mbox{in } \Omega \\
 & & \nabla \cdot \boldsymbol{u} & = & 0 & \mbox{in } \Omega
 \end{array}
\right. ,
\end{equation} 
where $\boldsymbol{u}: \bar{\Omega} \rightarrow \mathbb{R}^d$ is the velocity field, $p:\bar{\Omega} \rightarrow \mathbb{R}$ the pressure, $\nu >
0$ the viscosity which is considered to be constant and $\boldsymbol{f} \in [L^2(\Omega)]^d$ is a given function. We define the stress tensor $\boldsymbol{\sigma} := \nu \nabla \boldsymbol{u} - p \boldsymbol{I}$ and the flux as $\boldsymbol{\sigma_n} := \boldsymbol{\sigma} \ \boldsymbol{n}$. 
%%%%%%%%%%%%%%%%%%%%%%%%%%%%%%%%%%%%%%%%%%%%%%%%%%
%% boundary conditions
For $\boldsymbol{u_D} \in [H^{\frac{1}{2}}(\Gamma)]^d$ and ${g} \in L^2(\Gamma)$ we consider three types of boundary conditions:
\begin{itemize}
  \item Dirichlet (non-slip) 
\begin{equation}
 \label{eq:Dirichlet}
 \boldsymbol{u} = \boldsymbol{u_D}  \mbox{ on } \Gamma;
\end{equation}
  \item tangential-velocity and normal-flux (TVNF) 
  \begin{equation}
  \left\{
  \begin{array}{rcll}
  \label{eq:TVNF}
   u_t & = & 0 & \mbox{ on } \Gamma \\
   {\sigma_{nn}} & = & {g} & \mbox{ on } \Gamma
  \end{array}
  \right. ;
  \end{equation}
  \item normal-velocity and tangential-flux (NVTF) 
\begin{equation}
 \label{eq:NVTF}
 \left\{
 \begin{array}{rcll}
  {u_n} & = & 0 & \mbox{ on } \Gamma \\
{\sigma_{nt}} & = & {g} & \mbox{ on } \Gamma
 \end{array}
 \right. .
\end{equation}
\end{itemize}
The third type of boundary condition has already been considered for the Stokes problem in~\cite{de2014simple}.
% which together with \eqref{eq:stokes} define two boundary value
% problems. 

\subsection{Nearly incompressible elasticity equation}
\label{sec:elasticity}

From a mathematical point of view, the nearly incompressible elasticity problem is very similar to the Stokes equations. The difference is that instead of considering  the gradient $\nabla \boldsymbol{v}$, the symmetric gradient $\boldsymbol{\varepsilon}(\boldsymbol{v}) := \frac{1}{2} \left(\nabla \boldsymbol{v}+ \nabla^T \boldsymbol{v}\right)$ is used. We want to solve the following $d$-dimensional, $d=2,3$, problem
\begin{equation}
\label{eq:elasticity}
\left\{
\begin{array}{rclclr}
 -2 \mu \nabla \cdot \boldsymbol{\varepsilon}(\boldsymbol{u}) & + & \nabla p & = & \boldsymbol{f} & \mbox{in } \Omega \\
 & - &  \nabla \cdot \boldsymbol{u} & = & \frac{1}{\lambda} p & \mbox{in } \Omega
 \end{array}
\right.
\end{equation} 
where $\boldsymbol{u}: \bar{\Omega} \rightarrow \mathbb{R}^d$ is the displacement field, $p:\Omega \rightarrow \mathbb{R}$ the pressure, $\boldsymbol{f} \in [L^2(\Omega)]^d$ is a given function, $\lambda$ and $\mu$  are the Lam\'{e} coefficients defined by
\begin{align*} 
	\lambda = \frac{E \nu}{(1+ \nu)(1 - 2 \nu)}, && \mu = \frac{E}{2(1+\nu)},
\end{align*}
where $E$ is the Young modulus and $\nu$ the Poisson ratio. We define the stress tensor as $\boldsymbol{\sigma}^{sym} := 2 \mu \boldsymbol{\varepsilon}(\boldsymbol{u}) - p \boldsymbol{I}$ and its normal component as $\boldsymbol{\sigma}^{sym}_{\boldsymbol{n}}:= \boldsymbol{\sigma}^{sym} \boldsymbol{n}$.
%%%%%%%%%%%%%%%%%%%%%%%%%%%%%%%%%%%%%%%%%%%%%%%%%%
%% boundary conditions
For  ${g} \in L^2(\Gamma)$ we consider three types of boundary conditions:
\begin{itemize}
  \item mixed such that $\Gamma = \Gamma_D \cup \Gamma_N $ and $\Gamma_D \cap \Gamma_N = \emptyset$
\begin{equation}
 \label{eq:mixed}
 \left\{
 \begin{array}{rcll}
 \boldsymbol{u} & = & 0 & \mbox{ on } \Gamma_D \\
 \boldsymbol{\sigma_n^{sym}} & = & 0 & \mbox{ on } \Gamma_N
 \end{array}
 \right. ;
\end{equation}
  \item tangential-displacement and normal-normal-stress (TDNNS)
  \begin{equation}
  \left\{
  \begin{array}{rcll}
  \label{eq:TDNNS}
   u_t & = & 0 & \mbox{ on } \Gamma \\
   {\sigma_{nn}^{sym}} & = & {g} & \mbox{ on } \Gamma 
  \end{array}
  \right. ;
  \end{equation}
  \item normal-displacement and tangential-normal-stress (NDTNS)
\begin{equation}
 \label{eq:NDTNS}
 \left\{
 \begin{array}{rcll}
 {u_n} & = & 0 & \mbox{ on } \Gamma \\
  {\sigma_{nt}^{sym}} & = & {g} & \mbox{ on } \Gamma
 \end{array}
 \right. .
\end{equation}
\end{itemize}
The second type of boundary condition has already been considered for linear elasticity equation in~\cite{MR2826472}.

%%%%%%%%%%%%%%%%%%%%%%%%%%%%%%%%%%%%%%%%%%%%%%%%%%%%%%%%%%%%%%%%%%%%%%%%%%%%%%%%%%%%%%%%%%%%%%%%%%%%%%%%%%%%
% Numerical methods
%%%%%%%%%%%%%%%%%%%%%%%%%%%%%%%%%%%%%%%%%%%%%%%%%%%%%%%%%%%%%%%%%%%%%%%%%%%%%%%%%%%%%%%%%%%%%%%%%%%%%%%%%%%%

\section{The numerical methods}
\label{sec:methods}

% In this section we discuss the numerical discretisations that we will be using in the numerical experiments. We start with the discrete notation.

%%%%%%%%%%%%%%%%%%%%%%%%%%%%%%%%%%%%%%%%%%%%%%%%%%
%% Triangulation
Let $\left\{\mathcal{T}_h\right\}_{h > 0}$ be a regular family of triangulations of $\bar{\Omega}$ made of simplices. For each triangulation $\mathcal{T}_h$, $\mathcal{E}_h$ denotes the set of its facets (edges for $d=2$, faces for $d=3$). In addition, for each element $K \in \mathcal{T}_h$, $h_K := \mbox{diam}(K)$, and we denote $h := \max_{K \in \mathcal{T}_h} h_K$. 
We define the following broken Sobolev spaces on %the triangulation $\mathcal{T}_h$ and 
the set of all edges in $\mathcal{E}_h$ (for $d=2$)
\begin{align*}
L^2(\mathcal{E}_h) & :=  \left\{v: \ v|_E \in L^2(E) \ \forall \ E \in \mathcal{E}_h \right\}. %\\
% H^m(\mathcal{T}_h) & :=  \left\{v \in L^2(\Omega): \ {v}|_K \in H^m(K) \ \forall \ K \in \mathcal{T}_h \right\} \mbox{ for } m \in \mathbb{N},
\end{align*}
% with the corresponding broken norms. 
Moreover, for $D \subset \Omega$, $\mathbb{P}_k (D)$ denotes the space of polynomials of total degree smaller than, or equal to,  $k$ on the set $D$.

We now  present the two discretisations to be used in the numerical experiments.

\subsection{Taylor-Hood discretisation}
\label{sec:TH}

We first consider the Taylor-Hood discretisation using the following approximation spaces
\begin{align*}
	\boldsymbol{TH^k_h} & := \left\{\boldsymbol{v_h} \in [H^1(\Omega)]^d: \quad \boldsymbol{v_h}|_K \in [\mathbb{P}_k(K)]^d  \quad \forall \ K \in \mathcal{T}_h\right\}, \\
	R^{k-1}_h & := \left\{q_h \in C^0(\bar{\Omega}): \quad q_h|_K \in \mathbb{P}_{k-1}(K)  \quad \forall \ K \in \mathcal{T}_h\right\}
\end{align*}
where $k \geq 2$ (see \cite[Chapter II, Section 4.2]{MR851383}).

If~\eqref{eq:stokes} is supplied with the homogeneous boundary conditions~\eqref{eq:Dirichlet}, then the discrete problem reads:

\begin{center}
\textit{Find $\left(\boldsymbol{u_h}, p_h\right) \in \left(\left(\boldsymbol{TH^k_h} \cap [H_0^1(\Omega)]^d\right)\right) \times \left(R^{k-1}_h \cap L_0^2(\Omega)\right)$ \\ s.t. for all $\left(\boldsymbol{v_h}, q_h\right) \in \left(\left(\boldsymbol{TH^k_h} \cap [H_0^1(\Omega)]^d\right)\right) \times \left(R^{k-1}_h \cap L_0^2(\Omega)\right)$}                                                                                                                                                                                                                                                                                                                                                \end{center}
\begin{equation}
\label{eq:TH_Stokes_Dirichlet}
 \left\{
		\begin{array}{rclcl}
   \displaystyle\int_{\Omega} \nu \nabla \boldsymbol{u_h} : \nabla \boldsymbol{v_h} \dx & - & \displaystyle\int_{\Omega} p_h \nabla \cdot \boldsymbol{v_h} \dx & = & \displaystyle\int_{\Omega}\boldsymbol{f}\boldsymbol{v_h} \dx\\
		& - & \displaystyle\int_{\Omega} \nabla \cdot \boldsymbol{u_h} q_h \dx & = & 0 .
   \end{array}
 \right.
\end{equation}
In case of TVNF boundary conditions~\eqref{eq:TVNF}, we define $\boldsymbol{V_t} := \left\{\boldsymbol{v} \in [H^1(\Omega)]^d: \ v_t = 0 \mbox{ on } \Gamma\right\}$, and the discrete problem reads:

\begin{center}
\textit{Find $\left(\boldsymbol{u_h}, p_h\right) \in \left(\boldsymbol{TH^k_h} \cap \boldsymbol{V_t}\right) \times R^{k-1}_h$ \\ s.t. for all $\left(\boldsymbol{v_h}, q_h\right) \in \left(\boldsymbol{TH^k_h} \cap \boldsymbol{V_t}\right) \times R^{k-1}_h$}                                                                                                                                                                                                                                                \end{center}
\begin{equation}
\label{eq:TH_Stokes_TVNF}
 \left\{
		\begin{array}{rclcl}
   \displaystyle\int_{\Omega} \nu \nabla \boldsymbol{u_h} : \nabla \boldsymbol{v_h} \dx & - & \displaystyle\int_{\Omega} p_h \nabla \cdot \boldsymbol{v_h} \dx & = & \displaystyle\int_{\Omega}\boldsymbol{f}\boldsymbol{v_h} \dx + \displaystyle\int_{\Gamma} g \left(\boldsymbol{v_h}\right)_t \ds\\
		& - & \displaystyle\int_{\Omega} \nabla \cdot \boldsymbol{u_h} q_h \dx & = & 0 .
   \end{array}
 \right.
\end{equation}
If NVTF boundary conditions~\eqref{eq:NVTF} are used, then we define the following space \\ $\boldsymbol{V_n} := \left\{ \boldsymbol{v} \in [H^1(\Omega)]^d: \ v_n = 0 \mbox{ on } \Gamma \right\}$, and the discrete problem reads:

\begin{center}
\textit{Find $\left(\boldsymbol{u_h}, p_h\right) \in \left(\boldsymbol{TH^k_h} \cap \boldsymbol{V_n}\right) \times \left(R^{k-1}_h \cap L_0^2(\Omega)\right)$ \\ s.t. for all $\left(\boldsymbol{v_h}, q_h\right) \in \left(\boldsymbol{TH^k_h} \cap \boldsymbol{V_n}\right) \times \left(R^{k-1}_h \cap L_0^2(\Omega)\right)$}                                                                                                                                                                                                                                                                                      \end{center}
\begin{equation}
\label{eq:TH_Stokes_NVTF}
 \left\{
		\begin{array}{rclcl}
   \displaystyle\int_{\Omega} \nu \nabla \boldsymbol{u_h} : \nabla \boldsymbol{v_h} \dx & - & \displaystyle\int_{\Omega} p_h \nabla \cdot \boldsymbol{v_h} \dx & = & \displaystyle\int_{\Omega}\boldsymbol{f}\boldsymbol{v_h} \dx + \displaystyle\int_{\Gamma} g \left(\boldsymbol{v_h}\right)_n \ds\\
		& - & \displaystyle\int_{\Omega} \nabla \cdot \boldsymbol{u_h} q_h \dx & = & 0 .
   \end{array}
 \right.
\end{equation}

In similar way, if the problem~\eqref{eq:elasticity} is supplied with the boundary conditions~\eqref{eq:mixed}, then the discrete problem reads

\begin{center}
\textit{Find $\left(\boldsymbol{u_h}, p_h\right) \in \left(\boldsymbol{TH^k_h} \cap [H_{\Gamma_D}^1(\Omega)]^d\right) \times R^{k-1}_h $ \\ s.t. for all $\left(\boldsymbol{v_h}, q_h\right) \in \left(\boldsymbol{TH^k_h} \cap [H_{\Gamma_D}^1(\Omega)]^d\right) \times R^{k-1}_h$}                                                                                                                                                                                                                                                                                        \end{center}
\begin{equation}
\label{eq:TH_elasticity_Dirichlet}
 \left\{
		\begin{array}{rclcl}
   \displaystyle\int_{\Omega} 2 \mu \boldsymbol{\varepsilon}\left(\boldsymbol{u_h}\right) : \boldsymbol{\varepsilon}\left(\boldsymbol{v_h}\right) \dx & - & \displaystyle\int_{\Omega} p_h \nabla \cdot \boldsymbol{v_h} \dx & = & \displaystyle\int_{\Omega}\boldsymbol{f}\boldsymbol{v_h} \dx\\
		-\displaystyle\int_{\Omega} \nabla \cdot \boldsymbol{u_h} q_h \dx & - & \frac{1}{\lambda} \displaystyle\int_{\Omega} p_h q_h \dx & = & 0 .
   \end{array}
 \right.
\end{equation}

The rest of the discrete problems associated with~\eqref{eq:elasticity} that is supplied with TDNNS boundary conditions~\eqref{eq:TDNNS} or NDTNS boundary conditions~\eqref{eq:NDTNS} are similar to \eqref{eq:TH_Stokes_TVNF} or~\eqref{eq:TH_Stokes_NVTF}, respectively.

\subsection{Hybrid discontinuous Galerkin discretisation}
\label{sec:hdg}

We restrict the discussion of this to two dimensional case $d=2$. This method has been presented and analysed in~\cite{barrenechea2016hybrid}. The velocity is approximated using the Brezzi-Douglas-Marini spaces (see~\cite[Section~2.3.1]{MR3097958}) of degree $k$ given by
\begin{align*}
 \boldsymbol{BDM_h^k} & := \left\{\boldsymbol{v_h} \in H\left(div, \Omega\right): \ \boldsymbol{v_h}|_K \in \left[\mathbb{P}_k\left(K\right)\right]^2 \ \forall \ {K \in \mathcal{T}_h}\right\}, \\
 \boldsymbol{BDM_{h,\tilde{\Gamma}}^k} & := \left\{\boldsymbol{v_h} \in H\left(div, \Omega\right): \ \boldsymbol{v_h}|_K \in \left[\mathbb{P}_k\left(K\right)\right]^2 \ \forall \ {K \in \mathcal{T}_h} \land \ {\left(\boldsymbol{v_h}\right)_n} = 0 \mbox{ on } \tilde{\Gamma}\right\},
\end{align*}
where $\tilde{\Gamma} \subset \partial \Omega$. If $\tilde{\Gamma} = \partial \Omega$, then $ \boldsymbol{BDM_{h,\tilde{\Gamma}}^k}$ is denoted $\boldsymbol{BDM_{h,0}^k}$.
% In addition, for $1 \leq m \leq k+1$ we denote $\Pi^k: [H^m(\Omega)]^2 \rightarrow \boldsymbol{BDM_h^k}$  the BDM projection defined in~\cite[Section~2.5]{MR3097958}.

The pressure is approximated in the space
\begin{align*}
 Q_h^{k-1} & := \left\{q_h \in L^2\left(\Omega\right): \ q_h|_K \in \mathbb{P}_{k-1}\left(K\right) \ \forall \ {K \in \mathcal{T}_h}\right\} . 
\end{align*} 
Finally, a Lagrange multiplier, aimed at approximating the tangential component of the velocity is introduced. The space swhere this multipliers is sought are given by %In order to propose a discretisation with fewer degrees of freedom, we discretise the Lagrange multiplier $\tilde{u}$ using the spaces
\begin{align*}
 M_{h}^{k-1} &:= \left\{\tilde{v}_h \in L^2\left(\mathcal{E}_h\right): \ \tilde{v}_h|_E \in \mathbb{P}_{k-1}\left(E \right) \ \forall \ {E \in \mathcal{E}_h} \right\}, \\ 
 M_{h,\tilde{\Gamma}}^{k-1} &:= \left\{\tilde{v}_h \in M_{h}^{k-1}: \tilde{v}_h = 0 \mbox{ on } \tilde{\Gamma} \right\},
\end{align*} 
where $\tilde{\Gamma} \subset \partial \Omega$. If $\tilde{\Gamma} = \partial \Omega$, then $M_{h,\tilde{\Gamma}}^{k-1}$ is denoted $M_{h,0}^{k-1}$.
 Furthermore, we introduce for all $E \in \mathcal{E}_h$ the $L^2(E)$-projection $\Phi^{k-1}_E: L^2\left(E\right) \rightarrow \mathbb{P}_{k-1}\left(E\right)$ defined by %as follows. For every $\tilde{w} \in L^2\left(E\right)$, $\Phi^{k-1}_E(\tilde{w})$ is the unique element of $\mathbb{P}_{k-1}\left(E\right)$ satisfying 
\begin{equation}
\label{eq:L2edge_projection}
\int_{E} \Phi^{k-1}_E(\tilde{w}) {\tilde{v}_h} \ds = \int_{E} \tilde{w} {\tilde{v}_h} \ds \quad \forall \ {\tilde{v}_h} \in \mathbb{P}_{k-1}\left(E\right),
\end{equation}
and we denote $\Phi^{k-1}: L^2\left(\mathcal{E}_h\right) \rightarrow M_{h}^{k-1}$ defined as $\Phi^{k-1}|_E := \Phi^{k-1}_E$ for all $E \in \mathcal{E}_h$. 
% Let us denote $\boldsymbol{V_h} := \boldsymbol{BDM_h^k} \times M_{h}^{k-1}$. 

If~\eqref{eq:stokes} is supplied with the homogeneous boundary conditions~\eqref{eq:Dirichlet}, then the discrete problem reads:

%%%%%%%%%%%%%%%%%%%%%%%%%%%%%%%%%%%%%%%%%%%%%%%%%%
%% Discrete problem
\begin{center}
\textit{Find $\left(\boldsymbol{u_h}, \tilde{u}_h, p_h\right) \in \boldsymbol{BDM_{h,0}^k} \times M_{h,0}^{k-1} \times \left(Q^{k-1}_h \cap L_0^2(\Omega)\right)$ \\ s.t. for all $\left(\boldsymbol{v_h}, \tilde{v}_h, q_h\right) \in \boldsymbol{BDM_{h,0}^k} \times M_{h,0}^{k-1} \times\left(Q^{k-1}_h \cap L_0^2(\Omega)\right)$,}                                                                                                                                                                                                                                                                                                                                                  \end{center}
\begin{equation}
\label{eq:hdg_Stokes_Dirichlet}
 \left\{
		\begin{array}{rclcl}
   a \left(\left(\boldsymbol{u_h}, \tilde{u}_h\right),\left(\boldsymbol{v_h}, \tilde{v}_h\right)\right) & + & b\left(\left(\boldsymbol{v_h}, \tilde{v}_h\right), p_h\right) & = & \displaystyle\int_{\Omega}\boldsymbol{f}\boldsymbol{v_h} \dx \\
		 & & b\left(\left(\boldsymbol{u_h}, \tilde{u}_h\right), q_h\right) & = & 0
   \end{array}
 \right. ,
\end{equation}
where %we define the velocity bilinear form $a: \boldsymbol{BDM_h^k} \times M_{h}^{k-1} \times \boldsymbol{BDM_h^k} \times M_{h}^{k-1} \rightarrow \mathbb{R}$ as
\begin{align}
\label{eq:a_form}
  \nonumber a \left(\left(\boldsymbol{w_h}, {\tilde{w}_h}\right), \left(\boldsymbol{v_h}, \tilde{v}_h\right)\right) := & \sum_{K \in \mathcal{T}_h} \left(\int_K \nu \nabla \boldsymbol{w_h} : \nabla \boldsymbol{v_h} \dx \right. \\ \nonumber
& - \int_{\partial K} \nu \left(\boldsymbol{\partial_n {w_h}}\right)_t \big(\left(\boldsymbol{v_h}\right)_t - \tilde{v}_h \big) \ds \\
& - \int_{\partial K} \nu \big(\left(\boldsymbol{w_h}\right)_t - \tilde{w}_h \big) \left(\boldsymbol{\partial_n {v_h}}\right)_t \ds \\ \nonumber
 &  \left. + \nu \frac{\tau}{h_K} \int_{\partial K} \Phi^{k-1}\big(\left(\boldsymbol{w_h}\right)_t - \tilde{w}_h \big) \Phi^{k-1}\big(\left(\boldsymbol{v_h}\right)_t - \tilde{v}_h \big) \ds \right) ,
\end{align}
and $\tau > 0$ is a stabilisation parameter, 
%%%%%%%%%%%%%%%%%%%%%%%%%%%%%%%%%%%%%%%%%%%%%%%%%%
%% Bilinear form for pressure part
and %$b: \boldsymbol{BDM_h^k} \times M_{h}^{k-1} \times Q_h^{k-1} \rightarrow \mathbb{R}$ as
\begin{equation}
\label{eq:b_form}
 b\left(\left(\boldsymbol{v_h}, \tilde{v}_h\right), q_h\right) := - \sum_{K \in \mathcal{T}_h} \int_K q_h \nabla \cdot \boldsymbol{v_h} \dx.
\end{equation}
If TVNF boundary conditions~\eqref{eq:TVNF} are used, then the discrete problem reads:

\begin{center}
\textit{Find $\left(\boldsymbol{u_h}, \tilde{u}_h, p_h\right) \in \boldsymbol{BDM_{h}^k} \times M_{h,0}^{k-1} \times Q^{k-1}_h$ \\ s.t. for all $\left(\boldsymbol{v_h}, \tilde{v}_h, q_h\right) \in \boldsymbol{BDM_{h}^k} \times M_{h,0}^{k-1} \times Q^{k-1}_h$,}                                                                                                                                                                                                                                                                     \end{center}
\begin{equation}
\label{eq:hdg_Stokes_TVNF}
 \left\{
		\begin{array}{rclcl}
   a \left(\left(\boldsymbol{u_h}, \tilde{u}_h\right),\left(\boldsymbol{v_h}, \tilde{v}_h\right)\right) & + & b\left(\left(\boldsymbol{v_h}, \tilde{v}_h\right), p_h\right) & = & \displaystyle\int_{\Omega}\boldsymbol{f}\boldsymbol{v_h} \dx + \int_{\Gamma} g \left(\boldsymbol{v_h}\right)_n \ds\\
		 & & b\left(\left(\boldsymbol{u_h}, \tilde{u}_h\right), q_h\right) & = & 0
   \end{array}
 \right. .
\end{equation}
In case of NVTF boundary conditions~\eqref{eq:NVTF}, the discrete problem reads:

\begin{center}
\textit{Find $\left(\boldsymbol{u_h}, \tilde{u}_h, p_h\right) \in \boldsymbol{BDM_{h,0}^k} \times M_h^{k-1} \times \left(Q^{k-1}_h \cap L_0^2(\Omega)\right)$ \\ s.t. for all $\left(\boldsymbol{v_h}, \tilde{v}_h, q_h\right) \in \boldsymbol{BDM_{h,0}^k} \times M_h^{k-1} \times \left(Q^{k-1}_h \cap L_0^2(\Omega)\right)$,}                                                                                                                                                                                                                                                                                                       \end{center}
\begin{equation}
\label{eq:hdg_Stokes_NVTF}
 \left\{
		\begin{array}{rclcl}
   a \left(\left(\boldsymbol{u_h}, \tilde{u}_h\right),\left(\boldsymbol{v_h}, \tilde{v}_h\right)\right) & + & b\left(\left(\boldsymbol{v_h}, \tilde{v}_h\right), p_h\right) & = & \displaystyle\int_{\Omega}\boldsymbol{f}\boldsymbol{v_h} \dx + \int_{\Gamma} g \tilde{v}_h \ds\\
		 & & b\left(\left(\boldsymbol{u_h}, \tilde{u}_h\right), q_h\right) & = & 0
   \end{array}
 \right. .
\end{equation}

In similar way, if the problem~\eqref{eq:elasticity} is supplied with the mixed boundary conditions~\eqref{eq:mixed}, then the discrete problem reads:

\begin{center}
\textit{Find $\left(\boldsymbol{u_h}, \tilde{u}_h, p_h\right) \in \boldsymbol{BDM_{h,\Gamma_D}^k} \times M_{h,\Gamma_D}^{k-1} \times Q^{k-1}_h$ \\ s.t. for all $\left(\boldsymbol{v_h}, \tilde{v}_h, q_h\right) \in \boldsymbol{BDM_{h,\Gamma_D}^k} \times M_{h,\Gamma_D}^{k-1} \times Q^{k-1}_h$}                                                                                                                                                                                                                                                                                                              \end{center}
\begin{equation}
\label{eq:hdg_elasticity}
 \left\{
		\begin{array}{rclcl}
   a_{s} \left(\left(\boldsymbol{u_h}, \tilde{u}_h\right),\left(\boldsymbol{v_h}, \tilde{v}_h\right)\right) & + & b\left(\left(\boldsymbol{v_h}, \tilde{v}_h\right), p_h\right) & = & \displaystyle\int_{\Omega}\boldsymbol{f}\boldsymbol{v_h} \dx \\
		 b\left(\left(\boldsymbol{u_h}, \tilde{u}_h\right), q_h\right) & + & c(p_h, q_h) & = & 0,
   \end{array}
 \right.
\end{equation}
where
\begin{align}
\label{eq:a_sym_form}
  \nonumber a_{s} \left(\left(\boldsymbol{w_h}, {\tilde{w}_h}\right), \left(\boldsymbol{v_h}, \tilde{v}_h\right)\right) := & \sum_{K \in \mathcal{T}_h} \left(\int_K 2 \mu \boldsymbol{\varepsilon}(\boldsymbol{w_h}) : \boldsymbol{\varepsilon}(\boldsymbol{v_h}) \dx \right. \\ \nonumber
& - \int_{\partial K} 2 \mu \left(\boldsymbol{\varepsilon_n} (\boldsymbol{w_h})\right)_t \big(\left(\boldsymbol{v_h}\right)_t - \tilde{v}_h \big) \ds \\
& - \int_{\partial K} 2 \mu \big(\left(\boldsymbol{w_h}\right)_t - \tilde{w}_h \big) \left(\boldsymbol{\varepsilon_n} (\boldsymbol{v_h})\right)_t \ds \\ \nonumber
 & \  \left. + 2 \mu \frac{\tau}{h_K} \int_{\partial K} \Phi^{k-1}\big(\left(\boldsymbol{w_h}\right)_t - \tilde{w}_h \big) \Phi^{k-1}\big(\left(\boldsymbol{v_h}\right)_t - \tilde{v}_h \big) \ds \right) ,
\end{align}
$b$ is defined by~\eqref{eq:b_form}, and 
\begin{equation*}
	c(r_h, q_h) := - \frac{1}{\lambda} \int_{\Omega} r_h q_h \ds.
\end{equation*}

The rest of the discrete problems associated with~\eqref{eq:elasticity} that is supplied with TDNNS boundary conditions~\eqref{eq:TDNNS} or NDTNS boundary conditions~\eqref{eq:NDTNS} are similar to~\eqref{eq:hdg_Stokes_TVNF} or~\eqref{eq:hdg_Stokes_NVTF}, respectively.

%%%%%%%%%%%%%%%%%%%%%%%%%%%%%%%%%%%%%%%%%%%%%%%%%%%%%%%%%%%%%%%%%%%%%%%%%%%%%%%%%%%%%%%%%%%%%%%%%%%%%%%%%%%%
% Domain decomposition
%%%%%%%%%%%%%%%%%%%%%%%%%%%%%%%%%%%%%%%%%%%%%%%%%%%%%%%%%%%%%%%%%%%%%%%%%%%%%%%%%%%%%%%%%%%%%%%%%%%%%%%%%%%%

\section{The domain decomposition preconditioners}
\label{sec:dd}

Let us assume that we have to solve the following linear system $\mathbf{A} \boldsymbol{U} = \boldsymbol{F}$
where $\mathbf{A}$ is the matrix arising from discretisation of the Stokes or linear elasticity 
equation on the domain $\Omega$, $\boldsymbol{U}$ is the vector of
unknowns, and $\boldsymbol{F}$ is the right hand side. To accelerate the
performance of an iterative Krylov method~\cite[Chapter~3]{MR3450068} applied to this system we
will consider domain decomposition preconditioners which are naturally
parallel. They are based on an overlapping decomposition of the
computational domain.

\begin{figure}[!ht]
\centering
    \subfloat[L-shaped domain]{
      \includegraphics[width=0.46\textwidth]{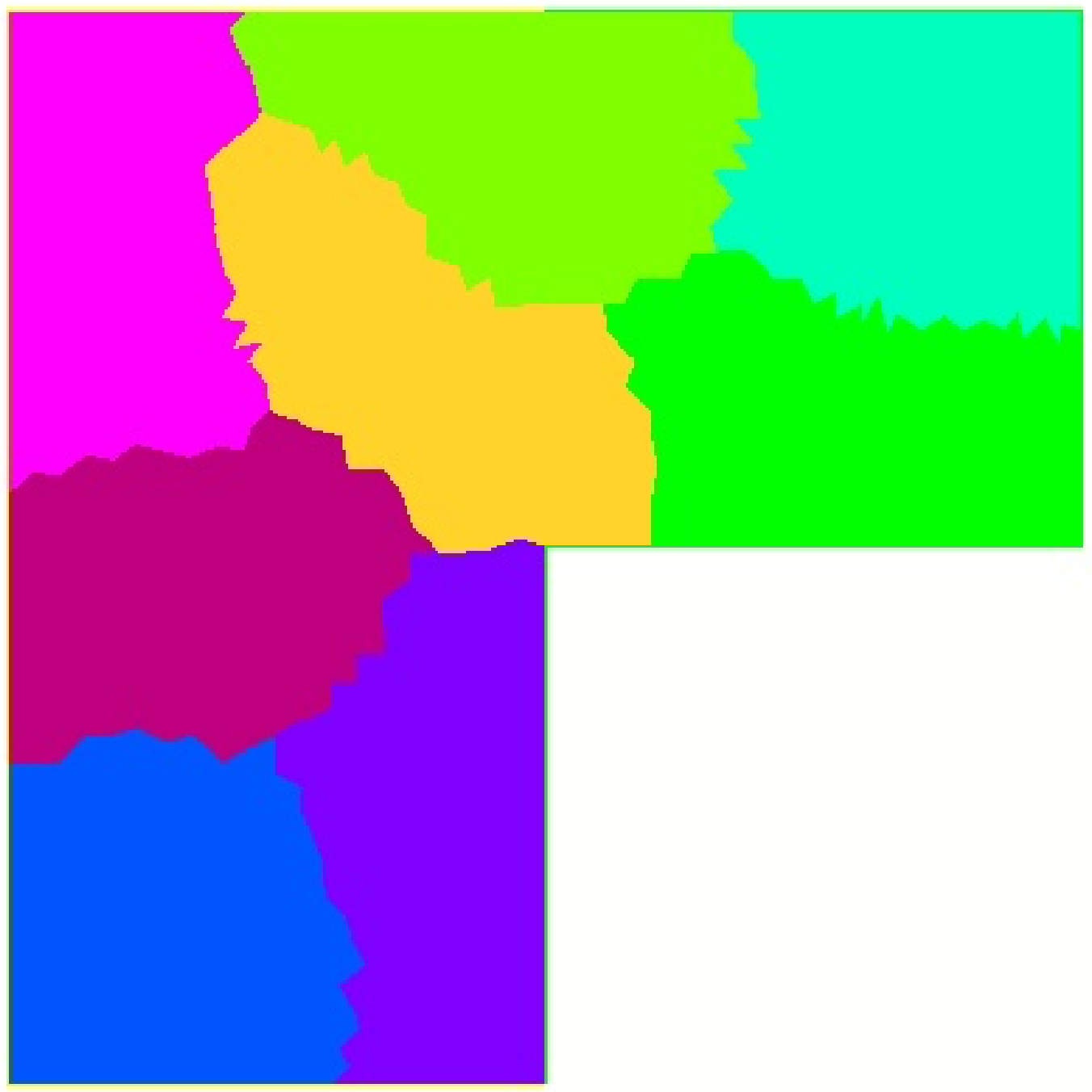}
      }  \hfill
    \subfloat[T-shaped domain]{
      \includegraphics[width=0.36\textwidth]{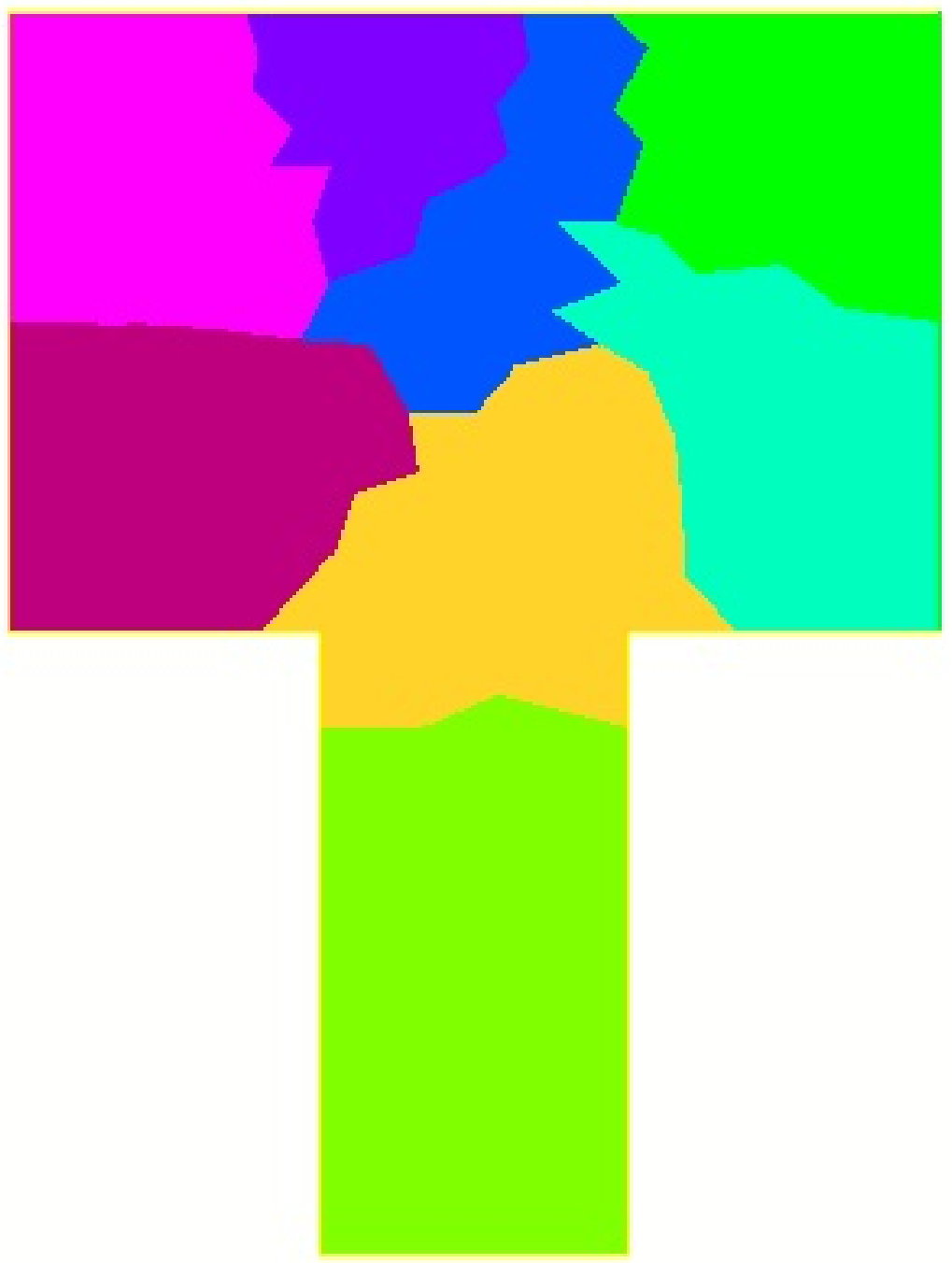}
}
\caption{Partition of the domain for 8 subdomains by METIS.}
\label{fig:Partition}
  \end{figure}

Let $\{\mathcal{T}_{h,i}\}_{i=1}^N$ be a partition of the
triangulation $\mathcal{T}_h$ (see examples in Figure~\ref{fig:Partition}). For an integer value $l \geq 0$, we
define an overlapping decomposition
$\{\mathcal{T}_{h,i}^{l}\}_{i=1}^N$ such that $\mathcal{T}_{h,i}^{l}$
is a set of all triangles from $\mathcal{T}_{h,i}^{l-1}$ and all
triangles from $\mathcal{T}_h \setminus \mathcal{T}_{h,i}^{l-1}$ that
have non-empty intersection with $\mathcal{T}_{h,i}^{l-1}$, and
$\mathcal{T}_{h,i}^0 = \mathcal{T}_{h,i}$. With this definition the width of
the overlap will be of $2l$. Furthermore, if $W_h$
stands for the finite element space associated to $\mathcal{T}_h$,
$W_{h,i}^{l}$ is the local finite element spaces on
$\mathcal{T}_{h,i}^{l}$, which is a triangulation of $\Omega_i$. 

Let $\mathcal{N}$ be the set of indices of degrees of
freedom of $W_h$ and $\mathcal{N}_i^{l}$ the set of indices of degrees of freedom of $W_{h,i}^{l}$ for $l \geq 0$. Moreover, we define the restriction operator $\mathbf{R}_i: W_h \rightarrow W_{h,i}^{l}$ as a rectangular matrix $|\mathcal{N}_i^{l}| \times |\mathcal{N}|$ such that if $\boldsymbol{V}$ is the vector of degrees of freedom of $v_h \in W_h$, then $\mathbf{R}_i \boldsymbol{V}$ is the vector of degrees of freedom of $W_{h,i}^{l}$ in $\Omega_i$. The extension operator from $W_{h,i}^{l}$ to $W_h$ and its associated matrix are both given by $\mathbf{R}_i^T$. In addition we introduce a partition of unity $\mathbf{D}_i$ as a diagonal matrix $|\mathcal{N}_i^{l}| \times |\mathcal{N}_i^{l}|$ such that
\begin{equation}
\mathbf{Id} = \sum_{i=1}^N \mathbf{R}_i^T \mathbf{D}_i \mathbf{R}_i,
\end{equation}
where $\mathbf{Id} \in \mathbb{R}^{|\mathcal{N}| \times |\mathcal{N}|}$ is the identity matrix.

We first recall the Modified Restricted Additive Schwarz (MRAS) preconditioner introduced in~\cite{barrenechea2016hybrid} for the Stokes equation. This preconditioner is given by
\begin{equation}
\label{eq:MRAS}
\mathbf{M}^{-1}_{MRAS} = \sum_{i=1}^N \mathbf{R}_i^T \mathbf{D}_i \mathbf{B}_i^{-1} \mathbf{R}_i,
\end{equation}
where $\mathbf{B}_i$ is the matrix associated to a discretisation of Stokes equation~\eqref{eq:stokes} in $\Omega_i$ where we impose either TVNF~\eqref{eq:TVNF} or NVTF~\eqref{eq:NVTF} boundary conditions in $\Omega_i$. In case of a discretisation of elasticity equation~\eqref{eq:elasticity} in $\Omega_i$, we impose either TDNNS~\eqref{eq:TDNNS} or NDTNS~\eqref{eq:NDTNS} boundary conditions in $\Omega_i$.

We new introduce a symmetrised variant of~\eqref{eq:MRAS} called Symmetrised Modified Restricted Additive Schwarz (SMRAS), given by
\begin{equation}
\label{eq:SMRAS}
\mathbf{M}^{-1}_{SMRAS} = \sum_{i=1}^N \mathbf{R}_i^T \mathbf{D}_i \mathbf{B}_i^{-1} \mathbf{D}_i \mathbf{R}_i.
\end{equation}

\subsection{Two-level methods}
\label{sec:two_level}

A two-level version of the SMRAS and MRAS preconditioners will be based on a spectral coarse space obtained by solving the following local generalised eigenvalue problems

\begin{center}
\textit{Find $\left(\boldsymbol{V}_{jk}, \lambda_{jk}\right) \in \mathbb{R}^{|\mathcal{N}_j|} \setminus \{0\} \times \mathbb{R}$ s.t.}                                                                                                                                           \end{center}
\begin{equation}
\label{eq:coarse_second_GenEO}
	\mathbf{\tilde{A}}_j \boldsymbol{V}_{jk} = \lambda_{jk} \mathbf{B}_j \boldsymbol{V}_{jk},
\end{equation}
where $\mathbf{\tilde{A}}_j$ are local matrices associated to a discretisation of local Neumann boundary value problem in $\Omega_j$. Let $\theta > 0$ be a user-defined threshold. We define $Z_{GenEO} \subset \mathbb{R}^{|\mathcal{N}|}$ as the vector space spanned by the family of vectors $\left(\mathbf{R}_j^T \mathbf{D}_j \boldsymbol{V}_{jk} \right)_{\lambda_{jk} < \theta}$, $1 \leq j \leq N$, corresponding to eigenvalues smaller than $\theta$. The value of $\theta$ is chosen such that for a given problem and preconditioner, the behaviour of the method should be robust in the sense that, its convergence should not depend, or depends very weakly, on the number of subdomains. 

We are now ready to introduce the two-level method with coarse space
% \begin{equation}
% 	\label{eq:coarse_GenEO_CS}
% Z_{GenEO} := Z^{\gamma}_{geneo} \oplus 
$Z_{GenEO}$. 
% \end{equation}
Let $\mathbf{P}_0$ be the $\mathbf{A}$-orthogonal projection onto the coarse space $Z_{GenEO}$. The two-level SMRAS preconditioner is defined as
\begin{equation}
	\label{eq:coarse_two_level_SMRAS}
\mathbf{M}^{-1}_{SMRAS,2} = \mathbf{P}_0 \mathbf{A}^{-1} + (\mathbf{Id} - \mathbf{P}_0)\mathbf{M}^{-1}_{SMRAS}(\mathbf{Id} - \mathbf{P}_0^T).
\end{equation}
Furthermore, if $\mathbf{R}_0$ is a matrix whose rows are a basis of the coarse space $Z_{GenEO}$, then
\begin{equation*}
	\mathbf{P}_0 \mathbf{A}^{-1} = \mathbf{R}_0^T \left(\mathbf{R}_0 \mathbf{A R}_0^T \right)^{-1} \mathbf{R}_0.
\end{equation*}
In similar way, we can introduce the two-level MRAS preconditioner
\begin{equation}
	\label{eq:coarse_two_level_MRAS}
\mathbf{M}^{-1}_{MRAS,2} = \mathbf{P}_0 \mathbf{A}^{-1} + (\mathbf{Id} - \mathbf{P}_0)\mathbf{M}^{-1}_{MRAS}(\mathbf{Id} - \mathbf{P}_0^T).
\end{equation}

%%%%%%%%%%%%%%%%%%%%%%%%%%%%%%%%%%%%%%%%%%%%%%%%%%%%%%%%%%%%%%%%%%%%%%%%%%%%%%%%%%%%%%%%%%%%%%%%%%%%%%%%%%%%
% 2D numerics
%%%%%%%%%%%%%%%%%%%%%%%%%%%%%%%%%%%%%%%%%%%%%%%%%%%%%%%%%%%%%%%%%%%%%%%%%%%%%%%%%%%%%%%%%%%%%%%%%%%%%%%%%%%%

\section{Numerical results for two dimensional problems}
\label{sec:2D_numerics}

In this section we assess the performance of the preconditioners defined in Section~\ref{sec:two_level}. We will compare the newly introduced ones with that of ORAS and SORAS introduced in~\cite{haferssas:hal-01278347}. These kind of preconditioners are associated with the Robin interface conditions and require an optimised parameter as it can be seen in~\eqref{eq:coarse_elasticity_Robin_BC} below. The big advantage of SMRAS and MRAS preconditioners from the previous section is that they are parameter-free.
% In this section we will compare a two-level SMRAS preconditioner with the coarse space $Z_{GenEO}$ presented above with the two-level SORAS preconditioner with the coarse space $Z_{GenEO}$ presented in~\cite{haferssas:hal-01278347}. Moreover, we compare the results with the one-level methods.
We consider the partial differential equation model for nearly incompressible elasticity and Stokes flow as problems of similar mixed formulation. Each of these problems is discretised by using the Taylor-Hood methods from Section~\ref{sec:TH} and the hdG discretisation from Section~\ref{sec:hdg}.

Our experiments will be based on the classical weak scaling test. This test is built as follows. A domain $\bar{\Omega}$ is split into a triangulation $\mathcal{T}_h$. For each of element $K \in \mathcal{T}_h$, $h_K = \mbox{diam}(K)$, and we denote the mesh size by $h := \max_{K \in \mathcal{T}_h} h_K$. Then, this triangulation is split into overlapping subdomains of size H, in such a way $\frac{H}{h}$ remains constant. In the absence of a second level in the preconditioner, if the number of subdomains grows, then the convergence gets slower. A coarse space provides a global information and leads to a more robust behaviour.
% , where $H$ is maximum diameter of the subdomains and $h := max_{K \in \mathcal{T}_h} h_K$.

The simplest way to build a coarse space is to consider the zero energy modes. More precisely they are the eigenvectors associated with the zero eigenvalues of~\eqref{eq:coarse_second_GenEO} on a floating subdomains. Hence, by a floating subdomain we mean a subdomain without Dirichlet boundary condition on any part of the boundary. Then the matrix on the left hand side of~\eqref{eq:coarse_second_GenEO} is singular and there are zero eigenvalues.
% For a floating subdomain, this is a domain whose boundary conditions do not include a Dirichlet type, the eigenvalue problem~\eqref{eq:coarse_second_GenEO} has zero eigenvalues. In fact, the left hand side is a matrix of a disccretisation of the local problem~\eqref{eq:stokes}, or~\eqref{eq:elasticity}, with Neumann boundary conditions. The eigenvectors associated to the zero eigenvalues are known as zero energy modes, and they are the most simple choice to build a coarse space. 
%In case of the linear elasticity equation we have three in two dimensional space and six in three dimensional space such zeros associated with the rigid body modes, that are translations and rotations. The Stokes problem has two in two dimensional space and three in three dimensional space zero energy modes that corresponds to the translations in each of directions. 
These zero energy modes are the rigid body motions (three in two dimensions, six in three dimensions) for the elasticity problem, and the constants (two in two dimensions, three in three dimensions) for the Stokes equations. Unfortunately for some cases, this choice is not sufficient, so we have collected the smallest $M$ eigenvalues for each subdomain and build a coarse space by including the eigenvectors associated to them. The different values of $M$ are presented in the table in brackets.

All experiments have been made by using FreeFem++~\cite{MR3043640}, which is a free software specialised in variational discretisations of partial differential equations. We use GMRES~\cite{MR848568} as an iterative solver. Generalized eigenvalue problems to generate the coarse space are solved using ARPACK~\cite{MR1621681}. %In all cases, they are used in conjunction with a Krylov iterative solver like GMRES. We define the coarse space by collecting the eigenvectors resulting from the solutions of the local eigenvalue problems. The easiest way to build a coarse space is to incorporate zero energy modes in it. By the zero energy mode we understand the eigenvector associated with the zero eigenvalue on the floating subdomain, that is subdomain that do not have any Dirichlet boundary conditions. 
The overlapping decomposition into subdomains can be uniform (Unif) or generated by METIS (MTS)~\cite{karypis1998software}. In each of the examples in this section we consider decomposition with two layers of mesh size $h$ in the overlap. Tables show the number of iterations needed to achieve a relative $l^2$ norm of the error smaller than $10^{-6}$, $ \frac{\|\boldsymbol{U} - \boldsymbol{U}_n\|_{l^2}}{\|\boldsymbol{U} - \boldsymbol{U}_0\|_{l^2}} < 10^{-6}$, where $\boldsymbol{U}$ is the solution of global problem given by direct solver and $\boldsymbol{U}_m$ denotes $m$-th iteration of the iterative solver. In addition, $DOF$ stands for number of degrees of freedom and $N$  for the number of subdomains in all tables.

%%%%%%%%%%%%%%%%%%%%%%%%%%%%%%%%%%%%%%%%%%%%%%%%%%%%%%%%%%%%%%% Taylor-Hood %%%%%%%%%%%%%%%%%%%%%%%%%%%%%%%%%%%%%%%%%%%%%%%%%%%%%%%%%%%%%%%%
\subsection{Taylor-Hood discretisation}
\label{sec:coarse_TH}

In this section we consider the Taylor-Hood discretisation from Section~\ref{sec:TH} with different values of $k \geq 2$ for nearly incompressible elasticity and Stokes equations. 

\subsubsection{Nearly incompressible elasticity}
Since we consider the preconditioners with various interface conditions we need to comment the way of imposing them. ORAS and SORAS preconditioners follow~\cite{haferssas:hal-01278347} and use Robin interface conditions. This means, the weak formulation of the linear elasticity problem contains the following term
\begin{equation}
\label{eq:coarse_elasticity_Robin_BC}
\int_{\partial \Omega_i \setminus \Gamma} \boldsymbol{\sigma}^{sym}_{\boldsymbol{n}}\left(\boldsymbol{v_h}\right)_n \ds + \int_{\partial \Omega_i \setminus  \Gamma} 2 \alpha \frac{\mu (2 \mu + \lambda)}{\lambda +3 \mu} \boldsymbol{u_h} \boldsymbol{v_h} \ds
\end{equation}
where again, following~\cite{haferssas:hal-01278347} we choose $\alpha = 10$. 
Fortunately, the MRAS and SMRAS preconditioners are parameter-free. %We already know that for the Stokes equation these preconditioners are associated with NVTF and TVNF interface conditions. In the case of nearly incompressible elasticity we refer to them as normal-displacement and tangential-normal-stress (NDTNS) and tangential-displacement and normal-normal-stress (TDNNS) interface conditions. The second type of the boundary conditions has been already introduced for linear elasticity equation in~\cite{MR2826472}.
For all numerical experiments associated in this section we use zero as an initial guess for the GMRES iterative solver.  Moreover, the overlapping decomposition into subdomains is generated by METIS. 
% \newpage
\begin{figure}[!ht]
\centering
    \subfloat[Domain and boundary\label{fig:L_shape}]{\resizebox {0.4\columnwidth} {!}{
	  \begin{tikzpicture}
	  \draw[draw=red]  (0,-1) -- (-1,-1) -- (-1,0) node[red,right] {\tiny $\Gamma_D$} -- (-1,1) -- (0,1) ;
	  \draw[blue, dashed] (0,1) -- (1,1) -- (1,0) -- (0.5,0) node[above] {\tiny $\Gamma_N$}  -- (0,0)  -- (0,-1) ;
    \end{tikzpicture}
    }}
    \hfill
    \subfloat[Discrete solution\label{fig:L_shape_solution}]{
      \includegraphics[width=0.4\textwidth]{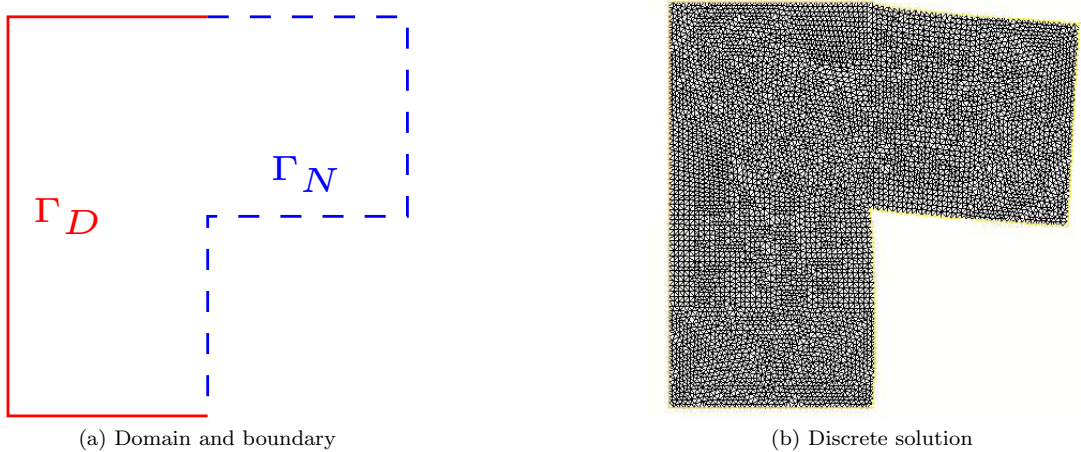}
}
\caption{L-shaped domain problem.}
  \end{figure}

\begin{figure}[!ht]
\centering
    \subfloat[Robin interface conditions]{
      \includegraphics[width=0.3\textwidth]{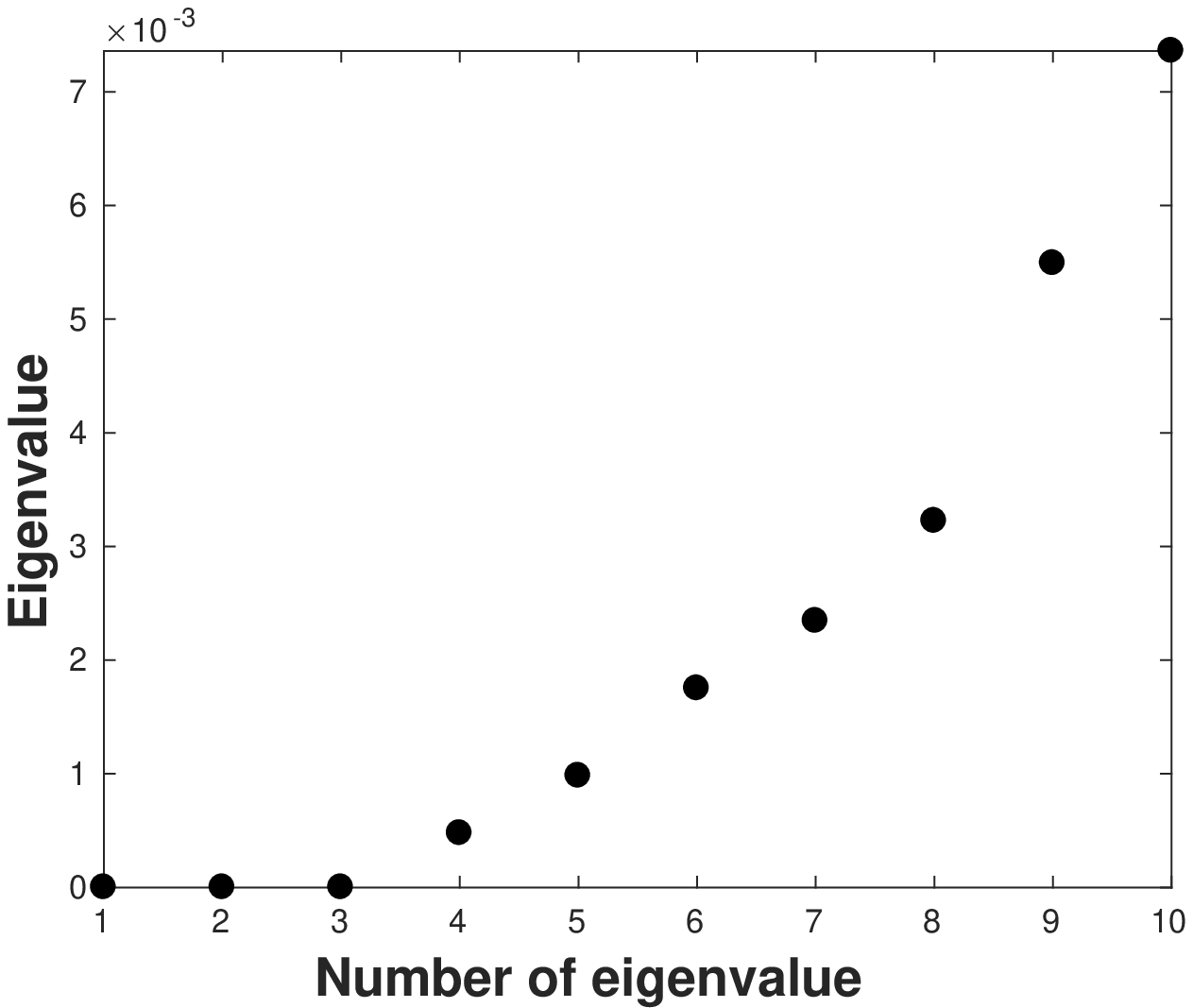}
    }
     \subfloat[NDTNS interface conditions]{
      \includegraphics[width=0.3\textwidth]{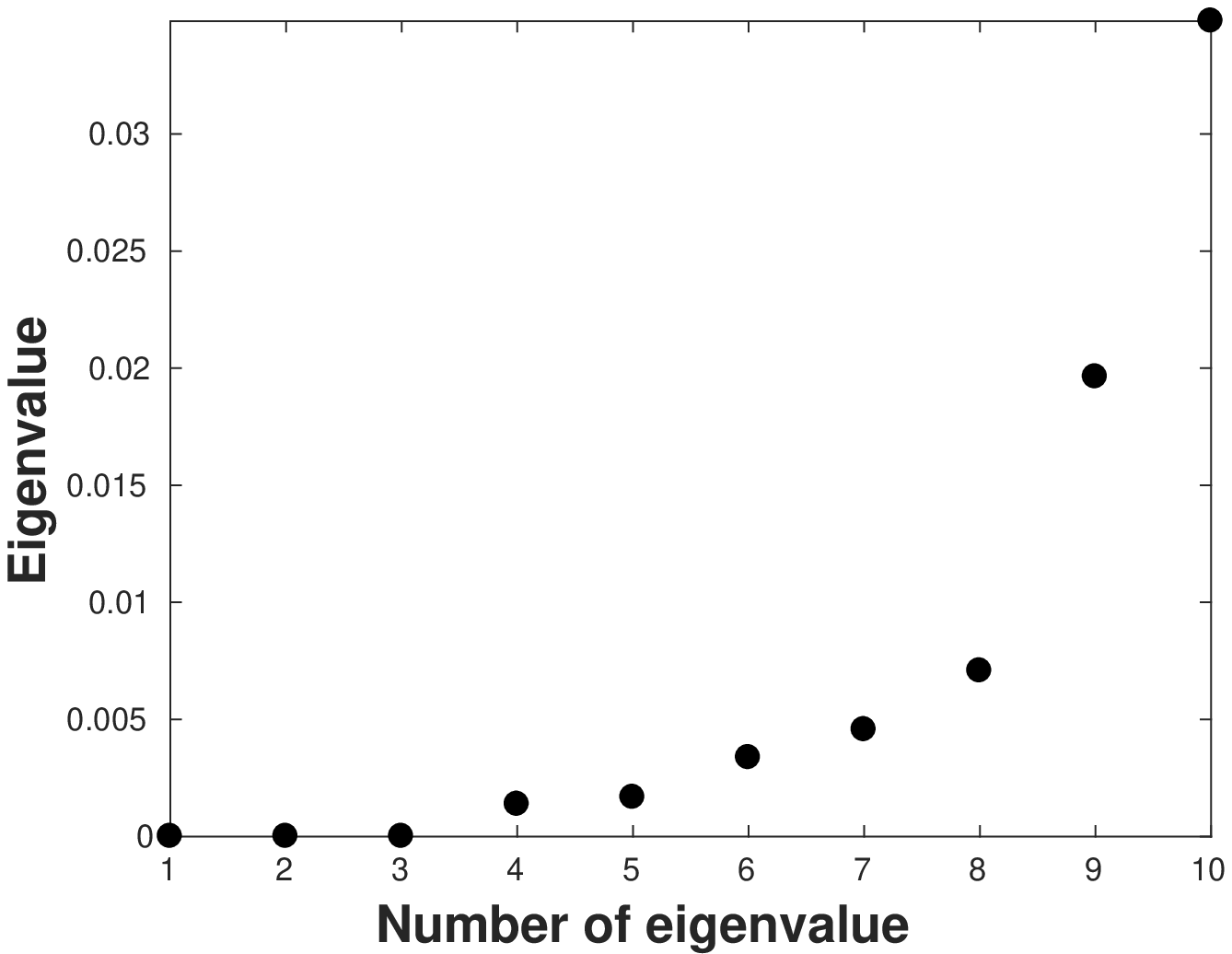}
    }
     \subfloat[TDNNS interface conditions]{
      \includegraphics[width=0.3\textwidth]{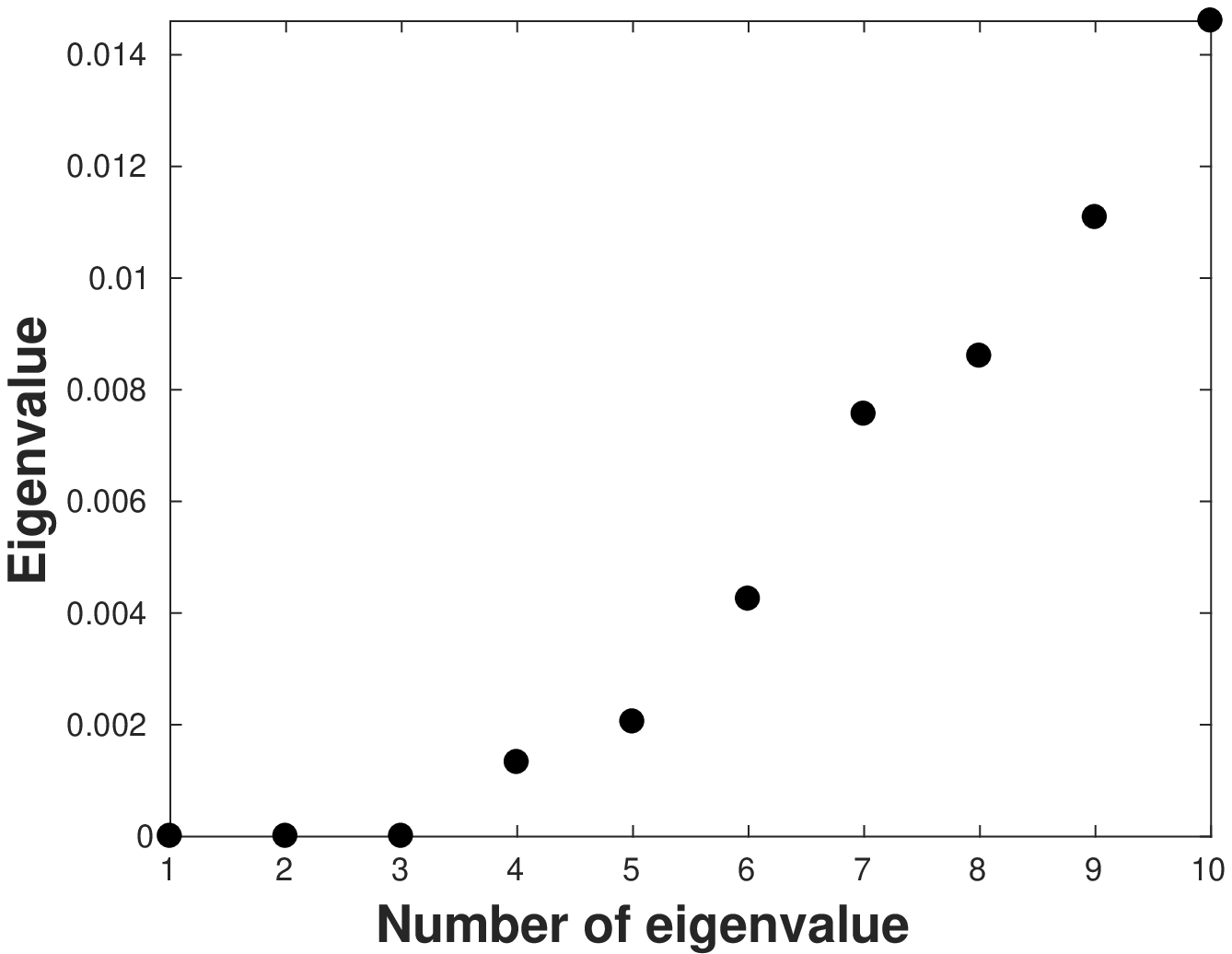}
    }
 	\caption{Eigenvalues on one of the floating subdomains in case of METIS decomposition and Taylor-Hood discretisation ($\boldsymbol{TH^3_h}, R^{2}_h$) - the heterogeneous beam.}
	\label{fig:TH_beam_eigenvalues}
  \end{figure}

\begin{testcase}[The L-shaped domain problem]
\label{exp:L_shape}
{\rm We consider the L-shaped domain $\Omega = (-1,1)^2 \setminus \left\{(0,1) \times (-1,0)\right\}$ clamped on the left side and partly from a top and bottom as it is depicted in Figure~\ref{fig:L_shape}. This example is similar to the one in~\cite{MR3407237}. The associated boundary value problem is
\begin{equation}
\label{eq:L_shape_elasticity}
\left\{
\begin{array}{rclclr}
 -2 \mu \nabla \cdot \boldsymbol{\varepsilon}(\boldsymbol{u}) & + & \nabla p & = & (0,-1)^T & \mbox{in } \Omega \\
 & -&  \nabla \cdot \boldsymbol{u} & = & \frac{1}{\lambda} p & \mbox{in } \Omega \\
  	& & \boldsymbol{u}(x,y) & = &(0,0)^T & \mbox{on } \Gamma_D \\
	& & \boldsymbol{\sigma}^{sym}_{\boldsymbol{n}}(x,y)  & = & (0,0)^T & \mbox{on } \Gamma_N
 \end{array}
\right. .
\end{equation} 
The physical parameters are $E = 10^5$ and $\nu = 0.4999$ (nearly incompressible). In Figure~\ref{fig:L_shape_solution} we plot the mesh of the bent domain. %This solution is obtained by solving the problem numerically and multiplying by huge constant to expose the bend. 

% Table~\ref{tab:hdG_Cavity_elasticity} shows the number of iterations needed to
% achieve a relative $L^2$ norm of the error (with respect to the one domain
% solution) smaller than $10^{-6}$.
We choose $k=3$ for the Taylor-Hood discretisation. In Figure~\ref{fig:TH_Cavity_elasticity_eigenvalues} we plot the eigenvalues of one floating subdomain. The clustering of small eigenvalues of the generalised eigenvalue problem defined in~\eqref{eq:coarse_second_GenEO} suggests the number of eigenvectors to be added to the coarse space. The three zero eigenvalues correspond to the zero energy modes.

% Table~\ref{tab:TH_Cavity_elasticity} shows the number of iterations needed to achieve a relative $L^2$ norm of the error (with respect to the one domain solution) smaller than $10^{-6}$.
\begin{table}[!ht]
\caption{Comparison of preconditioners for Taylor-Hood discretisation ($\boldsymbol{TH^3_h}, R^{2}_h$)  - the L-shaped domain problem.}
\label{tab:TH_L_shape_elasticity}
\begin{center}
\begin{adjustbox}{max width=\textwidth}
\begin{tabular}{ c c | c c c c c c }
  & & \multicolumn{6}{ c }{\textbf{One-level}} \\
  \textbf{DOF} & \textbf{N} & {\textbf{ORAS}} & {\textbf{SORAS}} & {\textbf{NDTNS-MRAS}} & {\textbf{NDTNS-SMRAS}} & {\textbf{TDNNS-MRAS}} & {\textbf{TDNNS-SMRAS}} \\
  \hline
		\textbf{124 109} & \textbf{4} & 26 & 60 & 26 & 60 & 30 & 59 \\
		\textbf{478 027} & \textbf{16} & 57 & 131 & 69 & 143 & 65 & 140 \\
		\textbf{933 087} & \textbf{32} & 84 & 180 & 109 & 221 & 104 & 211 \\
		\textbf{1 899 125 } & \textbf{64} & 130 & 293 & 181 & 362 & 161 & 312 \\
		\textbf{3 750 823} & \textbf{128} & 209 & 412 & 302 & 568 & 251 & 510 \\
		%\textbf{7 715 413} & \textbf{256} &  &  &  &  &  &  \\ \\
  & & \multicolumn{6}{ c }{\textbf{Two-level (3 eigenvectors)}} \\
  \textbf{DOF} & \textbf{N} & {\textbf{ORAS}} & {\textbf{SORAS}} & {\textbf{NDTNS-MRAS}} & {\textbf{NDTNS-SMRAS}} & {\textbf{TDNNS-MRAS}} & {\textbf{TDNNS-SMRAS}} \\
  \hline
		\textbf{124 109} & \textbf{4} & 18 & 40 & 19 & 36 & 24 & 41 \\
		\textbf{478 027} & \textbf{16} & 37 & 52 & 40 & 57 & 46 & 56 \\
		\textbf{933 087} & \textbf{32} & 49 & 57 & 56 & 67 & 53 & 66 \\
		\textbf{1 899 125 } & \textbf{64} & 65 & 64 & 70 & 75 & 61 & 74 \\
		\textbf{3 750 823} & \textbf{128} & 83 & 64 & 74 & 77 & 75 & 72 \\
		%\textbf{7 715 413} & \textbf{256} &  &  &  &  &  &  \\ \\
  & & \multicolumn{6}{ c }{\textbf{Two-level (5 eigenvectors)}} \\
  \textbf{DOF} & \textbf{N} & {\textbf{ORAS}} & {\textbf{SORAS}} & {\textbf{NDTNS-MRAS}} & {\textbf{NDTNS-SMRAS}} & {\textbf{TDNNS-MRAS}} & {\textbf{TDNNS-SMRAS}} \\
  \hline
		\textbf{124 109} & \textbf{4} & 15 & 32 & 17 & 35 & 24 & 37 \\
		\textbf{478 027} & \textbf{16} & 31 & 41 & 31 & 47 & 42 & 47 \\
		\textbf{933 087} & \textbf{32} & 40 & 48 & 38 & 52 & 53 & 51 \\
		\textbf{1 899 125 } & \textbf{64} & 49 & 51 & 45 & 53 & 64 & 56 \\
		\textbf{3 750 823} & \textbf{128} & 69 & 54 & 49 & 54 & 70 & 53 \\
		%\textbf{7 715 413} & \textbf{256} &  &  &  &  &  &  \\ \\
  & & \multicolumn{6}{ c }{\textbf{Two-level (7 eigenvectors)}} \\
  \textbf{DOF} & \textbf{N} & {\textbf{ORAS}} & {\textbf{SORAS}} & {\textbf{NDTNS-MRAS}} & {\textbf{NDTNS-SMRAS}} & {\textbf{TDNNS-MRAS}} & {\textbf{TDNNS-SMRAS}} \\
  \hline
		\textbf{124 109} & \textbf{4} & 14 & 33 & 16 & 30 & 24 & 35 \\
		\textbf{478 027} & \textbf{16} & 26 & 41 & 25 & 38 & 42 & 44 \\
		\textbf{933 087} & \textbf{32} & 31 & 43 & 25 & 42 & 49 & 46 \\
		\textbf{1 899 125 } & \textbf{64} & 39 & 47 & 30 & 39 & 59 & 50 \\
		\textbf{3 750 823} & \textbf{128} & 58 & 49 & 30 & 43 & 61 & 50 \\
		%\textbf{7 715 413} & \textbf{256} &  &  &  &  &  &  \\ \\
 \end{tabular}
\end{adjustbox}
\end{center}
\end{table}

\begin{figure}[!ht]
\centering
    \subfloat[Robin interface conditions]{
      \includegraphics[width=0.3\textwidth]{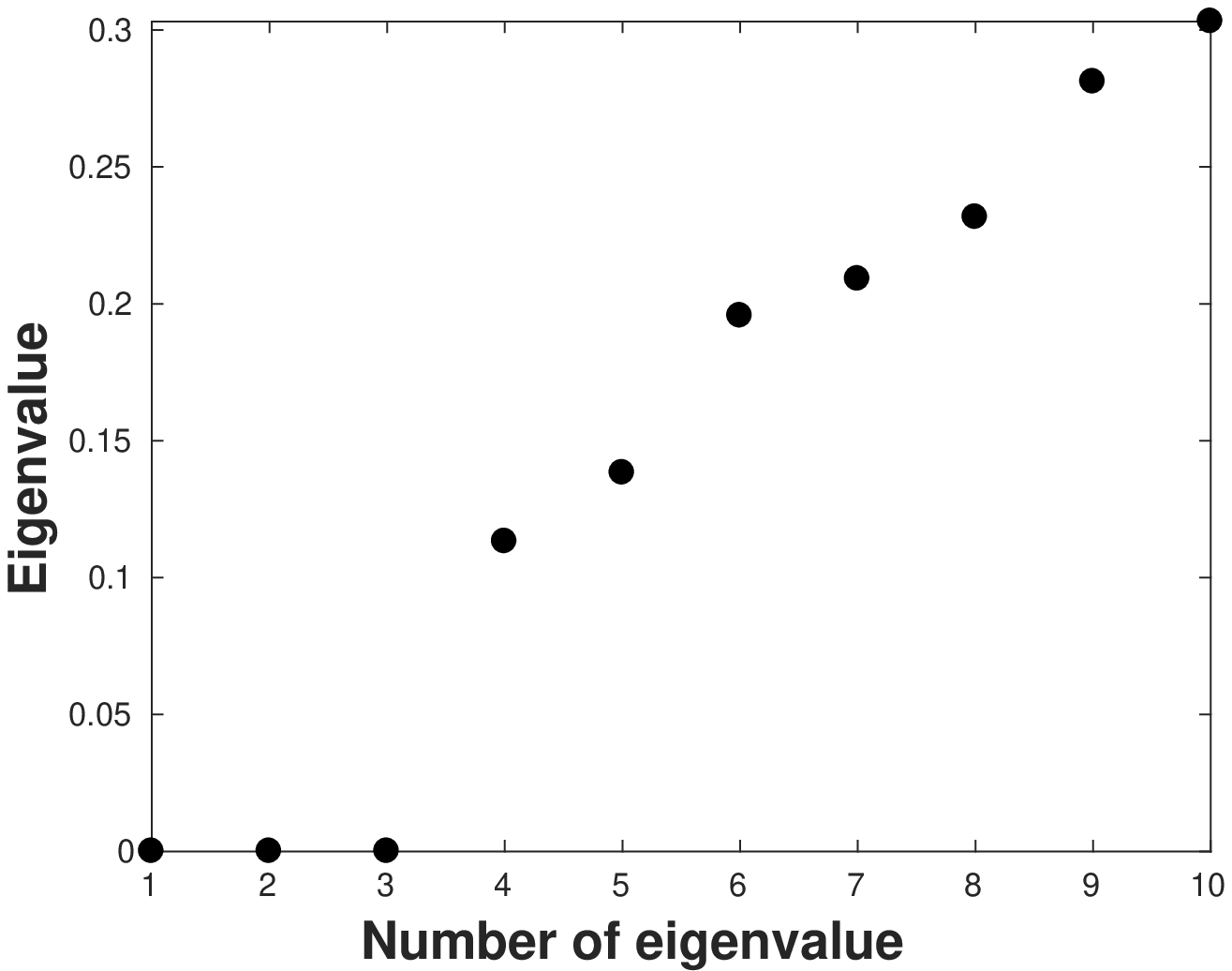}
    }
     \subfloat[NDTNS interface conditions]{
      \includegraphics[width=0.3\textwidth]{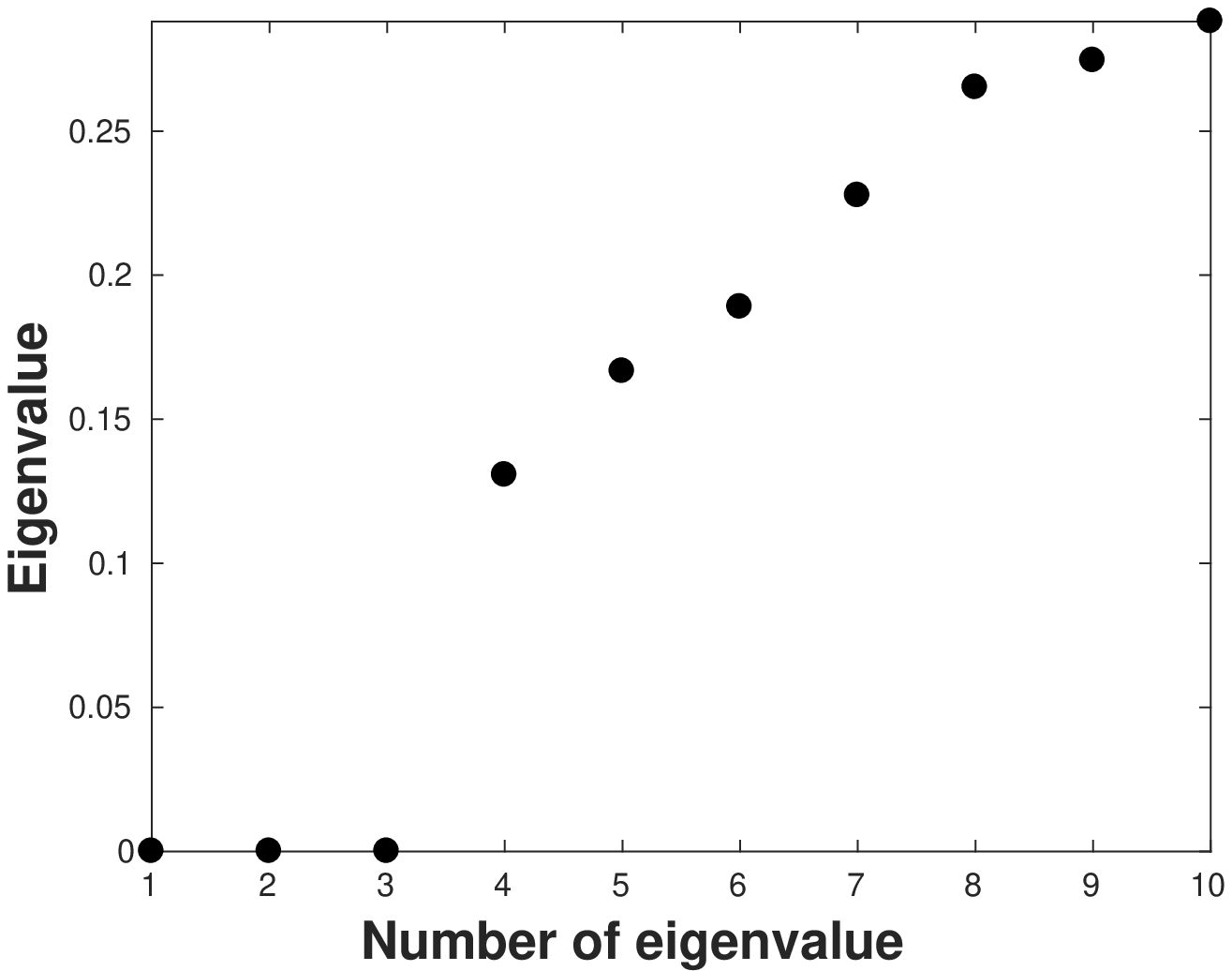}
    }
     \subfloat[TDNNS interface conditions]{
      \includegraphics[width=0.3\textwidth]{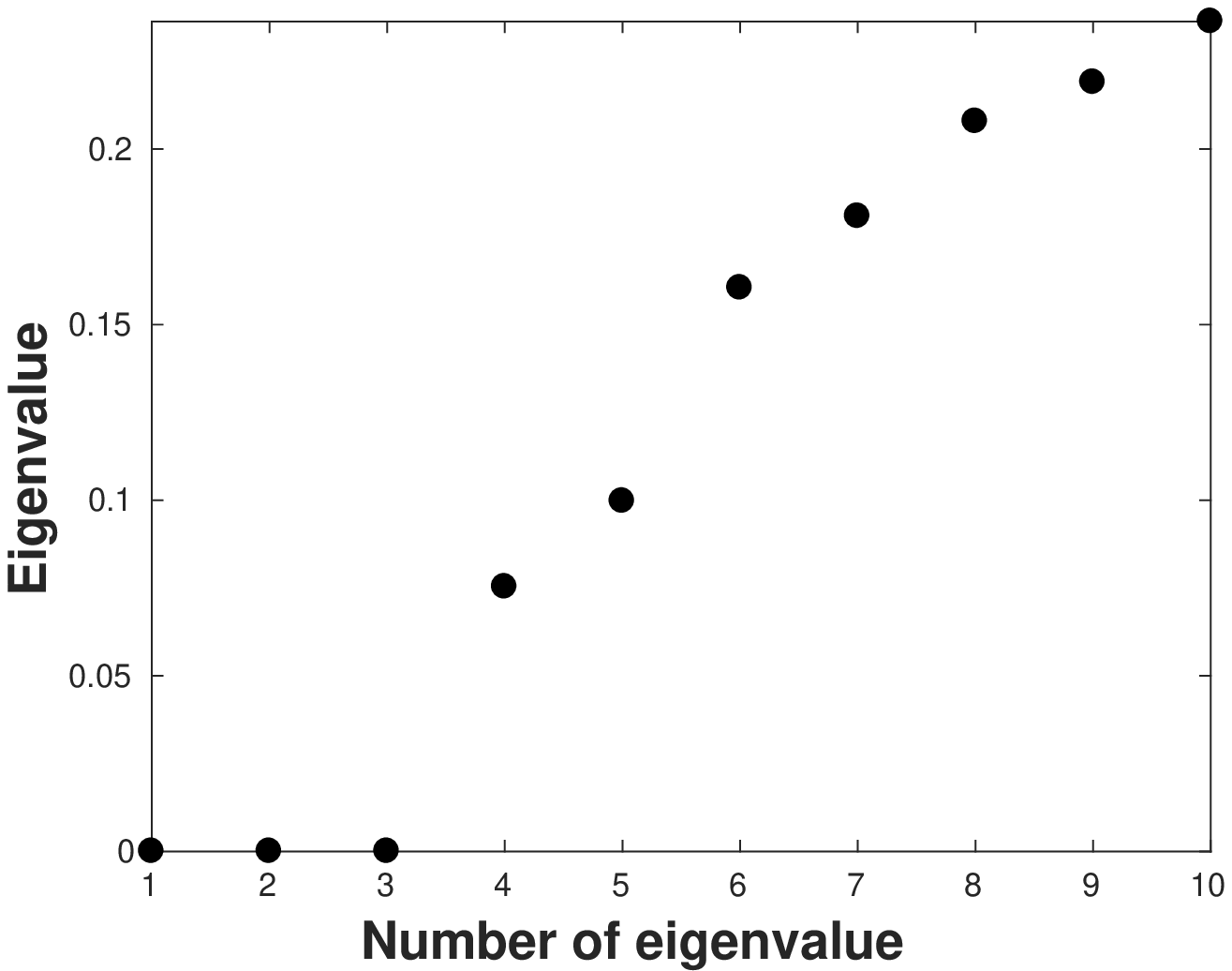}
    }
 	\caption{Eigenvalues on one of the floating subdomains in case of uniform decomposition and Taylor-Hood discretisation ($\boldsymbol{TH^3_h}, R^{2}_h$)  - the L-shaped domain problem.}
	\label{fig:TH_Cavity_elasticity_eigenvalues}
  \end{figure}
  
The results of Table~\ref{tab:TH_L_shape_elasticity} show a clear improvement in the scalability of the two-level preconditioners over the one-level ones. In fact, using five eigenvectors per subdomain, the number of iterations is virtually unaffected by the number of subdomains. All two-level preconditioners show a comparable performance. For this case, increasing the dimension of the coarse space beyond $5 \times N$ eigenvectors does not seem to improve the results dramatically.
}
\end{testcase}
% \newpage
\begin{testcase}[The heterogeneous beam problem]
\label{exp:beam}
{\rm We consider a heterogeneous beam with ten layers of steel and rubber. Five layers are made from steel with the physical parameters $E = 210 \cdot 10^9$ and $\nu = 0.3$, and other five are made from rubber with the physical parameters $E = 10^8$ and $\nu = 0.4999$ as it is depicted in Figure~\ref{fig:beam_layers}. A similar example was considered in~\cite{haferssas:hal-01278347}.
\begin{figure}[!ht]
\centering
    \subfloat[Steel and rubber layers\label{fig:beam_layers}]{\resizebox {0.45\columnwidth} {!}{
	  \begin{tikzpicture}
	  \draw[draw=white] (0,-0.2) rectangle (5,-0.1);
	  \draw[draw=white] (0,-0.1) rectangle (5,0);
	  \filldraw[fill=blue, draw=blue] (0,0) rectangle (5,0.1);
	  \filldraw[fill=red, draw=red] (0,0.1) rectangle (5,0.2);
	  \filldraw[fill=blue, draw=blue] (0,0.2) rectangle (5,0.3);
	  \filldraw[fill=red, draw=red] (0,0.3) rectangle (5,0.4);
	  \filldraw[fill=blue, draw=blue] (0,0.4) rectangle (5,0.5);
	  \filldraw[fill=red, draw=red] (0,0.5) rectangle (5,0.6);
	  \filldraw[fill=blue, draw=blue] (0,0.6) rectangle (5,0.7);
	  \filldraw[fill=red, draw=red] (0,0.7) rectangle (5,0.8);
	  \filldraw[fill=blue, draw=blue] (0,0.8) rectangle (5,0.9);
	  \filldraw[fill=red, draw=red] (0,0.9) rectangle (5,1);
    \end{tikzpicture}
    }}
    \hfill
    \subfloat[Discrete solution\label{fig:beam_solution}]{
      \includegraphics[width=0.45\textwidth]{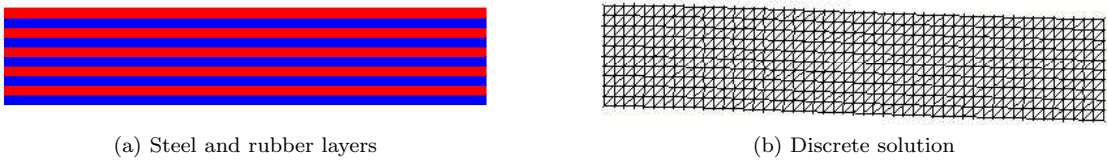}
}
\caption{Heterogeneous beam.}
  \end{figure}
The computational domain is the rectangle $\Omega = (0,5) \times (0,1)$. The beam is clamped on its left side, hence we consider the following problem
\begin{equation}
\label{eq:Heterogeneous_elasticity}
\left\{
\begin{array}{rclclr}
 -2 \mu \nabla \cdot \boldsymbol{\varepsilon}(\boldsymbol{u}) & + & \nabla p & = & (0,-1)^T& \mbox{in } \Omega \\
 & - &  \nabla \cdot \boldsymbol{u} & = & \frac{1}{\lambda} p & \mbox{in } \Omega \\
  	& & \boldsymbol{u}(x,y)  & = & (0,0)^T & \mbox{on } \partial \Omega \cap \{x=0\} \\
	& &  \boldsymbol{\sigma}^{sym}_{\boldsymbol{n}}(x,y)  & = & (0,0)^T & \mbox{on } \partial \Omega \setminus \{x=0\}
 \end{array}
\right. .
\end{equation} 
In Figure~\ref{fig:beam_solution} we plot the mesh of the bent beam. %The same as before we solved the problem numerically and multiplied by huge constant to expose the bend.
% Table~\ref{tab:TH_beam} shows the number of iterations needed to achieve a relative $L^2$ norm of the error (with respect to the one domain
% solution) smaller than $10^{-6}$.
Because of the heterogeneity of the problem, we do not notice a clear clustering of the eigenvalues (see Figure~\ref{fig:TH_beam_eigenvalues}). In such case it is well known that coarse space including only three zero energy modes is not sufficient~\cite{MR3033238}. That is why we consider a coarse space built using 5 or 7 eigenvectors per subdomain.
\begin{table}[!ht]
\caption{Comparison of preconditioners for Taylor-Hood discretisation ($\boldsymbol{TH^3_h}, R^{2}_h$) - the heterogeneous beam.}
\label{tab:TH_beam}
\begin{center}
\begin{adjustbox}{max width=\textwidth}
\begin{tabular}{ c c || c | c | c | c | c | c }
  & & \multicolumn{6}{ c }{\textbf{One-level}} \\
  \textbf{DOF} & \textbf{N} & {\textbf{ORAS}} & {\textbf{SORAS}} & {\textbf{NDTNS-MRAS}} & {\textbf{NDTNS-SMRAS}} & {\textbf{TDNNS-MRAS}} & {\textbf{TDNNS-SMRAS}} \\
%   & & Unif & MTS & Unif & MTS & Unif & MTS & Unif & MTS & Unif & MTS & Unif & MTS \\
  \hline
		\textbf{44 963} & \textbf{8} & 168 & 301 & 160 & 267 & 177 & 264 \\
		\textbf{87 587} & \textbf{16} & 226 & 490 & 245 & 462 & 229 & 424 \\
		\textbf{177 923} & \textbf{32} & 373 & 711 & 447 & 684 & 440 & 672 \\
		\textbf{347 651} & \textbf{64} & 615 & $>$1000 & 728 & $>$1000 & 746 & $>$1000 \\
		\textbf{707 843} & \textbf{128} & 973 & $>$1000 & $>$1000 & $>$1000 & $>$1000 & $>$1000 \\
		\textbf{1 385 219} & \textbf{256} & $>$1000 & $>$1000 & $>$1000 & $>$1000 & $>$1000 & $>$1000 \\ 
%   & & \multicolumn{6}{ c }{\textbf{Two-level (3 eigenvectors)}} \\
%   \textbf{DOF} & \textbf{N} & {\textbf{ORAS}} & {\textbf{SORAS}} & {\textbf{NDTNS-MRAS}} & {\textbf{NDTNS-SMRAS}} & {\textbf{TDNNS-MRAS}} & {\textbf{TDNNS-SMRAS}} \\
% %   & & Unif & MTS & Unif & MTS & Unif & MTS & Unif & MTS & Unif & MTS & Unif & MTS \\
%   \hline
% 		\textbf{44 963} & \textbf{8} & 151 & 216 & 146 & 192 & 166 & 181 \\
% 		\textbf{87 587} & \textbf{16} & 191 & 342 & 218 & 319 & 202 & 268 \\
% 		\textbf{177 923} & \textbf{32} & 306 & 467 & 371 & 444 & 386 & 449 \\
% 		\textbf{347 651} & \textbf{64} & 437 & 616 & 531 & 598 & 615 & 611 \\
% 		\textbf{707 843} & \textbf{128} & 604 & 765 & 677 & 760 & 860 & 804 \\
% 		\textbf{1 385 219} & \textbf{256} & 756 & 782 & 689 & 831 & $>$1000 & 805 \\ \\
  & & \multicolumn{6}{ c }{\textbf{Two-level (5 eigenvectors)}} \\
  \textbf{DOF} & \textbf{N} & {\textbf{ORAS}} & {\textbf{SORAS}} & {\textbf{NDTNS-MRAS}} & {\textbf{NDTNS-SMRAS}} & {\textbf{TDNNS-MRAS}} & {\textbf{TDNNS-SMRAS}} \\
%   & & Unif & MTS & Unif & MTS & Unif & MTS & Unif & MTS & Unif & MTS & Unif & MTS \\
  \hline
		\textbf{44 963} & \textbf{8} & 109 & 160 & 136 & 147 & 148 & 136 \\
		\textbf{87 587} & \textbf{16} & 136 & 204 & 192 & 200 & 181 & 184 \\
		\textbf{177 923} & \textbf{32} & 193 & 291 & 296 & 275 & 326 & 276 \\
		\textbf{347 651} & \textbf{64} & 260 & 304 & 363 & 282 & 491 & 299 \\
		\textbf{707 843} & \textbf{128} & 412 & 356 & 420 & 369 & 601 & 346 \\
		\textbf{1 385 219} & \textbf{256} & 379 & 414 & 448 & 400 & 711 & 317 \\ 
  & & \multicolumn{6}{ c }{\textbf{Two-level (7 eigenvectors)}} \\
  \textbf{DOF} & \textbf{N} & {\textbf{ORAS}} & {\textbf{SORAS}} & {\textbf{NDTNS-MRAS}} & {\textbf{NDTNS-SMRAS}} & {\textbf{TDNNS-MRAS}} & {\textbf{TDNNS-SMRAS}} \\
%   & & Unif & MTS & Unif & MTS & Unif & MTS & Unif & MTS & Unif & MTS & Unif & MTS \\
  \hline
		\textbf{44 963} & \textbf{8} & 76 & 118 & 124 & 115 & 133 & 103 \\
		\textbf{87 587} & \textbf{16} & 106 & 146 & 166 & 138 & 159 & 123 \\
		\textbf{177 923} & \textbf{32} & 157 & 202 & 203 & 185 & 302 & 214 \\
		\textbf{347 651} & \textbf{64} & 178 & 191 & 225 & 170 & 326 & 182 \\
		\textbf{707 843} & \textbf{128} & 140 & 114 & 153 & 112 & 266 & 122 \\
		\textbf{1 385 219} & \textbf{256} & 119 & 86 & 118 & 77 & 259 & 94 \\ 
 \end{tabular}
\end{adjustbox}
\end{center}
\end{table}
}
\end{testcase}
% \newpage

As in the previous example, the introduction of a coarse space provides a significant improvement in the number of iterations needed for convergence. Due to the high heterogeneity of this problem, more eigenvectors per subdomain are needed to obtain scalable results. We notice an important improvement of the convergence when using two-level methods (see Table~\ref{tab:TH_beam}). Although we get a stable number of iterations only when considering a coarse space which is sufficiently big.
% \newpage
\subsubsection{Stokes equation}
\label{sec:coarse_Stokes_TH}
% We start with the Taylor-Hood discretisation
% \begin{align*}
% 	\boldsymbol{TH^3_h} & = \left\{\boldsymbol{v_h} \in [H^1(\Omega)]^2: \quad \boldsymbol{v_h}|_K \in [\mathbb{P}_3(K)]^2 \quad \forall K \in \mathcal{T}_h\right\}, \\
% 	R^2_{h,0} & := \left\{q_h \in L^2(\Omega): \quad q_h \in \mathbb{P}_2(\Omega) \land \ \int_{\Omega} q_h \dx = 0\right\}.
% \end{align*}
We now turn to the Stokes discrete problem given in Sections~\ref{sec:TH}. Once again in case of ORAS and SORAS we choose $\alpha = 10$ as in~\cite{haferssas:hal-01278347} for the Robin interface conditions~\eqref{eq:coarse_elasticity_Robin_BC}. In the first case we consider a random initial guess for the GMRES iterative solver. Later with the second example we use zero as an initial guess.

\begin{figure}[!ht]
\centering
    \subfloat[Velocity field]{%      \includegraphics[width=0.44\textwidth]{./Cavity.eps}
      \includegraphics[width=0.49\textwidth]{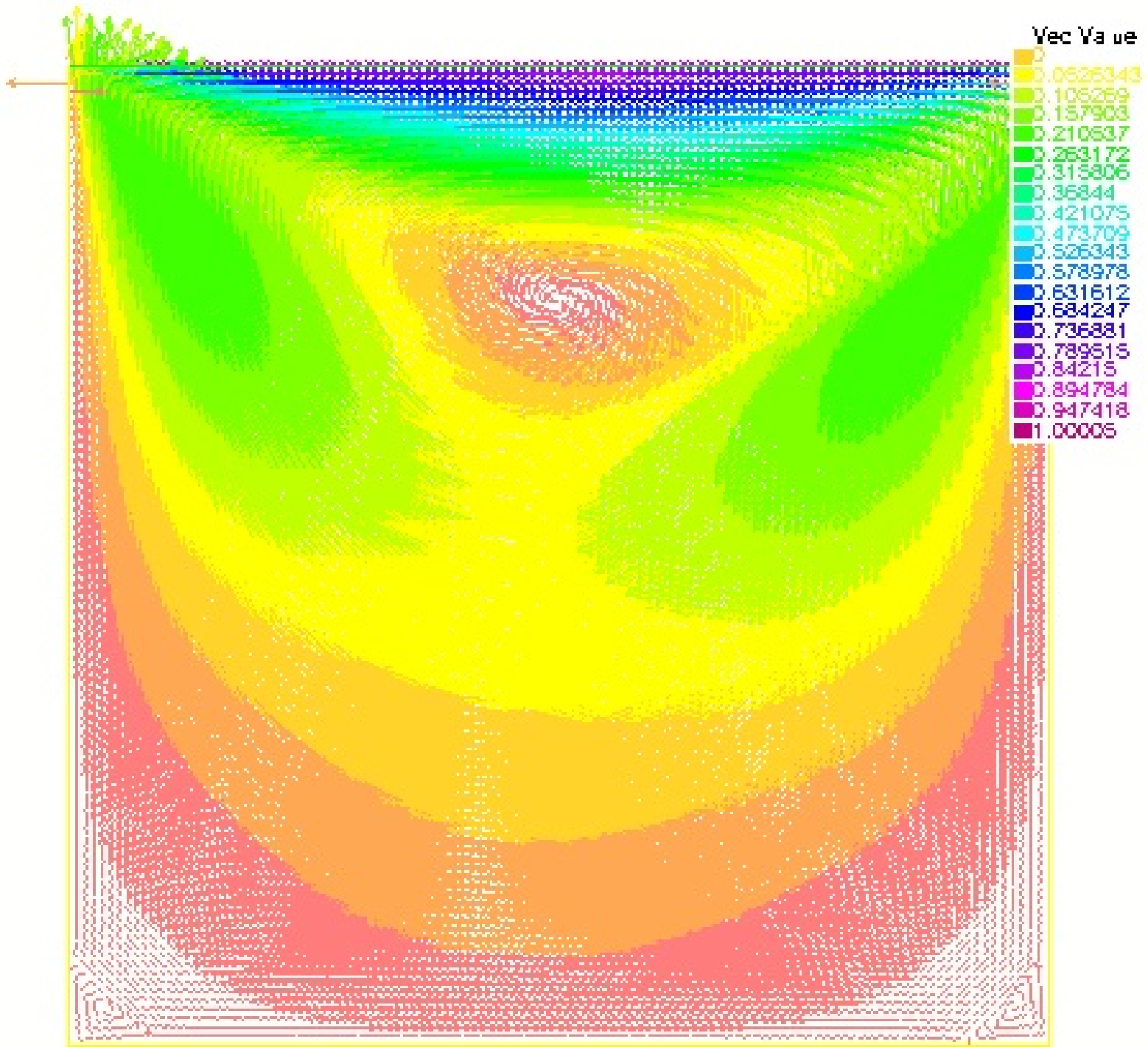}
    }
    \hfill
    \subfloat[Pressure]{%      \includegraphics[width=0.37\textwidth]{./Poisseuille_T_shape.eps}
      \includegraphics[width=0.49\textwidth]{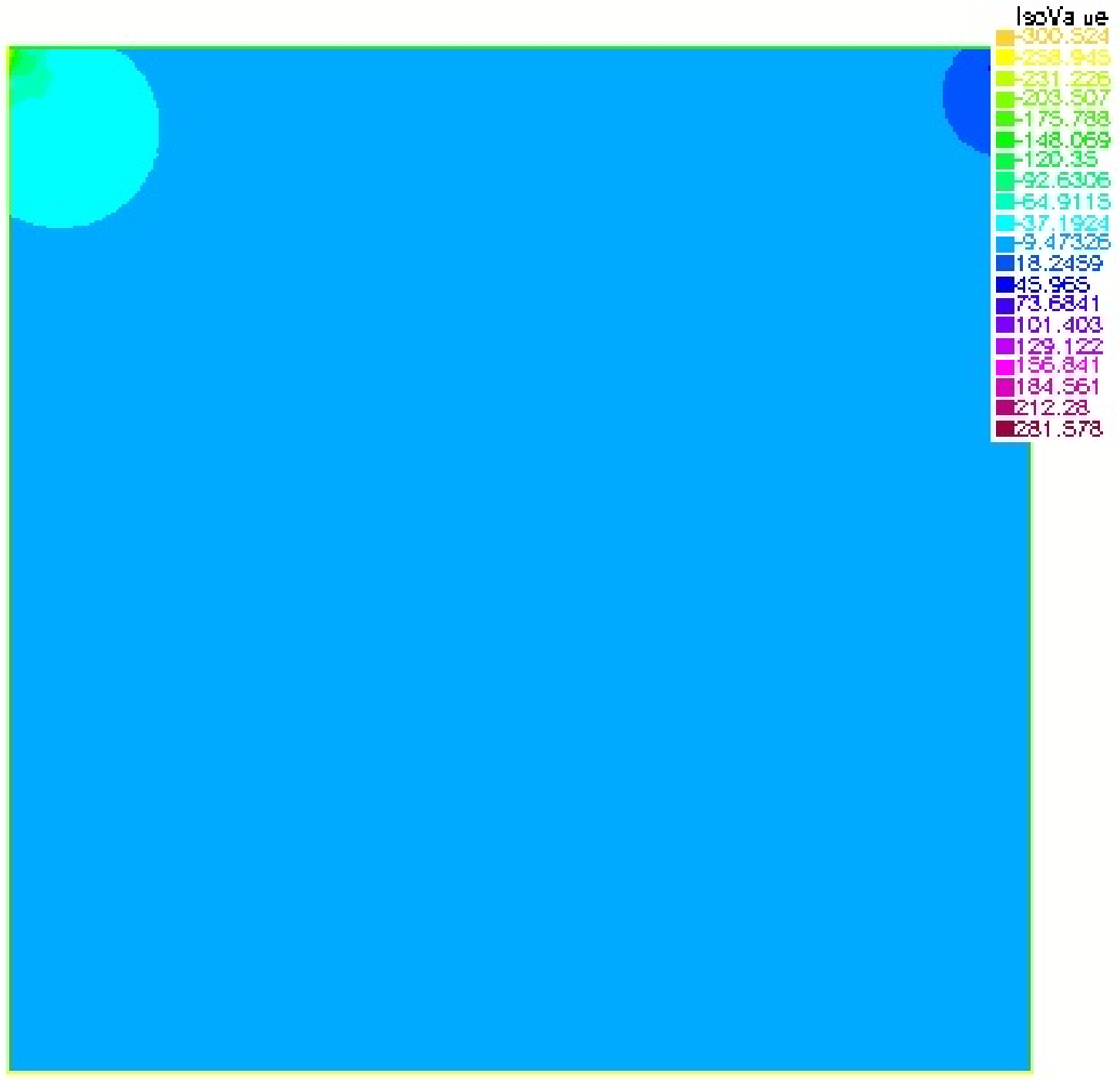}
    }
 	\caption{Numerical solution of the driven cavity problem - the driven cavity problem.}
	\label{fig:Cavity}
  \end{figure}

\begin{testcase}[The driven cavity problem]
\label{exp:cavity}
{\rm The test case is the driven cavity. We consider the following problem on the unit square $\Omega = (0,1)^2$
\begin{equation}
\label{eq:Cavity}
\left\{
\begin{array}{rclclr}
 - \Delta \boldsymbol{u} & + & \nabla p & = & \boldsymbol{f} & \mbox{in } \Omega \\
 & -&  \nabla \cdot \boldsymbol{u} & = & 0 & \mbox{in } \Omega \\
  	& & \boldsymbol{u}(x,y) & = &(1,0)^T & \mbox{on } \partial \Omega \cap \{y=1\} \\
	& &  \boldsymbol{u}(x,y)& = & (0,0)^T & \mbox{on } \partial \Omega \setminus \{y=1\}
 \end{array}
\right. .
\end{equation}
In Figure~\ref{fig:Cavity} we plot the vector field and pressure, after solving numerically the problem.
  
We start with the two energy modes only (see Figure~\ref{fig:TH_Cavity_Stokes_eigenvalues}). This already provides some improvement. Then, we add more eigenvectors to see if they bring improvement. %Unlike for linear elasticity problem, in the case of the Stokes equations we can observe only two zero eigenvalues. These two zeros are associated with two constants as eigenvectors.
  
\begin{figure}[!ht]
\centering
    \subfloat[Robin interface conditions]{
      \includegraphics[width=0.3\textwidth]{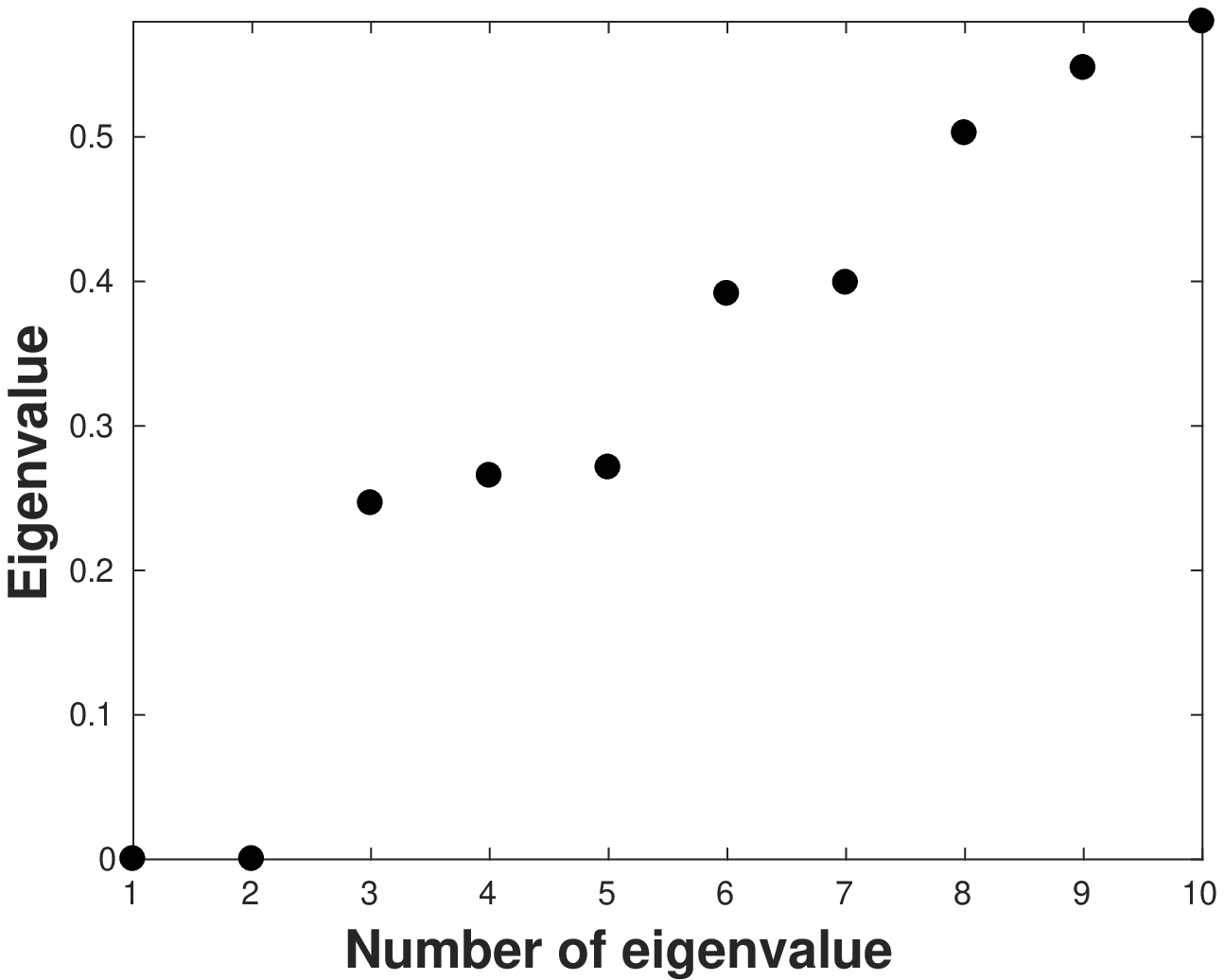}
    }
     \subfloat[NVTF interface conditions]{
      \includegraphics[width=0.3\textwidth]{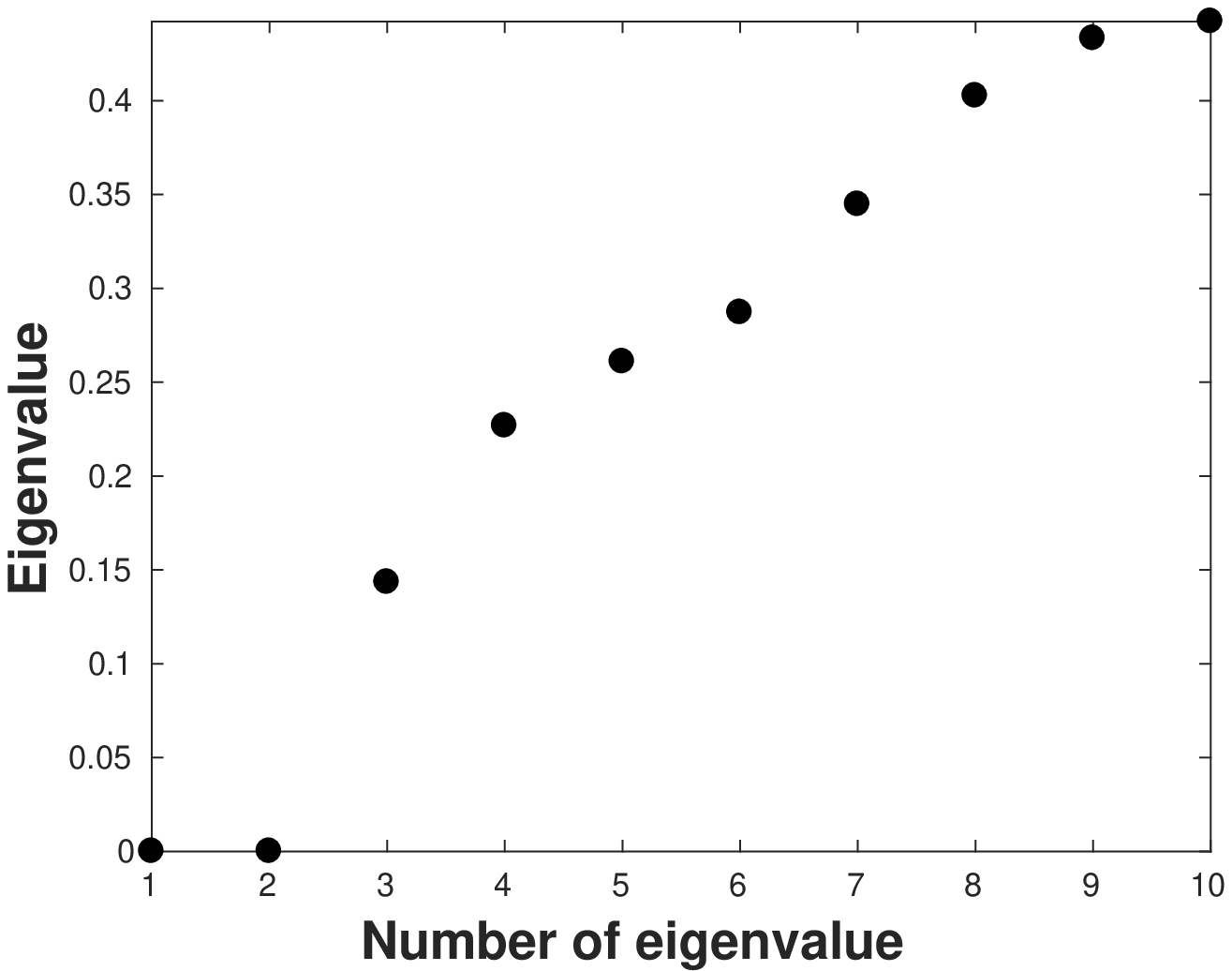}
    }
     \subfloat[TVNF interface conditions]{
      \includegraphics[width=0.3\textwidth]{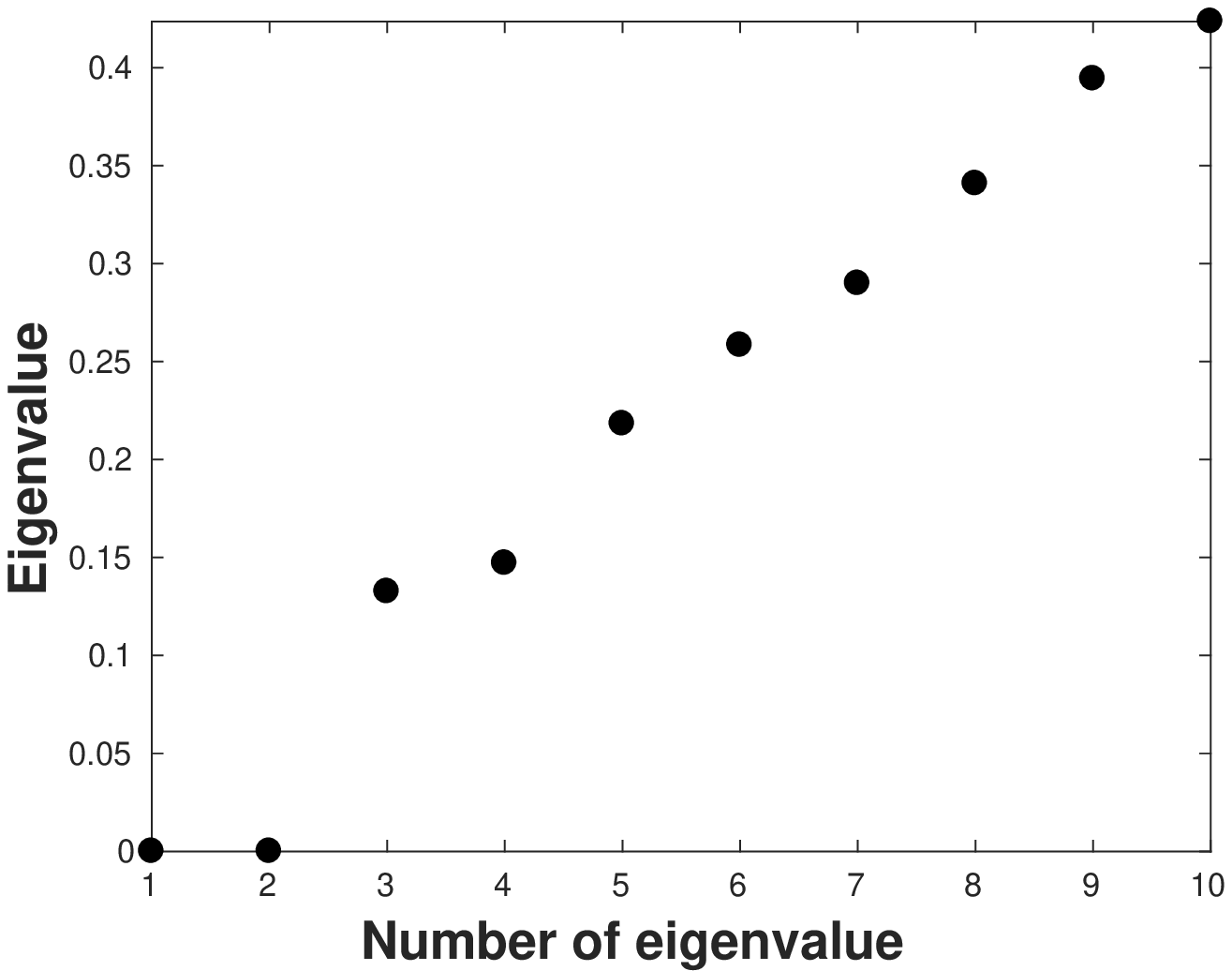}
    }
 	\caption{Eigenvalues on one of the floating subdomains in case of uniform decomposition and Taylor-Hood discretisation ($\boldsymbol{TH^2_h}, R^{1}_h$) - the driven cavity problem.}
	\label{fig:TH_Cavity_Stokes_eigenvalues}
  \end{figure}
% \newpage

\begin{table}[!ht]
\caption{Comparison of preconditioners for Taylor-Hood discretisation ($\boldsymbol{TH^2_h}, R^{1}_h$) - the driven cavity problem.}
\label{tab:lowest_TH_Cavity_Stokes}
\begin{center}
\begin{adjustbox}{max width=\textwidth}
\begin{tabular}{ c c || c c | c c | c c | c c | c c | c c }
%   & \multicolumn{4}{ c }{\textbf{Uniform decomposition (236 $\times$ 236)}} \\ 
%   \hline
  & & \multicolumn{12}{ c }{\textbf{One-level}} \\
  \textbf{DOF} & \textbf{N} & \multicolumn{2}{ c |}{\textbf{ORAS}} & \multicolumn{2}{ c |}{\textbf{SORAS}} & \multicolumn{2}{ c |}{\textbf{NVTF-MRAS}} & \multicolumn{2}{ c |}{\textbf{NVTF-SMRAS}} & \multicolumn{2}{ c |}{\textbf{TVNF-MRAS}} & \multicolumn{2}{ c }{\textbf{TVNF-SMRAS}} \\
  & & Unif & MTS & Unif & MTS & Unif & MTS & Unif & MTS & Unif & MTS & Unif & MTS \\
  \hline
		\textbf{91 003} & \textbf{4} & 12 & 17 & 24 & 34 & 22 & 22 & 34 & 40 & 22 & 25 & 30 & 40 \\
		\textbf{362 003} & \textbf{16} & 28 & 35 & 56 & 67 & 52 & 53 & 90 & 106 & 54 & 53 & 70 & 84 \\
		\textbf{813 003} & \textbf{36} & 39 & 75 & 92 & 103 & 85 & 91 & 165 & 185 & 91 & 88 & 118 & 136 \\
		\textbf{1 444 003} & \textbf{64} & 53 & 91 & 120 & 144 & 120 & 135 & 254 & 283 & 132 & 132 & 169 & 206 \\
		\textbf{2 728 003} & \textbf{121} & 80 & 278 & 180 & 212 & 182 & 280 & 412 & 580 & 199 & 213 & 251 & 439 \\
		\textbf{5 768 003} & \textbf{256} & $>$1000 & $>$1000 & 271 & 317 & 303 & 452 & 917 & 955 & 322 & 319 & 397 & 695 \\ 
  & & \multicolumn{12}{ c }{\textbf{Two-level (2 eigenvectors)}} \\
  \textbf{DOF} & \textbf{N} & \multicolumn{2}{ c |}{\textbf{ORAS}} & \multicolumn{2}{ c |}{\textbf{SORAS}} & \multicolumn{2}{ c |}{\textbf{NVTF-MRAS}} & \multicolumn{2}{ c |}{\textbf{NVTF-SMRAS}} & \multicolumn{2}{ c |}{\textbf{TVNF-MRAS}} & \multicolumn{2}{ c }{\textbf{TVNF-SMRAS}} \\
  & & Unif & MTS & Unif & MTS & Unif & MTS & Unif & MTS & Unif & MTS & Unif & MTS \\
  \hline
		\textbf{91 003} & \textbf{4} & 10 & 14 & 18 & 22 & 19 & 17 & 26 & 30 & 27 & 20 & 21 & 26 \\
		\textbf{362 003} & \textbf{16} & 20 & 25 & 32 & 37 & 33 & 34 & 50 & 62 & 60 & 40 & 42 & 51 \\
		\textbf{813 003} & \textbf{36} & 27 & 33 & 36 & 44 & 47 & 49 & 62 & 86 & 79 & 53 & 59 & 63 \\
		\textbf{1 444 003} & \textbf{64} & 31 & 42 & 38 & 53 & 104 & 66 & 85 & 114 & 85 & 52 & 62 & 79 \\
		\textbf{2 728 003} & \textbf{121} & 39 & 103 & 39 & 51 & 74 & 81 & 85 & 133 & 92 & 86 & 62 & 93 \\
		\textbf{5 768 003} & \textbf{256} & 300 & 849 & 46 & 54 & 109 & 108 & 146 & 132 & 91 & 78 & 63 & 90 \\ 
  & & \multicolumn{12}{ c }{\textbf{Two-level (5 eigenvectors)}} \\
  \textbf{DOF} & \textbf{N} & \multicolumn{2}{ c |}{\textbf{ORAS}} & \multicolumn{2}{ c |}{\textbf{SORAS}} & \multicolumn{2}{ c |}{\textbf{NVTF-MRAS}} & \multicolumn{2}{ c |}{\textbf{NVTF-SMRAS}} & \multicolumn{2}{ c |}{\textbf{TVNF-MRAS}} & \multicolumn{2}{ c }{\textbf{TVNF-SMRAS}} \\
  & & Unif & MTS & Unif & MTS & Unif & MTS & Unif & MTS & Unif & MTS & Unif & MTS \\
  \hline
		\textbf{91 003} & \textbf{4} & 9 & 12 & 13 & 16 & 16 & 15 & 18 & 20 & 25 & 20 & 16 & 18 \\
		\textbf{362 003} & \textbf{16} & 16 & 20 & 21 & 24 & 27 & 22 & 28 & 37 & 56 & 37 & 26 & 35 \\
		\textbf{813 003} & \textbf{36} & 23 & 27 & 25 & 26 & 33 & 30 & 39 & 40 & 65 & 41 & 28 & 37 \\
		\textbf{1 444 003} & \textbf{64} & 26 & 36 & 27 & 29 & 40 & 34 & 35 & 45 & 77 & 45 & 28 & 42 \\
		\textbf{2 728 003} & \textbf{121} & 35 & 41 & 29 & 32 & 43 & 38 & 34 & 48 & 84 & 72 & 29 & 47 \\
		\textbf{5 768 003} & \textbf{256} & 66 & 60 & 32 & 33 & 56 & 41 & 60 & 49 & 88 & 61 & 29 & 44 \\ 
%   & & \multicolumn{12}{ c }{\textbf{Two-level (7 eigenvectors)}} \\
%   \textbf{DOF} & \textbf{N} & \multicolumn{2}{ c |}{\textbf{ORAS}} & \multicolumn{2}{ c |}{\textbf{SORAS}} & \multicolumn{2}{ c |}{\textbf{NVTF-MRAS}} & \multicolumn{2}{ c |}{\textbf{NVTF-SMRAS}} & \multicolumn{2}{ c |}{\textbf{TVNF-MRAS}} & \multicolumn{2}{ c }{\textbf{TVNF-SMRAS}} \\
%   & & Unif & MTS & Unif & MTS & Unif & MTS & Unif & MTS & Unif & MTS & Unif & MTS \\
%   \hline
% 		\textbf{91 003} & \textbf{4} & 8 & 10 & 12 & 15 & 13 & 14 & 15 & 17 & 25 & 18 & 14 & 17 \\
% 		\textbf{362 003} & \textbf{16} & 15 & 19 & 19 & 21 & 25 & 22 & 23 & 30 & 61 & 39 & 22 & 28 \\
% 		\textbf{813 003} & \textbf{36} & 19 & 25 & 21 & 22 & 31 & 26 & 25 & 31 & 69 & 44 & 23 & 27 \\
% 		\textbf{1 444 003} & \textbf{64} & 23 & 30 & 22 & 25 & 36 & 30 & 27 & 34 & 83 & 43 & 24 & 31 \\
% 		\textbf{3 246 003} & \textbf{144} & 31 & 40 & 22 & 28 &  & 32 &  & 36 & 92 & 60 & 25 & 33 \\
% 		\textbf{5 768 003} & \textbf{256} &  &  &  &  &  &  &  &  &  &  &  &  \\ \\
 \end{tabular}
\end{adjustbox}
\end{center}
\end{table}

The conclusions remain the same as for the L-shaped domain problem for the nearly incompressible elasticity equation discretised by Taylor-Hood method ($\boldsymbol{TH^3_h}, R^{2}_h$) since Tables %~\ref{tab:TH_Cavity_Stokes},
 \ref{tab:lowest_TH_Cavity_Stokes} and~\ref{tab:TH_L_shape_elasticity} %, \ref{tab:lowest_TH_L_shape_elasticity}
 show similar results.
}
\end{testcase}
%\newpage
\begin{testcase}[The T-shaped domain problem]
\label{exp:T_shape}
{\rm Finally, we consider a T-shaped domain $\Omega = (0,1.5) \times (0,1) \cup (0.5,1) \times (-1,1)$, and we impose mixed boundary conditions given by
\begin{equation}
\label{eq:T_shape_BC}
 \boldsymbol{u}(x,y) = \left\{
  \begin{array}{c l}
  	(4y(1-y),0)^T & \mbox{if } x=0 \mbox{ or } x=1.5 \\
	(0,0)^T & \mbox{otherwise}.
  \end{array}
%   \begin{array}{c l}
%   	\boldsymbol{u}(x,y) = (4y(1-y),0)^T & \mbox{if } x=0 
% 	\sigma_{nn}(x,y) = 0, u_t(x,y)=0 & \mbox{if } x=1.5 \\
% 	\boldsymbol{u}(x,y) = (0,0)^T & \mbox{otherwise}.
%   \end{array}
\right.
\end{equation}
\begin{figure}[!ht]
\centering
    \subfloat[Velocity field $\boldsymbol{u}$]{%
      \includegraphics[width=0.49\textwidth]{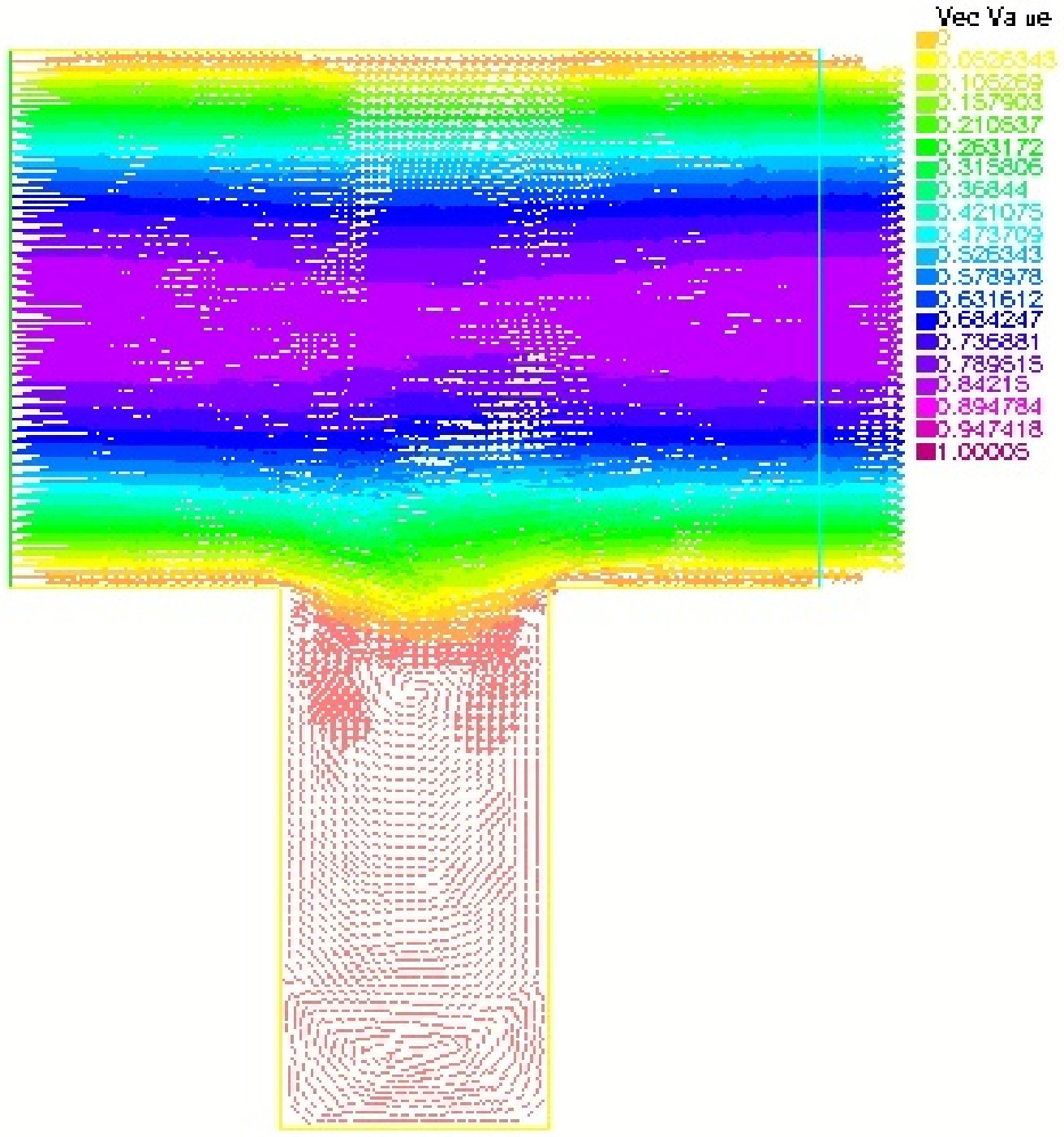}
    }\hfill
    \subfloat[Pressure $p$]{%
      \includegraphics[width=0.49\textwidth]{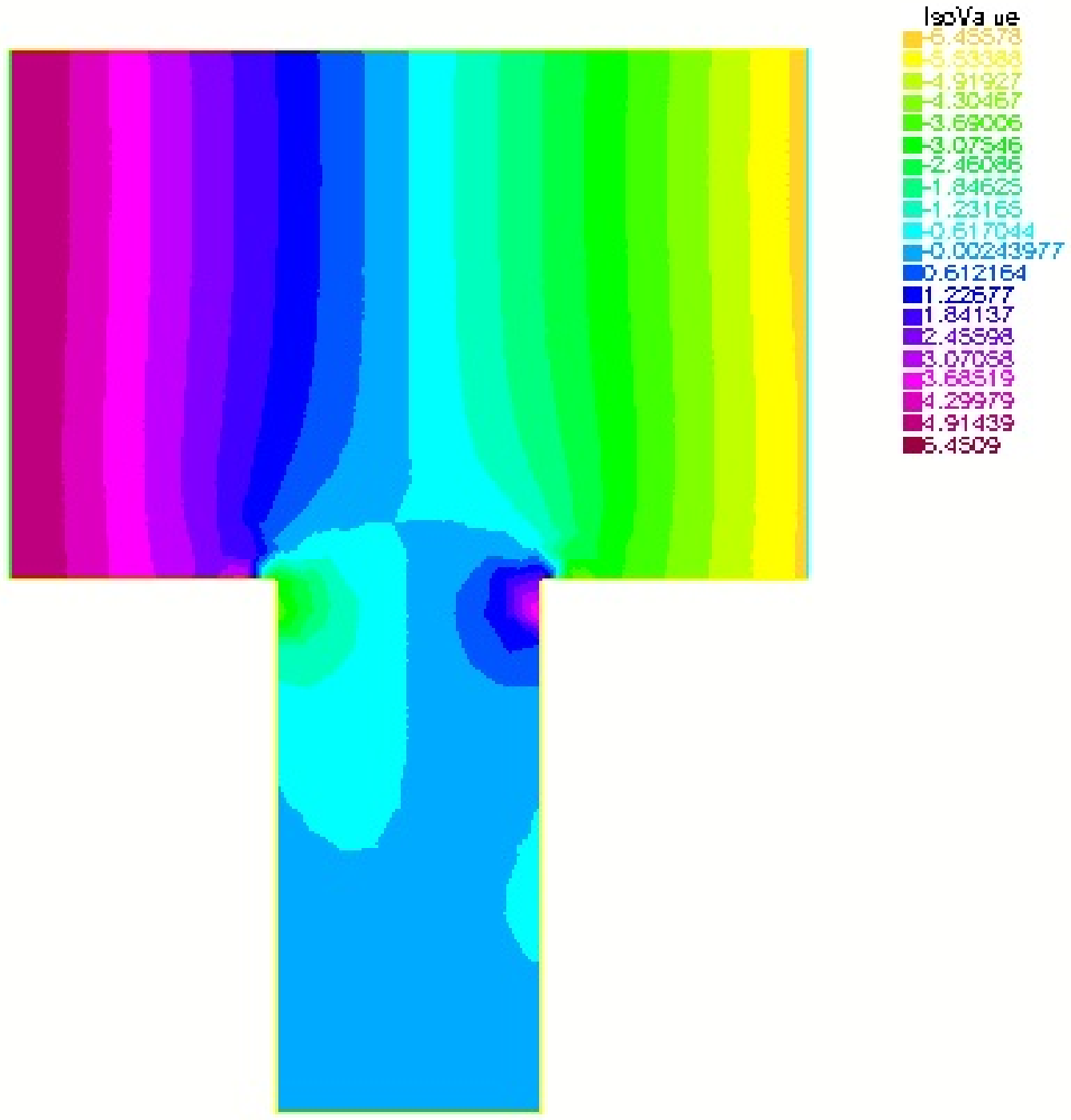}
    }
  \caption{Numerical solution - the T-shaped problem.}
  \label{fig:T_shape}
  \end{figure}
The numerical solution of this problem is depicted in Figure~\ref{fig:T_shape}. The overlapping decomposition into subdomains is generated by METIS.

% Table~\ref{tab:TH_T_shape_Stokes} shows the number of iterations needed to
% achieve a relative $L^2$ norm of the error (with respect to the one domain
% solution) smaller than $10^{-6}$. 
Once again a clustering of small eigenvalues of generalised eigenvalue problem defined in~\eqref{eq:coarse_second_GenEO} is a motivation of the size of the coarse space (see Figure~\ref{fig:TH_T_shape_Stokes_eigenvalues}).
\begin{table}[!ht]
\caption{Comparison of preconditioners for Taylor-Hood discretisation ($\boldsymbol{TH^3_h}, R^{2}_h$) - the T-shaped problem.}
\label{tab:TH_T_shape_Stokes}
\begin{center}
\begin{adjustbox}{max width=\textwidth}
\begin{tabular}{ c c | c c c c c c }
%   & \multicolumn{4}{ c }{\textbf{Uniformorm decomposition (236 $\times$ 236)}} \\ 
%   \hline
  & & \multicolumn{6}{ c }{\textbf{One-level}} \\
  \textbf{DOF} & \textbf{N} & {\textbf{ORAS}} & {\textbf{SORAS}} & {\textbf{NVTF-MRAS}} & {\textbf{NVTF-SMRAS}} & {\textbf{TVNF-MRAS}} & {\textbf{TVNF-SMRAS}} \\
  \hline
		\textbf{33 269} & \textbf{4} & 13 & 20 & 12 & 19 & 13 & 19 \\
		\textbf{138 316} & \textbf{16} & 36 & 51 & 33 & 52 & 31 & 45 \\
		\textbf{269 567} & \textbf{32} & 59 & 85 & 52 & 85 & 49 & 75 \\
		\textbf{553 103} & \textbf{64} & 92 & 132 & 83 & 136 & 78 & 115 \\
		\textbf{1 134 314} & \textbf{128} & 146 & 208 & 132 & 223 & 117 & 188 \\
		\textbf{2 201 908} & \textbf{256}  & 232 & 328 & 209 & 357 & 189 & 293 \\
  & & \multicolumn{6}{ c }{\textbf{Two-level (2 eigenvectors)}} \\
  \textbf{DOF} & \textbf{N} & {\textbf{ORAS}} & {\textbf{SORAS}} & {\textbf{NVTF-MRAS}} & {\textbf{NVTF-SMRAS}} & {\textbf{TVNF-MRAS}} & {\textbf{TVNF-SMRAS}} \\
  \hline
		\textbf{33 269} & \textbf{4} & 10 & 14 & 9 & 15 & 12 & 15 \\
		\textbf{138 316} & \textbf{16} & 21 & 27 & 19 & 24 & 22 & 24 \\
		\textbf{269 567} & \textbf{32} & 29 & 35 & 30 & 38 & 25 & 30 \\
		\textbf{553 103} & \textbf{64} & 35 & 45 & 34 & 43 & 33 & 35 \\
		\textbf{1 134 314} & \textbf{128} & 42 & 52 & 47 & 58 & 34 & 41 \\
		\textbf{2 201 908} & \textbf{256}  & 47 & 56 & 69 & 76 & 38 & 45 \\
  & & \multicolumn{6}{ c }{\textbf{Two-level (5 eigenvectors)}} \\
  \textbf{DOF} & \textbf{N} & {\textbf{ORAS}} & {\textbf{SORAS}} & {\textbf{NVTF-MRAS}} & {\textbf{NVTF-SMRAS}} & {\textbf{TVNF-MRAS}} & {\textbf{TVNF-SMRAS}} \\
  \hline
		\textbf{33 269} & \textbf{4} & 8 & 13 & 8 & 13 & 12 & 14 \\
		\textbf{138 316} & \textbf{16} & 15 & 16 & 14 & 16 & 20 & 18 \\
		\textbf{269 567} & \textbf{32} & 14 & 19 & 20 & 22 & 24 & 19 \\
		\textbf{553 103} & \textbf{64} & 16 & 20 & 18 & 19 & 29 & 20 \\
		\textbf{1 134 314} & \textbf{128} & 17 & 22 & 23 & 24 & 30 & 22 \\
		\textbf{2 201 908} & \textbf{256}  & 16 & 21 & 34 & 37 & 35 & 24 \\
%  & & \multicolumn{6}{ c }{\textbf{Two-level (7 eigenvectors)}} \\
%  \textbf{DOF} & \textbf{N} & {\textbf{ORAS}} & {\textbf{SORAS}} & {\textbf{NVTF-MRAS}} & {\textbf{NVTF-SMRAS}} & {\textbf{TVNF-MRAS}} & {\textbf{TVNF-SMRAS}} \\
%  \hline
%		\textbf{33 269} & \textbf{4} & 8 & 13 & 8 & 12 & 12 & 15 \\
%		\textbf{138 316} & \textbf{16} & 15 & 15 & 15 & 16 & 23 & 19 \\
%		\textbf{269 567} & \textbf{32} & 18 & 16 & 14 & 16 & 23 & 19 \\
%		\textbf{553 103} & \textbf{64} & 23 & 17 & 16 & 17 & 26 & 21 \\
%		\textbf{1 134 314} & \textbf{128} & 30 & 21 & 17 & 18 & 31 & 20 \\
%		\textbf{2 201 908} & \textbf{256} & 36 & 22 & 19 & 20 & 32 & 22 \\
 \end{tabular}
\end{adjustbox}
\end{center}
\end{table}

\begin{figure}[!ht]
\centering
    \subfloat[Robin interface conditions]{
      \includegraphics[width=0.3\textwidth]{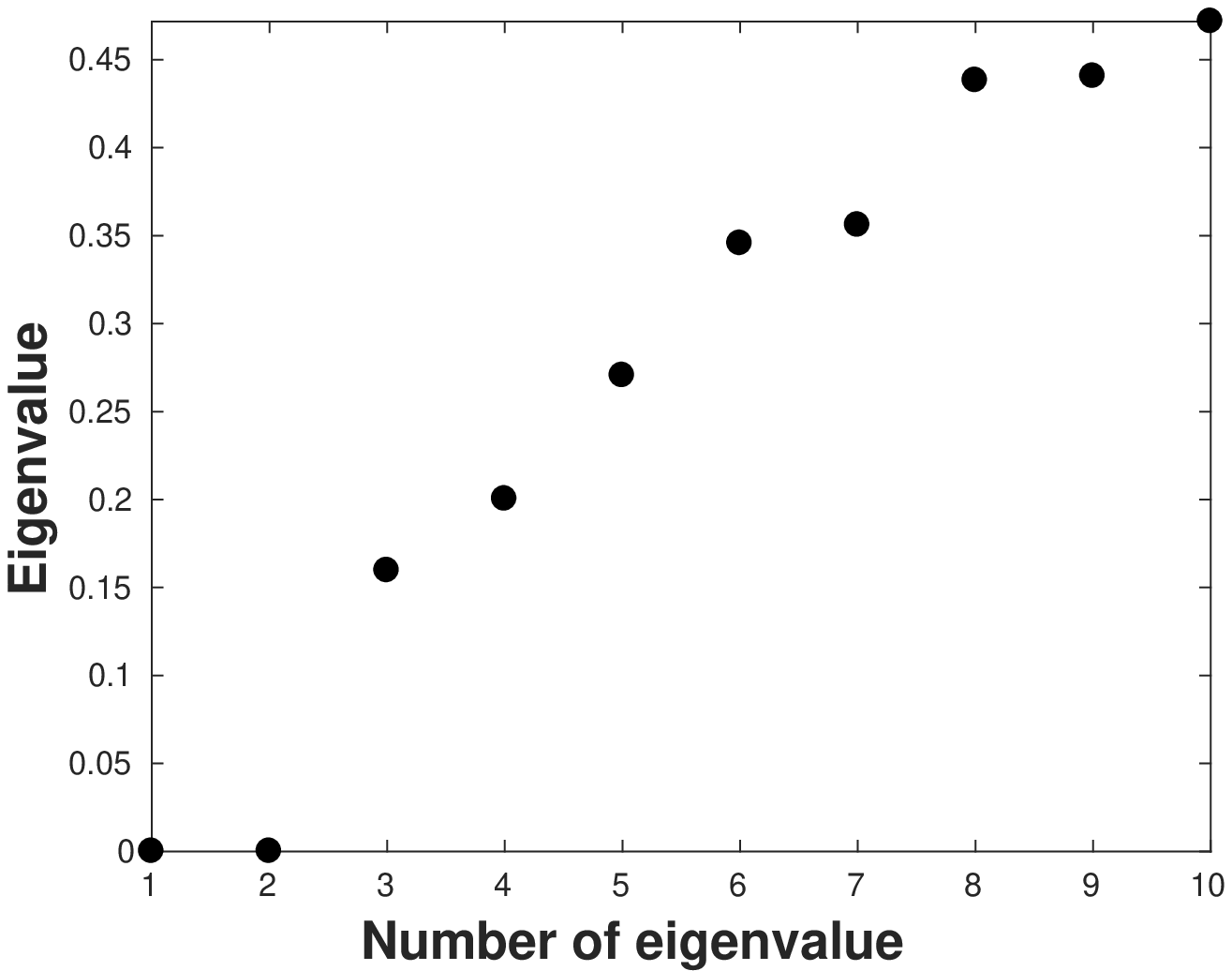}
    }
     \subfloat[NVTF interface conditions]{
      \includegraphics[width=0.3\textwidth]{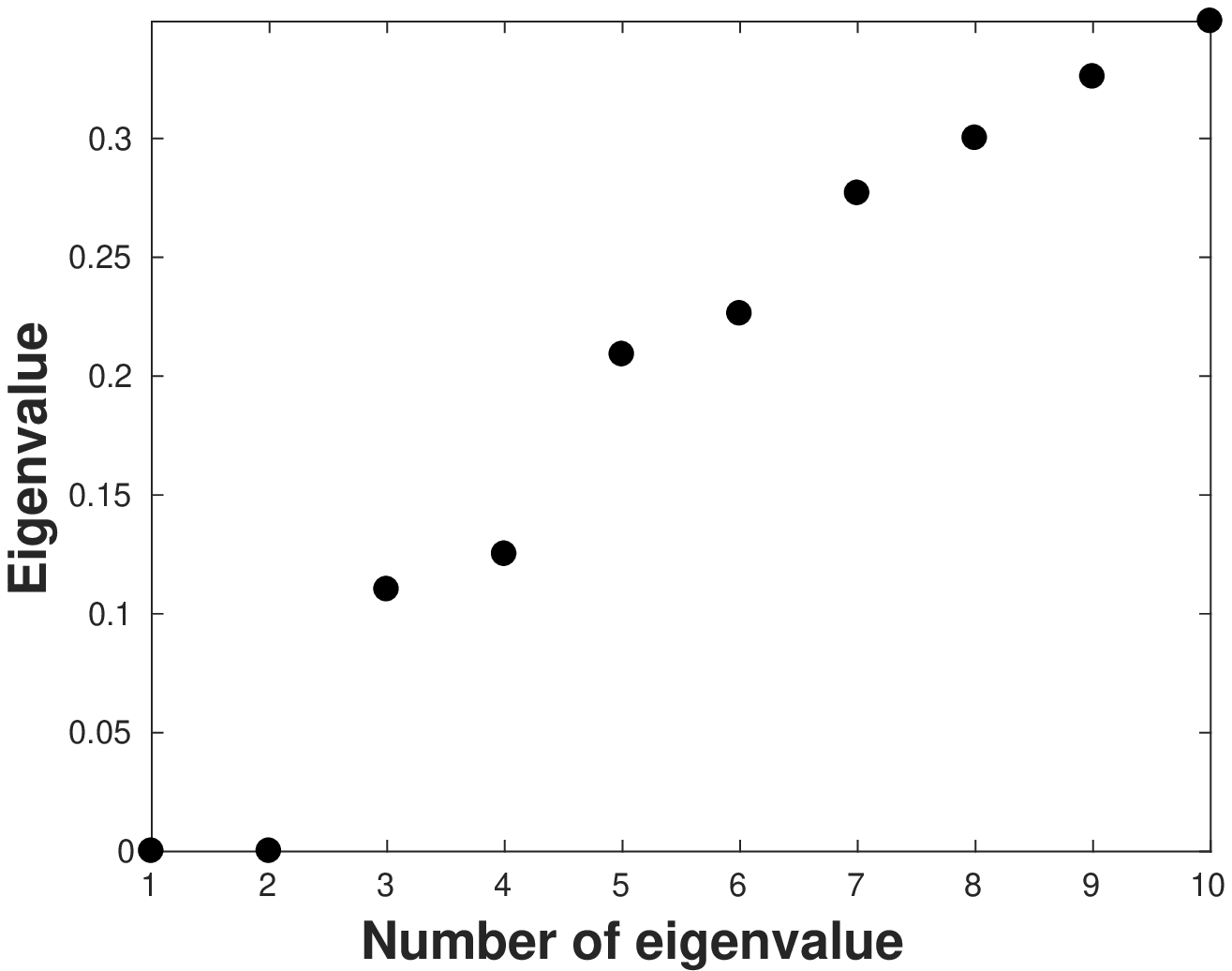}
    }
     \subfloat[TVNF interface conditions]{
      \includegraphics[width=0.3\textwidth]{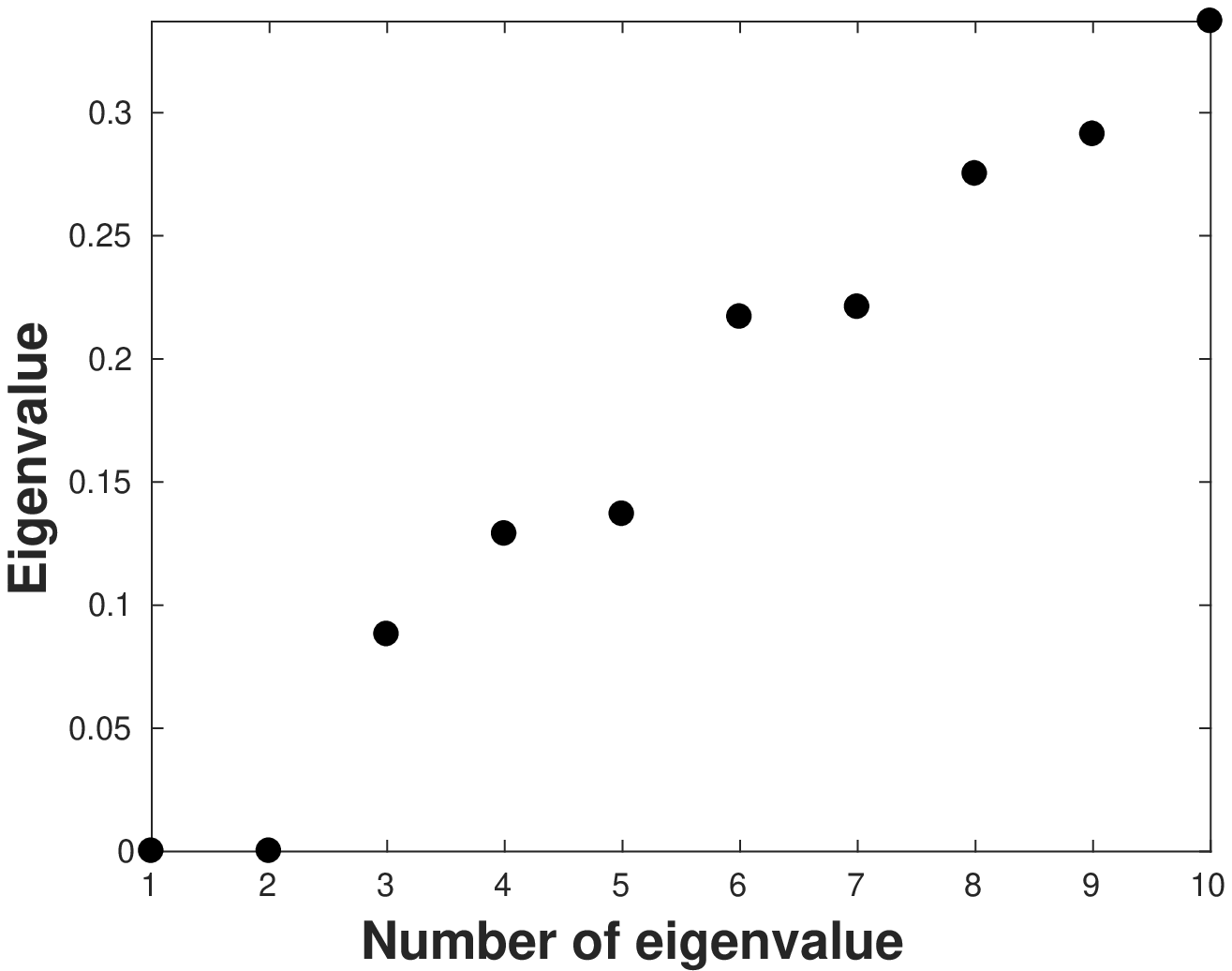}
    }
 	\caption{Eigenvalues on one of the floating subdomains in case of METIS decomposition and Taylor-Hood discretisation ($\boldsymbol{TH^3_h}, R^{2}_h$) - the T-shaped problem.}
	\label{fig:TH_T_shape_Stokes_eigenvalues}
  \end{figure}

The same as in all examples for Taylor-Hood discretisation we notice an important improvement of the convergence when using two-level methods. Although from Table~\ref{tab:TH_T_shape_Stokes} we can see that the coarse spaces containing five eigenvectors seem to be sufficient.
}
\end{testcase}
%  \newpage

% 
%%%%%%%%%%%%%%%%%%%%%%%%%%%%%%%%%%%%%%%%%%%%%%%%%%%%% hybrid discontinuous Galerkin %%%%%%%%%%%%%%%%%%%%%%%%%%%%%%%%%%%%%%%%%%%%%%%%%%%%%%%%%%%
\subsection{hdG discretisation}
\label{sec:coarse_hdg}
In this section we discretise nearly incompressible elasticity equation and the Stokes flow by using the lowest order hdG discretisation introduced in Section~\ref{sec:hdg}.  

\subsubsection{Nearly incompressible elasticity}
\label{sec:coarse_elasticity_hdg}

In case of ORAS and SORAS we consider the Robin interface conditions as in~\cite{haferssas:hal-01278347} with $\alpha = 10$. %We already know that for the Stokes equation these preconditioners are associated with NVTF and TVNF interface conditions. In the case of nearly incompressible elasticity we refer to them as normal-displacement and tangential-normal-stress (NDTNS) and tangential-displacement and normal-normal-stress (TDNNS) interface conditions. The second type of the boundary conditions has been already introduced for linear elasticity equation in~\cite{MR2826472}.
For all numerical experiments in this section we use zero as an initial guess for the GMRES iterative solver.  Moreover, the overlapping decomposition into subdomains is generated by METIS.

\begin{testcase}[The L-shaped domain problem]%{\ref{exp:L_shape}}
{\rm We consider the L-shaped domain discrete problem~\eqref{eq:L_shape_elasticity}. 

% Table~\ref{tab:hdG_Cavity_elasticity} shows the number of iterations needed to
% achieve a relative $L^2$ norm of the error (with respect to the one domain
% solution) smaller than $10^{-6}$.

% Table~\ref{tab:TH_Cavity_elasticity} shows the number of iterations needed to achieve a relative $L^2$ norm of the error (with respect to the one domain solution) smaller than $10^{-6}$.
\begin{table}[!ht]
\caption{Comparison of preconditioners for hdG discretisation - the L-shaped domain problem.}
\label{tab:hdg_L_shape_elasticity}
\begin{center}
\begin{adjustbox}{max width=\textwidth}
\begin{tabular}{ c c | c c c c c c }
  & & \multicolumn{6}{ c }{\textbf{One-level}} \\
  \textbf{DOF} & \textbf{N} & {\textbf{ORAS}} & {\textbf{SORAS}} & {\textbf{NDTNS-MRAS}} & {\textbf{NDTNS-SMRAS}} & {\textbf{TDNNS-MRAS}} & {\textbf{TDNNS-SMRAS}} \\
  \hline
		\textbf{238 692} & \textbf{8} & 61 & 158 & 64 & 174 & 77 & 177 \\
		\textbf{466 094} & \textbf{16} & 123 & 232 & 101 & 259 & 109 & 306 \\
		\textbf{948 921} & \textbf{32} & 267 & 331 & 160 & 415 & 179 & 473 \\
		\textbf{1 874 514} & \textbf{64} & 622 & 477 & 243 & 685 & 254 & 657 \\
		\textbf{3 856 425} & \textbf{128} & $>$1000 & 752 & 479 & $>$1000 & 523 & $>$1000 \\
		%\textbf{7 445 578} & \textbf{256} & $>$1000 &  &  &  &  &  \\ \\
  & & \multicolumn{6}{ c }{\textbf{Two-level (3 eigenvectors)}} \\
  \textbf{DOF} & \textbf{N} & {\textbf{ORAS}} & {\textbf{SORAS}} & {\textbf{NDTNS-MRAS}} & {\textbf{NDTNS-SMRAS}} & {\textbf{TDNNS-MRAS}} & {\textbf{TDNNS-SMRAS}} \\
  \hline
		\textbf{238 692} & \textbf{8} & 48 & 98 & 52 & 98 & 61 & 116 \\
		\textbf{466 094} & \textbf{16} & 89 & 99 & 71 & 123 & 75 & 148 \\
		\textbf{948 921} & \textbf{32} & 250 & 130 & 110 & 158 & 118 & 173 \\
		\textbf{1 874 514} & \textbf{64} & 535 & 135 & 135 & 155 & 129 & 159 \\
		\textbf{3 856 425} & \textbf{128} & $>$1000 & 152 & 172 & 176 & 181 & 192 \\
		%\textbf{7 445 578} & \textbf{256} &  &  &  &  &  &  \\ \\
  & & \multicolumn{6}{ c }{\textbf{Two-level (5 eigenvectors)}} \\
  \textbf{DOF} & \textbf{N} & {\textbf{ORAS}} & {\textbf{SORAS}} & {\textbf{NDTNS-MRAS}} & {\textbf{NDTNS-SMRAS}} & {\textbf{TDNNS-MRAS}} & {\textbf{TDNNS-SMRAS}} \\
  \hline
		\textbf{238 692} & \textbf{8} & 43 & 81 & 44 & 74 & 61 & 94 \\
		\textbf{466 094} & \textbf{16} & 77 & 82 & 51 & 92 & 63 & 103 \\
		\textbf{948 921} & \textbf{32} & 197 & 100 & 79 & 119 & 96 & 121 \\
		\textbf{1 874 514} & \textbf{64} & 429 & 103 & 102 & 122 & 110 & 138 \\
		\textbf{3 856 425} & \textbf{128} & $>$1000 & 118 & 122 & 129 & 141 & 167 \\
		%\textbf{7 445 578} & \textbf{256} &  &  &  &  &  &  \\ \\
  & & \multicolumn{6}{ c }{\textbf{Two-level (7 eigenvectors)}} \\
  \textbf{DOF} & \textbf{N} & {\textbf{ORAS}} & {\textbf{SORAS}} & {\textbf{NDTNS-MRAS}} & {\textbf{NDTNS-SMRAS}} & {\textbf{TDNNS-MRAS}} & {\textbf{TDNNS-SMRAS}} \\
  \hline
		\textbf{238 692} & \textbf{8} & 35 & 67 & 38 & 71 & 44 & 82 \\
		\textbf{466 094} & \textbf{16} & 61 & 79 & 47 & 80 & 51 & 95 \\
		\textbf{948 921} & \textbf{32} & 153 & 90 & 58 & 95 & 74 & 115 \\
		\textbf{1 874 514} & \textbf{64} & 423 & 95 & 71 & 93 & 72 & 111 \\
		\textbf{3 856 425} & \textbf{128} & 934 & 110 & 79 & 104 & 108 & 133 \\
		%\textbf{7 445 578} & \textbf{256} &  &  &  &  &  &  \\ \\
 \end{tabular}
 \end{adjustbox}
\end{center}
\end{table}
  
Table~\ref{tab:hdg_L_shape_elasticity} shows an important improvement of the convergence that is brought by the two-level methods. We cannot conclude that SMRAS preconditioners are much better than SORAS. Although we can note that coarse space improvement is visible for MRAS preconditioners and not for ORAS. For symmetric preconditioners (SMRAS and SORAS) five eigenvectors seem to lead already to satisfying results. While for the non-symmetric ones a bigger coarse space is required.  On the other hand, we state the fact that the new preconditioners are parameter-free, which makes them easier to use as no parameter is required. %We cannot conclude that MRAS or SMRAS are much better than ORAS or SORAS. Although by taking into account the fact that MRAS or SMRAS are parameter-free when imposing the interface conditions, we can say that they are easier to use. 
%Moreover, as expected, the iteration number does not increase with respect to the number of the subdomains. The size of coarse space seems to be sufficient as adding some more eigenvectors does not bring an improvement, while the computational time increases. 	
}
\end{testcase}
 \newpage
\begin{testcase}[The heterogeneous beam problem]%{\ref{exp:beam}}
{\rm We consider the heterogeneous beam with ten layers of steel and rubber that is defined as a problem~\eqref{eq:Heterogeneous_elasticity}. 
% Table~\ref{tab:TH_beam} shows the number of iterations needed to achieve a relative $L^2$ norm of the error (with respect to the one domain
% solution) smaller than $10^{-6}$.
%As it was observed when this problem was approximated using the Taylor-Hood discretisation, since the problem has strong heterogeneity, we are unable to see a clear clustering of the eigenvalues (see Figure~\ref{fig:TH_beam_eigenvalues}).
 %We use zero as an initial guess for the GMRES iterative solver. The overlapping decomposition into subdomains is generated by METIS. 

% Table~\ref{tab:hdg_beam} shows the number of iterations needed to
% achieve a relative $L^2$ norm of the error (with respect to the one domain
% solution) smaller than $10^{-6}$.
\begin{table}[!ht]
\caption{Comparison of preconditioners for hdG discretisation - the heterogeneous beam.}
\label{tab:hdg_beam}
\begin{center}
\begin{adjustbox}{max width=\textwidth}
\begin{tabular}{ c c || c | c | c | c | c | c }
  & & \multicolumn{6}{ c }{\textbf{One-level}} \\
  \textbf{DOF} & \textbf{N} & {\textbf{ORAS}} & {\textbf{SORAS}} & {\textbf{NDTNS-MRAS}} & {\textbf{NDTNS-SMRAS}} & {\textbf{TDNNS-MRAS}} & {\textbf{TDNNS-SMRAS}} \\
%   & & Unif & MTS & Unif & MTS & Unif & MTS & Unif & MTS & Unif & MTS & Unif & MTS \\
  \hline
		\textbf{46 777} & \textbf{8} & 196 & 440 & 189 & 402 & 186 & 463 \\
		\textbf{88 720} & \textbf{16} & 317 & 602 & 330 & 582 & 326 & 666 \\
		\textbf{179 721} & \textbf{32} & 537 & $>$1000 & 574 & $>$1000 & 587 & $>$1000 \\
		\textbf{353 440} & \textbf{64} & 899 & $>$1000 & 847 & $>$1000 & 846 & $>$1000 \\
		\textbf{704 329} & \textbf{128} & $>$1000 & $>$1000 & $>$1000 & $>$1000 & $>$1000 & $>$1000 \\
		\textbf{1 410 880} & \textbf{256} & $>$1000 & $>$1000 & $>$1000 & $>$1000 & $>$1000 & $>$1000 \\
%   & & \multicolumn{6}{ c }{\textbf{Two-level (3 eigenvectors)}} \\
%   \textbf{DOF} & \textbf{N} & {\textbf{ORAS}} & {\textbf{SORAS}} & {\textbf{NDTNS-MRAS}} & {\textbf{NDTNS-SMRAS}} & {\textbf{TDNNS-MRAS}} & {\textbf{TDNNS-SMRAS}} \\
% %   & & Unif & MTS & Unif & MTS & Unif & MTS & Unif & MTS & Unif & MTS & Unif & MTS \\
%   \hline
% 		\textbf{46 777} & \textbf{8} & 183 & 356 & 178 & 326 & 175 & 374 \\
% 		\textbf{88 720} & \textbf{16} & 267 & 436 & 293 & 432 & 289 & 472 \\
% 		\textbf{179 721} & \textbf{32} & 457 & 749 & 515 & 733 & 536 & 853 \\
% 		\textbf{353 440} & \textbf{64} & 705 & 870 & 720 & 918 & 718 & $>$1000 \\
% 		\textbf{704 329} & \textbf{128} & $>$1000 & $>$1000 & $>$1000 & $>$1000 & $>$1000 & $>$1000 \\
% 		\textbf{1 410 880} & \textbf{256} & $>$1000 & $>$1000 & $>$1000 & $>$1000 & $>$1000 & $>$1000 \\
  & & \multicolumn{6}{ c }{\textbf{Two-level (5 eigenvectors)}} \\
  \textbf{DOF} & \textbf{N} & {\textbf{ORAS}} & {\textbf{SORAS}} & {\textbf{NDTNS-MRAS}} & {\textbf{NDTNS-SMRAS}} & {\textbf{TDNNS-MRAS}} & {\textbf{TDNNS-SMRAS}} \\
%   & & Unif & MTS & Unif & MTS & Unif & MTS & Unif & MTS & Unif & MTS & Unif & MTS \\
  \hline
		\textbf{46 777} & \textbf{8} & 168 & 255 & 162 & 230 & 161 & 275 \\
		\textbf{88 720} & \textbf{16} & 244 & 313 & 273 & 299 & 262 & 346 \\
		\textbf{179 721} & \textbf{32} & 385 & 525 & 442 & 458 & 469 & 587 \\
		\textbf{353 440} & \textbf{64} & 514 & 444 & 551 & 526 & 590 & 558 \\
		\textbf{704 329} & \textbf{128} & 835 & 557 & 782 & 684 & 765 & 832 \\
		\textbf{1 410 880} & \textbf{256} & $>$1000 & 567 & $>$1000 & 694 & 844 & 821 \\
  & & \multicolumn{6}{ c }{\textbf{Two-level (7 eigenvectors)}} \\
  \textbf{DOF} & \textbf{N} & {\textbf{ORAS}} & {\textbf{SORAS}} & {\textbf{NDTNS-MRAS}} & {\textbf{NDTNS-SMRAS}} & {\textbf{TDNNS-MRAS}} & {\textbf{TDNNS-SMRAS}} \\
%   & & Unif & MTS & Unif & MTS & Unif & MTS & Unif & MTS & Unif & MTS & Unif & MTS \\
  \hline
		\textbf{46 777} & \textbf{8} & 148 & 197 & 149 & 192 & 158 & 231 \\
		\textbf{88 720} & \textbf{16} & 205 & 201 & 286 & 187 & 283 & 273 \\
		\textbf{179 721} & \textbf{32} & 318 & 337 & 385 & 301 & 433 & 419 \\
		\textbf{353 440} & \textbf{64} & 403 & 262 & 397 & 247 & 460 & 389 \\
		\textbf{704 329} & \textbf{128} & 490 & 168 & 447 & 182 & 558 & 443 \\
		\textbf{1 410 880} & \textbf{256} & $>$1000 & 116 & 387 & 138 & 473 & 298 \\
 \end{tabular}
\end{adjustbox}
\end{center}
\end{table}
}
\end{testcase}

We notice an improvement only when using a coarse space which is sufficiently big (see Table~\ref{tab:hdg_beam}). Furthermore, we get a stable number of iterations only for the symmetric preconditioners (SMRAS and SORAS), and the coarse space improvement in case of ORAS preconditioner is much less visible than in case of MRAS preconditioners. This may be due to the fact we have not chosen an optimized parameter in the Robin interface conditions~\eqref{eq:coarse_elasticity_Robin_BC}.
% \newpage

%%%%%%%%%%%%%%%%%%%%%%%%%%%%%%%%%%%%%%%%%%%%%%%%%%%%% hybrid discontinuous Galerkin %%%%%%%%%%%%%%%%%%%%%%%%%%%%%%%%%%%%%%%%%%%%%%%%%%%%%%%%%%%
\subsubsection{Stokes equation}
\label{sec:coarse_Stokes_hdg}
We now turn to the Stokes discrete problem given in~\ref{sec:hdg}. Once again in case of ORAS and SORAS we choose $\alpha = 10$ as in~\cite{haferssas:hal-01278347} for the Robin interface conditions~\eqref{eq:coarse_elasticity_Robin_BC}.  In the first case we consider a random initial guess for the GMRES iterative solver. Later with the second example we use zero as an initial guess.

\begin{testcase}[The driven cavity problem]%{\ref{exp:cavity}}
{\rm We consider the driven cavity defined as a problem~\eqref{eq:Cavity}. 
% Table~\ref{tab:TH_Cavity_Stokes} shows the number of iterations needed to
% achieve a relative $L^2$ norm of the error (with respect to the one domain
% solution) smaller than $10^{-6}$.
%In this case  the size of the coarse space is motivated by observing a clustering of small eigenvalues of the generalised eigenvalue problem defined in~\eqref{eq:coarse_second_GenEO} (see Figure~\ref{fig:hdg_Cavity_Stokes_eigenvalues}).  On the plots as expected we can two zero eigenvalues corresponding to the two zero energy modes.

% Table~\ref{tab:hdg_Cavity_Stokes} shows the number of iterations needed to
% achieve a relative $L^2$ norm of the error (with respect to the one domain
% solution) smaller than $10^{-6}$.
\begin{table}[!ht]
\caption{Comparison of preconditioners for hdG discretisation - the driven cavity problem.}
\label{tab:hdg_Cavity_Stokes}
\begin{center}
\begin{adjustbox}{max width=\textwidth}
\begin{tabular}{ c c || c c | c c | c c | c c | c c | c c }
%   & \multicolumn{4}{ c }{\textbf{Uniform decomposition (236 $\times$ 236)}} \\ 
%   \hline
  & & \multicolumn{12}{ c }{\textbf{One-level}} \\
  \textbf{DOF} & \textbf{N} & \multicolumn{2}{ c |}{\textbf{ORAS}} & \multicolumn{2}{ c |}{\textbf{SORAS}} & \multicolumn{2}{ c |}{\textbf{NVTF-MRAS}} & \multicolumn{2}{ c |}{\textbf{NVTF-SMRAS}} & \multicolumn{2}{ c |}{\textbf{TVNF-MRAS}} & \multicolumn{2}{ c }{\textbf{TVNF-SMRAS}} \\
  & & Unif & MTS & Unif & MTS & Unif & MTS & Unif & MTS & Unif & MTS & Unif & MTS \\
  \hline
		\textbf{93 656} & \textbf{4} & 17 & 18 & 37 & 38 & 24 & 22 & 44 & 44 & 32 & 25 & 48 & 50 \\
		\textbf{373 520} & \textbf{16} & 76 & 122 & 75 & 84 & 52 & 54 & 107 & 111 & 68 & 67 & 122 & 126 \\
		\textbf{839 592} & \textbf{36} & 152 & 327 & 120 & 133 & 91 & 96 & 194 & 200 & 112 & 115 & 206 & 210 \\
		\textbf{1 491 872} & \textbf{64} & 261 & 587 & 162 & 176 & 130 & 143 & 294 & 303 & 159 & 158 & 292 & 286 \\
		\textbf{2 819 432} & \textbf{121} & 364 & $>$1000 & 229 & 256 & 199 & 213 & 504 & 649 & 238 & 251 & 628 & 643 \\
		\textbf{5 963 072} & \textbf{256} & 592 & $>$1000 & 367 & 398 & 326 & 477 & $>$1000 & $>$1000 & 392 & 404 & 995 & 740 \\
  & & \multicolumn{12}{ c }{\textbf{Two-level (2 eigenvectors)}} \\
  \textbf{DOF} & \textbf{N} & \multicolumn{2}{ c |}{\textbf{ORAS}} & \multicolumn{2}{ c |}{\textbf{SORAS}} & \multicolumn{2}{ c |}{\textbf{NVTF-MRAS}} & \multicolumn{2}{ c |}{\textbf{NVTF-SMRAS}} & \multicolumn{2}{ c |}{\textbf{TVNF-MRAS}} & \multicolumn{2}{ c }{\textbf{TVNF-SMRAS}} \\
  & & Unif & MTS & Unif & MTS & Unif & MTS & Unif & MTS & Unif & MTS & Unif & MTS \\
  \hline
		\textbf{93 656} & \textbf{4} & 12 & 14 & 30 & 28 & 18 & 18 & 33 & 32 & 40 & 23 & 38 & 37 \\
		\textbf{373 520} & \textbf{16} & 81 & 80 & 47 & 57 & 36 & 40 & 61 & 73 & 100 & 49 & 85 & 82 \\
		\textbf{839 592} & \textbf{36} & 236 & 228 & 61 & 60 & 57 & 65 & 97 & 104 & 132 & 66 & 112 & 107 \\
		\textbf{1 491 872} & \textbf{64} & 395 & 463 & 67 & 71 & 79 & 85 & 139 & 129 & 142 & 70 & 128 & 122 \\
		\textbf{2 819 432} & \textbf{121} & 840 & $>$1000 & 73 & 86 & 113 & 127 & 188 & 178 & 157 & 86 & 127 & 139 \\
		\textbf{5 963 072} & \textbf{256} & $>$1000 & $>$1000 & 80 & 87 & 171 & 179 & 283 & 287 & 167 & 108 & 132 & 148 \\ 
  & & \multicolumn{12}{ c }{\textbf{Two-level (5 eigenvectors)}} \\
  \textbf{DOF} & \textbf{N} & \multicolumn{2}{ c |}{\textbf{ORAS}} & \multicolumn{2}{ c |}{\textbf{SORAS}} & \multicolumn{2}{ c |}{\textbf{NVTF-MRAS}} & \multicolumn{2}{ c |}{\textbf{NVTF-SMRAS}} & \multicolumn{2}{ c |}{\textbf{TVNF-MRAS}} & \multicolumn{2}{ c }{\textbf{TVNF-SMRAS}} \\
  & & Unif & MTS & Unif & MTS & Unif & MTS & Unif & MTS & Unif & MTS & Unif & MTS \\
  \hline
		\textbf{93 656} & \textbf{4} & 10 & 12 & 25 & 24 & 14 & 16 & 23 & 22 & 52 & 22 & 29 & 26 \\
		\textbf{373 520} & \textbf{16} & 27 & 35 & 38 & 37 & 27 & 29 & 38 & 41 & 117 & 39 & 53 & 53 \\
		\textbf{839 592} & \textbf{36} & 135 & 84 & 45 & 41 & 35 & 37 & 51 & 50 & 145 & 49 & 64 & 61 \\
		\textbf{1 491 872} & \textbf{64} & 278 & 212 & 49 & 45 & 44 & 42 & 58 & 55 & 157 & 59 & 64 & 64 \\
		\textbf{2 819 432} & \textbf{121} & 607 & 584 & 56 & 49 & 46 & 56 & 58 & 62 & 162 & 81 & 65 & 75 \\
		\textbf{5 963 072} & \textbf{256} & $>$1000 & $>$1000 & 62 & 55 & 52 & 64 & 57 & 69 & 166 & 75 & 65 & 75 \\ 
 \end{tabular}
\end{adjustbox}
\end{center}
\end{table}
The conclusions remain the same as in the case of nearly incompressible elasticity equation for the L-shaped domain. Although Table~\ref{tab:hdg_Cavity_Stokes} shows that coarse spaces containing five eigenvectors seem to decrease the number of iteration even in the case of MRAS preconditioners that are not fully scalable.
}
\end{testcase}
% \newpage
\begin{testcase}[The T-shaped domain problem]%{\ref{exp:T_shape}}
{\rm Finally, we consider a T-shaped domain $\Omega = (0,1.5) \times (0,1) \cup (0.5,1) \times (-1,1)$, and we impose mixed boundary conditions~\eqref{eq:T_shape_BC}. The numerical solution of this problem is depicted in Figure~\ref{fig:T_shape}.
% Table~\ref{tab:hdg_T_shape_Stokes} shows the number of iterations needed to
% achieve a relative $L^2$ norm of the error (with respect to the one domain
% solution) smaller than $10^{-6}$.
%In this case we motivate the size of the coarse space by observing a clustering of small eigenvalues of generalised eigenvalue problem defined in~\eqref{eq:coarse_second_GenEO} (see Figure~\ref{fig:hdg_T_shape_Stokes_eigenvalues}).
\begin{table}[!ht]
\caption{Comparison of preconditioners for hdG discretisation - the T-shaped problem.}
\label{tab:hdg_T_shape_Stokes}
\begin{center}
\begin{adjustbox}{max width=\textwidth}
\begin{tabular}{ c c | c c c c c c }
  & & \multicolumn{6}{ c }{\textbf{One-level}} \\
  \textbf{DOF} & \textbf{N} & {\textbf{ORAS}} & {\textbf{SORAS}} & {\textbf{NVTF-MRAS}} & {\textbf{NVTF-SMRAS}} & {\textbf{TVNF-MRAS}} & {\textbf{TVNF-SMRAS}} \\
  \hline
		\textbf{38 803} & \textbf{4} & 22 & 45 & 36 & 49 & 22 & 51 \\
		\textbf{154 606} & \textbf{16} & 111 & 98 & 83 & 172 & 83 & 182 \\
		\textbf{311 369} & \textbf{32} & 265 & 144 & 133 & 262 & 130 & 266 \\
		\textbf{616 772} & \textbf{64} & 568 & 238 & 212 & 410 & 195 & 412 \\
		\textbf{1 246 136} & \textbf{128} & $>$1000 & 494 & 333 & 665 & 313 & 602 \\
		\textbf{2 451 365} & \textbf{256} & $>$1000 & 712 & 464 & $>$1000 & 477 & 889 \\ 
  & & \multicolumn{6}{ c }{\textbf{Two-level (2 eigenvectors)}} \\
  \textbf{DOF} & \textbf{N} & {\textbf{ORAS}} & {\textbf{SORAS}} & {\textbf{NVTF-MRAS}} & {\textbf{NVTF-SMRAS}} & {\textbf{TVNF-MRAS}} & {\textbf{TVNF-SMRAS}} \\
  \hline
		\textbf{38 803} & \textbf{4} & 16 & 35 & 31 & 37 & 21 & 38 \\
		\textbf{154 606} & \textbf{16} & 113 & 69 & 73 & 75 & 38 & 75 \\
		\textbf{311 369} & \textbf{32} & 254 & 99 & 103 & 176 & 93 & 162 \\
		\textbf{616 772} & \textbf{64} & 510 & 153 & 171 & 273 & 121 & 140 \\
		\textbf{1 246 136} & \textbf{128} & $>$1000 & 221 & 242 & 252 & 155 & 138 \\
		\textbf{2 451 365} & \textbf{256} & $>$1000 & 286 & 343 & 515 & 189 & 231 \\ 
  & & \multicolumn{6}{ c }{\textbf{Two-level (5 eigenvectors)}} \\
  \textbf{DOF} & \textbf{N} & {\textbf{ORAS}} & {\textbf{SORAS}} & {\textbf{NVTF-MRAS}} & {\textbf{NVTF-SMRAS}} & {\textbf{TVNF-MRAS}} & {\textbf{TVNF-SMRAS}} \\
  \hline
		\textbf{38 803} & \textbf{4} & 14 & 30 & 27 & 27 & 28 & 30 \\
		\textbf{154 606} & \textbf{16} & 155 & 54 & 54 & 45 & 25 & 44 \\
		\textbf{311 369} & \textbf{32} & 159 & 55 & 72 & 59 & 29 & 52 \\
		\textbf{616 772} & \textbf{64} & 426 & 88 & 106 & 83 & 37 & 76 \\
		\textbf{1 246 136} & \textbf{128} & 955 & 113 & 115 & 99 & 43 & 72 \\
		\textbf{2 451 365} & \textbf{256} & $>$1000 & 182 & 138 & 101 & 54 & 73 \\ 
%  & & \multicolumn{6}{ c }{\textbf{Two-level (7 eigenvectors)}} \\
%  \textbf{DOF} & \textbf{N} & {\textbf{ORAS}} & {\textbf{SORAS}} & {\textbf{NVTF-MRAS}} & {\textbf{NVTF-SMRAS}} & {\textbf{TVNF-MRAS}} & {\textbf{TVNF-SMRAS}} \\
%  \hline
%		\textbf{38 803} & \textbf{4} & 13 & 26 & 25 & 24 & 14 & 28 \\
%		\textbf{154 606} & \textbf{16} & 85 & 46 & 46 & 38 & 24 & 38 \\
%		\textbf{311 369} & \textbf{32} & 168 & 51 & 35 & 44 & 29 & 45 \\
%		\textbf{616 772} & \textbf{64} & 382 & 82 & 75 & 62 & 40 & 59 \\
%		\textbf{1 246 136} & \textbf{128} & 896 & 129 & 86 & 74 & 40 & 64 \\
%		\textbf{2 451 365} & \textbf{256} & $>$1000 & 158 & 47 & 71 & 53 & 61 \\ \\
 \end{tabular}
\end{adjustbox}
\end{center}
\end{table}

In this case scalable results can be only observed for the preconditioners associated with the non standard interface conditions (MRAS and SMRAS), and when using a coarse space which is sufficiently big (see Table~\ref{tab:hdg_T_shape_Stokes}). In the case of ORAS or SORAS, one possibility is to choose a different parameter $\alpha$, but the proof of this, and even the question of whether this would have a positive impact, are open problems. 
}
\end{testcase}
%  \newpage

%%%%%%%%%%%%%%%%%%%%%%%%%%%%%%%%%%%%%%%%%%%%%%%%%%%%%%%%%%%%%%%%%%%%%%%%%%%%%%%%%%%%%%%%%%%%%%%%%%%%%%%%%%%%
% 3D numerics
%%%%%%%%%%%%%%%%%%%%%%%%%%%%%%%%%%%%%%%%%%%%%%%%%%%%%%%%%%%%%%%%%%%%%%%%%%%%%%%%%%%%%%%%%%%%%%%%%%%%%%%%%%%%

\section{Numerical results for three dimensional problems}
\label{sec:3D_numerics}

In this section we again assess the performance of the preconditioners as in Section~\ref{sec:2D_numerics}, but this time in case of three dimensional problems. We consider the partial differential equation model for nearly
incompressible elasticity and Stokes flow as three dimensional problems of similar mixed formulation. Each of
these problems is discretised by using the Taylor-Hood methods from Section~\ref{sec:TH}. In addition, we use the same tools as in Section~\ref{sec:2D_numerics}. For both test cases  we use zero as an initial guess.

\subsection{Taylor-Hood discretisation}
\label{sec:3D_coarse_TH}

In this section we consider the Taylor-Hood discretisation from Section~\ref{sec:TH} with $k = 2$ for nearly incompressible elasticity and Stokes equations. 

\subsubsection{Nearly incompressible elasticity}
In three dimensional space, ORAS and SORAS preconditioners also require an optimized parameter. We follow~\cite{haferssas:hal-01278347} and use Robin interface conditions~\eqref{eq:coarse_elasticity_Robin_BC} with $\alpha = 10$.

\begin{testcase}[The homogeneous beam problem]
\label{exp:homo_beam}
{\rm We consider a homogeneous beam with the physical parameters $E = 10^8$ and $\nu = 0.4999$. The computational domain is the rectangle $\Omega = (0,5) \times (0,1) \times (0,1)$. The beam is clamped on one side, hence we consider the following problem
\begin{equation}
\label{eq:Homogenuous_elasticity}
\left\{
\begin{array}{rclclr}
 -2 \mu \nabla \cdot \boldsymbol{\varepsilon}(\boldsymbol{u}) & + & \nabla p & = & (0,0,-1)^T& \mbox{in } \Omega \\
 & - &  \nabla \cdot \boldsymbol{u} & = & \frac{1}{\lambda} p & \mbox{in } \Omega \\
  	& & \boldsymbol{u}(x,y)  & = & (0,0,0)^T & \mbox{on } \partial \Omega \cap \{x=0\} \\
	& &  \boldsymbol{\sigma}^{sym}_{\boldsymbol{n}}(x,y)  & = & (0,0,0)^T & \mbox{on } \partial \Omega \setminus \{x=0\}
 \end{array}
\right. .
\end{equation}
\begin{table}[!ht]
\caption{Comparison of preconditioners for Taylor-Hood discretisation ($\boldsymbol{TH^2_h}, R^{1}_h$)  - the homogeneous beam.}
\label{tab:TH_homo_beam}
\begin{center}
\begin{adjustbox}{max width=\textwidth}
\begin{tabular}{ c c | c c c c c c }
   & & \multicolumn{6}{ c }{\textbf{One-level}} \\
  \textbf{DOF} & \textbf{N} & {\textbf{ORAS}} & {\textbf{SORAS}} & {\textbf{NDTNS-MRAS}} & {\textbf{NDTNS-SMRAS}} & {\textbf{TDNNS-MRAS}} & {\textbf{TDNNS-SMRAS}} \\
  \hline
		\textbf{32 446} & \textbf{8}  & 21 & 45 & 29 & 36 & 27 & 37 \\
		\textbf{73 548} & \textbf{16}  & 31 & 70 & 38 & 64 & 26 & 67 \\
		\textbf{139 794} & \textbf{32}  & 43 & 99 & 74 & 94 & 66 & 91 \\
		\textbf{299 433} & \textbf{64}  & 55 & 143 & 161 & 140 & 149 & 139 \\
		\textbf{549 396} & \textbf{128} & 78 & 192 & 314 & 192 & 229 & 199 \\
% 		\textbf{1 058 184} & \textbf{256} &  &  &  &  &  &  \\ \\
  & & \multicolumn{6}{ c }{\textbf{Two-level (6 eigenvectors)}} \\
  \textbf{DOF} & \textbf{N} & {\textbf{ORAS}} & {\textbf{SORAS}} & {\textbf{NDTNS-MRAS}} & {\textbf{NDTNS-SMRAS}} & {\textbf{TDNNS-MRAS}} & {\textbf{TDNNS-SMRAS}} \\
  \hline 
		\textbf{32 446} & \textbf{8}  & 10 & 17 & 13 & 18 & 12 & 17  \\
		\textbf{73 548} & \textbf{16}  & 11 & 22 & 16 & 25 & 16 & 22 \\
		\textbf{139 794} & \textbf{32}  & 13 & 26 & 25 & 28 & 17 & 26 \\
		\textbf{299 433} & \textbf{64}  & 15 & 27 & 19 & 27 & 24 & 28 \\
		\textbf{549 396} & \textbf{128} & 17 & 28 & 20 & 25 & 21 & 26 \\
% 		\textbf{1 058 184} & \textbf{256} &  &  &  &  &  &  \\ \\
& & \multicolumn{6}{ c }{\textbf{Two-level (8 eigenvectors)}} \\
  \textbf{DOF} & \textbf{N} & {\textbf{ORAS}} & {\textbf{SORAS}} & {\textbf{NDTNS-MRAS}} & {\textbf{NDTNS-SMRAS}} & {\textbf{TDNNS-MRAS}} & {\textbf{TDNNS-SMRAS}} \\
  \hline
		\textbf{32 446} & \textbf{8}  & 9 & 16 & 12 & 17 & 12 & 16 \\
		\textbf{73 548} & \textbf{16}  & 10 & 19 & 15 & 24 & 14 & 20 \\
		\textbf{139 794} & \textbf{32}  & 11 & 21 & 17 & 23 & 17 & 21 \\
		\textbf{299 433} & \textbf{64}  & 14 & 24 & 17 & 24 & 21 & 23 \\
		\textbf{549 396} & \textbf{128} & 16 & 27 & 18 & 23 & 20 & 22 \\
% 		\textbf{1 058 184} & \textbf{256} &  &  &  &  &  &  \\ \\
 \end{tabular}
\end{adjustbox}
\end{center}
\end{table}
  
The results of Table~\ref{tab:TH_homo_beam} show a clear improvement in the scalability of the two-level preconditioners over the one-level ones. In fact, using only zero energy modes, the number of iterations is virtually unaffected by the number of subdomains. All two-level preconditioners show a comparable performance. For this case, increasing the dimension of the coarse space beyond $6 \times N$ eigenvectors does not seem to improve the results dramatically.
}
\end{testcase}

\subsubsection{Stokes equation}
We now turn to the Stokes discrete problem given in~\ref{sec:TH}. Once again in case of ORAS and SORAS we choose $\alpha = 10$ as in~\cite{haferssas:hal-01278347} for the Robin interface conditions~\eqref{eq:coarse_elasticity_Robin_BC}.

\begin{testcase}[The driven cavity problem]
{\rm The test case is the three-dimensional version of the driven cavity problem. We consider the following problem on the unit cube $\Omega = (0,1)^3$
\begin{equation}
\label{eq:3D_Cavity}
\left\{
\begin{array}{rclclr}
 - \Delta \boldsymbol{u} & + & \nabla p & = & \boldsymbol{f} & \mbox{in } \Omega \\
 & -&  \nabla \cdot \boldsymbol{u} & = & 0 & \mbox{in } \Omega \\
  	& & \boldsymbol{u}(x,y) & = &(1,0,0)^T & \mbox{on } \partial \Omega \cap \{y=1\} \\
	& &  \boldsymbol{u}(x,y)& = & (0,0,0)^T & \mbox{on } \partial \Omega \setminus \{y=1\}
 \end{array}
\right. .
\end{equation}

\begin{table}[!ht]
\caption{Comparison of preconditioners for Taylor-Hood discretisation ($\boldsymbol{TH^2_h}, R^{1}_h$) - the driven cavity problem.}
\label{tab:3D_TH_Cavity_Stokes}
\begin{center}
\begin{adjustbox}{max width=\textwidth}
\begin{tabular}{ c c | c c c c c c }
%   & \multicolumn{4}{ c }{\textbf{Uniform decomposition (236 $\times$ 236)}} \\ 
%   \hline
  & & \multicolumn{6}{ c }{\textbf{One-level}} \\
  \textbf{DOF} & \textbf{N} & {\textbf{ORAS}} & {\textbf{SORAS}} & {\textbf{NVTF-MRAS}} & {\textbf{NVTF-SMRAS}} & {\textbf{TVNF-MRAS}} & {\textbf{TVNF-SMRAS}} \\
  \hline
		\textbf{38 229} & \textbf{8} & 12 & 24 & 12 & 22 & 11 & 23 \\
		\textbf{76 542} & \textbf{16} & 18 & 34 & 18 & 31 & 15 & 31 \\
		\textbf{158 818} & \textbf{32} & 23 & 45 & 20 & 45 & 19 & 45 \\
		\textbf{325 293} & \textbf{64} & 28 & 60 & 36 & 64 & 25 & 60 \\
		\textbf{643 137} & \textbf{128} & 37 & 79 & 64 & 91 & 33 & 88 \\
% 		\textbf{1 217 704} & \textbf{256} &  &  &  &  &  &  \\ \\
  & & \multicolumn{6}{ c }{\textbf{Two-level (3 eigenvectors)}} \\
 \textbf{DOF} & \textbf{N} & {\textbf{ORAS}} & {\textbf{SORAS}} & {\textbf{NVTF-MRAS}} & {\textbf{NVTF-SMRAS}} & {\textbf{TVNF-MRAS}} & {\textbf{TVNF-SMRAS}} \\
  \hline
		\textbf{38 229} & \textbf{8} & 10 & 17 & 10 & 18 & 11 & 18 \\
		\textbf{76 542} & \textbf{16} & 11 & 20 & 11 & 19 & 14 & 19 \\
		\textbf{158 818} & \textbf{32} & 13 & 24 & 13 & 24 & 16 & 23 \\
		\textbf{325 293} & \textbf{64} & 15 & 27 & 15 & 27 & 19 & 26 \\
		\textbf{643 137} & \textbf{128} & 18 & 31 & 17 & 32 & 22 & 31 \\
% 		\textbf{1 217 704} & \textbf{256} &  &  &  &  &  &  \\ \\
  & & \multicolumn{6}{ c }{\textbf{Two-level (7 eigenvectors)}} \\
  \textbf{DOF} & \textbf{N} & {\textbf{ORAS}} & {\textbf{SORAS}} & {\textbf{NVTF-MRAS}} & {\textbf{NVTF-SMRAS}} & {\textbf{TVNF-MRAS}} & {\textbf{TVNF-SMRAS}} \\
  \hline
		\textbf{38 229} & \textbf{8} & 9 & 16 & 9 & 16 & 12 & 17 \\
		\textbf{76 542} & \textbf{16} & 10 & 17 & 10 & 18 & 15 & 17 \\
		\textbf{158 818} & \textbf{32} & 11 & 19 & 11 & 20 & 17 & 20 \\
		\textbf{325 293} & \textbf{64} & 13 & 19 & 13 & 21 & 20 & 20 \\
		\textbf{643 137} & \textbf{128} & 15 & 21 & 16 & 22 & 22 & 22 \\
% 		\textbf{1 217 704} & \textbf{256} &  &  &  &  &  &  \\ \\
 \end{tabular}
\end{adjustbox}
\end{center}
\end{table}

The conclusions remain the same as for the homogeneous beam example for the nearly incompressible elasticity equation discretised by Taylor-Hood method ($\boldsymbol{TH^3_h}, R^{2}_h$) since Tables %~\ref{tab:TH_Cavity_Stokes},
 \ref{tab:3D_TH_Cavity_Stokes} and~\ref{tab:TH_homo_beam} %, \ref{tab:lowest_TH_L_shape_elasticity}
 show similar results.
}
\end{testcase}
 \newpage

\section{Conclusion}
\label{sec:conclusions}
 
We tested numerically two-level preconditioners with spectral coarse spaces for nearly incompressible elasticity and Stokes equations. We considered two finite element methods, namely, Taylor-Hood (Section~\ref{sec:TH}) and the hdG (Section~\ref{sec:hdg}) discretisations.

In the case of the homogeneous nearly incompressible elasticity the two-level methods coupled with SORAS preconditioner defined in~\cite{haferssas:hal-01278347} and SMRAS preconditioner defined by~\eqref{eq:SMRAS} allowed us to achieve good scalability results for both discretisations. Furthermore, for these symmetric preconditioners coarse spaces containing only zero energy modes seem to be enough for two and three dimensional problems. For the heterogeneous problem we also achieved scalability for two-level SORAS and SMRAS preconditioners, but, as expected, only in the case when the size of the coarse space is sufficiently big.

The improvement of the convergence in the case of the Stokes flow is visible only when the coarse space contains more eigenvectors than only constants. For the Taylor-Hood discretisation, taking sufficient big coarse space we were able to achieve good scalability for all preconditioners. It is remarkable that these good results occur even when using the hdG discretisation, despite the fact the optimized parameter to be used in SORAS and ORAS is not available.
% In the case of the hdG discretisation, coarse spaces applied to MRAS and SMRAS preconditioners gave satisfying results. Although, an optimal parameter for the Stokes and elasticity problems descretised by the hdG method is unavailable. 

We can conclude that the two-level preconditioners associated with non standard interface conditions are at least as good as the two-level ones in conjunction with Robin interface conditions using optimised parameters. It shows an important advantage of newly introduced preconditioners as they are parameter-free. 

Numerical tests have shown that the coarse spaces bring an important improvement in the convergence, but the size of the coarse space depends on the problem. Building as small as possible coarse spaces is important from computational point of view. Thus, it is necessary to investigate what could be an optimal criterion for choosing the eigenvectors for a coarse space.

%A theoretical study of multilevel methods requires proving the estimates that control the spectrum of the preconditioned operator. This is done for the systems where
As we mentioned before, the theoretical foundation of the two-level preconditioners has not been extended to saddle point problems. Hence, future research will be devoted to this topic.
% re is missing the theory for the two-level methods for saddle point problems, where the global matrix is non-symmetric. %However, it is not the case of saddle point problems such as nearly incompressible elasticity and Stokes. 
% Since the numerical results showed the good behaviour for such systems as Stokes and linear elasticity, the development of the theory %prove of these estimates 
% in the case of the non-symmetric problems is the interest of future studies.

\section*{Acknowledgements}

This research was supported supported by the Centre for Numerical Analysis and Intelligent Software (NAIS). We thank Fr\'{e}d\'{e}ric Nataf and Pierre-Henri Tournier for many helpful discussions and insightful comments, and Ryadh Haferssas and Fr\'{e}d\'{e}ric Hecht for their assistance with the FreeFem++ codes.

%%%%%%%%%%%%%%%%%%%%%%%%%%%%%%%%%%%%%%%%%%%%%%%%%%%%%%%%%%%%%%%%%%%%%%%%%%%%%%%%%%%%%%%%%%%%%%%%%%%%%%%%%%%%
% BILBIOGRAPHY
%%%%%%%%%%%%%%%%%%%%%%%%%%%%%%%%%%%%%%%%%%%%%%%%%%%%%%%%%%%%%%%%%%%%%%%%%%%%%%%%%%%%%%%%%%%%%%%%%%%%%%%%%%%%

\thispagestyle{empty}
  \bibliographystyle{alpha}
\bibliography{DD}

\newcommand{\etalchar}[1]{$^{#1}$}
\begin{thebibliography}{AdDBM{\etalchar{+}}14}

\bibitem[AdDBM{\etalchar{+}}14]{de2014simple}
B.~Ayuso~de Dios, F.~Brezzi, L.~D. Marini, J.~Xu, and L.~Zikatanov.
\newblock A simple preconditioner for a discontinuous {G}alerkin method for the
  {S}tokes problem.
\newblock {\em Journal of Scientific Computing}, 58(3):517--547, 2014.

\bibitem[BBD{\etalchar{+}}16]{barrenechea2016hybrid}
G.~R Barrenechea, M.~Bosy, V.~Dolean, F.~Nataf, and P.-H. Tournier.
\newblock Hybrid discontinuous {G}alerkin discretisation and preconditioning of
  the {S}tokes problem with non standard boundary conditions.
\newblock preprint, \url{https://arxiv.org/abs/1610.09207}, October 2016.

\bibitem[BBF13]{MR3097958}
D.~Boffi, F.~Brezzi, and M.~Fortin.
\newblock {\em Mixed finite element methods and applications}, volume~44 of
  {\em Springer Series in Computational Mathematics}.
\newblock Springer, Heidelberg, 2013.

\bibitem[BHMV99]{brezina1999iterative}
M.~Brezina, C.~I Heberton, J.~Mandel, and P.~Vanek.
\newblock An iterative method with convergence rate chosen a priori, ucd/ccm
  report 140.
\newblock Technical report, Center for Computational Mathematics, University of
  Colorado at Denver, 1999.

\bibitem[CDKN14]{MR3209915}
L.~Conen, V.~Dolean, R.~Krause, and F.~Nataf.
\newblock A coarse space for heterogeneous {H}elmholtz problems based on the
  {D}irichlet-to-{N}eumann operator.
\newblock {\em J. Comput. Appl. Math.}, 271:83--99, 2014.

\bibitem[CG09]{MR2485446}
B.~Cockburn and J.~Gopalakrishnan.
\newblock The derivation of hybridizable discontinuous {G}alerkin methods for
  {S}tokes flow.
\newblock {\em SIAM J. Numer. Anal.}, 47(2):1092--1125, 2009.

\bibitem[CGL09]{MR2485455}
B.~Cockburn, J.~Gopalakrishnan, and R.~Lazarov.
\newblock Unified hybridization of discontinuous {G}alerkin, mixed, and
  continuous {G}alerkin methods for second order elliptic problems.
\newblock {\em SIAM J. Numer. Anal.}, 47(2):1319--1365, 2009.

\bibitem[CS15]{MR3407237}
C.~Carstensen and M.~Schedensack.
\newblock Medius analysis and comparison results for first-order finite element
  methods in linear elasticity.
\newblock {\em IMA J. Numer. Anal.}, 35(4):1591--1621, 2015.

\bibitem[DJN15]{MR3450068}
V.~Dolean, P.~Jolivet, and F.~Nataf.
\newblock {\em An introduction to domain decomposition methods}.
\newblock Society for Industrial and Applied Mathematics (SIAM), Philadelphia,
  PA, 2015.
\newblock Algorithms, theory, and parallel implementation.

\bibitem[DNSS12]{MR3033238}
V.~Dolean, F.~Nataf, R.~Scheichl, and N.~Spillane.
\newblock Analysis of a two-level {S}chwarz method with coarse spaces based on
  local {D}irichlet-to-{N}eumann maps.
\newblock {\em Comput. Methods Appl. Math.}, 12(4):391--414, 2012.

\bibitem[EGLW12]{MR2916377}
Y.~Efendiev, J.~Galvis, R.~Lazarov, and J.~Willems.
\newblock Robust domain decomposition preconditioners for abstract symmetric
  positive definite bilinear forms.
\newblock {\em ESAIM Math. Model. Numer. Anal.}, 46(5):1175--1199, 2012.

\bibitem[GE10a]{MR2718268}
J.~Galvis and Y.~Efendiev.
\newblock Domain decomposition preconditioners for multiscale flows in
  high-contrast media.
\newblock {\em Multiscale Model. Simul.}, 8(4):1461--1483, 2010.

\bibitem[GE10b]{MR2728702}
J.~Galvis and Y.~Efendiev.
\newblock Domain decomposition preconditioners for multiscale flows in high
  contrast media: reduced dimension coarse spaces.
\newblock {\em Multiscale Model. Simul.}, 8(5):1621--1644, 2010.

\bibitem[GR86]{MR851383}
V.~Girault and P.~A. Raviart.
\newblock {\em Finite element methods for {N}avier-{S}tokes equations},
  volume~5 of {\em Springer Series in Computational Mathematics}.
\newblock Springer-Verlag, Berlin, 1986.
\newblock Theory and algorithms.

\bibitem[Hec12]{MR3043640}
F.~Hecht.
\newblock New development in {F}ree{F}em++.
\newblock {\em J. Numer. Math.}, 20(3-4):251--265, 2012.

\bibitem[HJN15]{haferssas:hal-01278347}
R.~Haferssas, P.~Jolivet, and F.~Nataf.
\newblock {An additive Schwarz method type theory for Lions' algorithm and
  Optimized Schwarz Methods}.
\newblock preprint, \url{https://hal.archives-ouvertes.fr/hal-01278347},
  December 2015.

\bibitem[KK98]{karypis1998software}
G.~Karypis and V.~Kumar.
\newblock A software package for partitioning unstructured graphs, partitioning
  meshes, and computing fill-reducing orderings of sparse matrices.
\newblock Technical report, University of Minnesota, Department of Computer
  Science and Engineering, Army HPC Research Center, Minneapolis, MN, 1998.

\bibitem[LNS15]{MR3432853}
S\'ebastien Loisel, Hieu Nguyen, and Robert Scheichl.
\newblock Optimized {S}chwarz and 2-{L}agrange multiplier methods for
  multiscale elliptic {PDE}s.
\newblock {\em SIAM J. Sci. Comput.}, 37(6):A2896--A2923, 2015.

\bibitem[LSY98]{MR1621681}
R.~B. Lehoucq, D.~C. Sorensen, and C.~Yang.
\newblock {\em A{RPACK} users' guide}, volume~6 of {\em Software, Environments,
  and Tools}.
\newblock Society for Industrial and Applied Mathematics (SIAM), Philadelphia,
  PA, 1998.
\newblock Solution of large-scale eigenvalue problems with implicitly restarted
  Arnoldi methods.

\bibitem[Nic87]{MR881370}
R.~A. Nicolaides.
\newblock Deflation of conjugate gradients with applications to boundary value
  problems.
\newblock {\em SIAM J. Numer. Anal.}, 24(2):355--365, 1987.

\bibitem[NXD10]{MR2738920}
F.~Nataf, H.~Xiang, and V.~Dolean.
\newblock A two level domain decomposition preconditioner based on local
  {D}irichlet-to-{N}eumann maps.
\newblock {\em C. R. Math. Acad. Sci. Paris}, 348(21-22):1163--1167, 2010.

\bibitem[PS11]{MR2826472}
A.~Pechstein and J.~Sch{\"o}berl.
\newblock Tangential-displacement and normal-normal-stress continuous mixed
  finite elements for elasticity.
\newblock {\em Math. Models Methods Appl. Sci.}, 21(8):1761--1782, 2011.

\bibitem[SDH{\etalchar{+}}14]{MR3175183}
N.~Spillane, V.~Dolean, P.~Hauret, F.~Nataf, C.~Pechstein, and R.~Scheichl.
\newblock Abstract robust coarse spaces for systems of {PDE}s via generalized
  eigenproblems in the overlaps.
\newblock {\em Numer. Math.}, 126(4):741--770, 2014.

\bibitem[SS86]{MR848568}
Y.~Saad and M.~H. Schultz.
\newblock G{MRES}: a {G}eneralized {M}inimal {R}esidual algorithm for solving
  nonsymmetric linear systems.
\newblock {\em SIAM J. Sci. Statist. Comput.}, 7(3):856--869, 1986.

\bibitem[TW05]{Toselli:2005:DDM}
A.~Toselli and O.~Widlund.
\newblock {\em Domain Decomposition Methods - Algorithms and Theory}, volume~34
  of {\em Springer Series in Computational Mathematics}.
\newblock Springer, 2005.

\end{thebibliography}

\end{document}